\documentclass[11pt]{article}

\usepackage[english]{babel}
\usepackage[T1]{fontenc}
\usepackage[utf8]{inputenc}
\usepackage{amsbsy}
\usepackage{amsmath}
\usepackage{amssymb}
\usepackage{amsfonts}
\usepackage{amsthm}
\usepackage{bm}
\usepackage{graphicx}
\usepackage[active]{srcltx}
\usepackage{upgreek}
\usepackage{xcolor}
\usepackage[all]{xy}
\usepackage{dsfont}

\usepackage{tikz}

\newcommand{\N}{\mathbb{N}}

\newcommand{\R}{\mathbb{R}}

\newcommand{\ga}{\mathfrak{a}}
\newcommand{\gb}{\mathfrak{b}}
\newcommand{\gc}{\mathfrak{c}}
\newcommand{\gd}{\mathfrak{d}}
\renewcommand{\ge}{\mathfrak{e}}

\newcommand{\gp}{\mathfrak{p}}
\newcommand{\gq}{\mathfrak{q}}

\newcommand{\gv}{\mathfrak{v}}

\newcommand{\energyset}{\mathcal{X}(\R)}
\newcommand{\energysethydro}{\mathcal{X}_{hy}(\R)}

\newcommand{\energysethydrok}{\mathcal{X}_{hy}^k(\R)}
\newcommand{\Nenergyset}{\mathcal{NX}(\R)}
\newcommand{\Nenergysethydro}{\mathcal{NX}_{hy}(\R)}
\newcommand{\Nenergysethydrok}{\mathcal{NX}_{hy}^k(\R)}

\newcommand{\Rca}{R_{\gc,\ga}}
\newcommand{\Operpca}{\mathcal{U}^\perp_{\gc,\ga}(L)}

\newcommand{\grandOde}[1]{\mathcal{O}\left( #1\right)}

\newcommand{\ii}{\infty}

\newcommand{\spess}{\sigma_{ess}}

\newcommand{\adm}{\mathrm{Adm}}
\newcommand{\pos}{\mathrm{Pos}}

\newcommand{\normX}[1]{\left\Vert #1\right\Vert_{\mathcal{X}}}
\newcommand{\normLii}[1]{\left\Vert #1\right\Vert_{L^\ii}}
\newcommand{\normLdeux}[1]{\left\Vert #1\right\Vert_{L^2}}
\newcommand{\normLun}[1]{\left\Vert #1\right\Vert_{L^1}}
\newcommand{\normLp}[1]{\left\Vert #1\right\Vert_{L^p}}

\newcommand{\psLdeux}[2]{\left\langle #1,#2 \right\rangle_{L^2}}
\newcommand{\psLdeuxLdeux}[2]{\left\langle #1,#2 \right\rangle_{L^2\times L^2}}
\newcommand{\vvvert}[1]{\big\vert\kern-0.25ex\big\vert\kern-0.25ex\big\vert #1\big\vert\kern-0.25ex\big\vert\kern-0.25ex\big\vert}

\newcommand{\spanned}[1]{\mathrm{Span}(#1)}
\newcommand{\normLdeuxLdeux}[1]{\left\Vert #1\right\Vert_{L^2\times L^2}}

\newcommand{\normR}[1]{\left| #1\right|}

\newtheorem{thm}{Theorem}[section]

\newtheorem{claim}[thm]{Claim}
\newtheorem{cor}[thm]{Corollary}
\newtheorem{lem}[thm]{Lemma}
\newtheorem{prop}[thm]{Proposition}

\newtheorem*{thm*}{Theorem}

\newtheorem{rem}[thm]{Remark}

\theoremstyle{remark}

\title{Orbital stability of a chain of dark solitons for general nonintegrable Schrödinger equations with non-zero condition at infinity.}

\begin{document}

\date{}
\author{
\renewcommand{\thefootnote}{\arabic{footnote}}
Jordan Berthoumieu\footnotemark[1]}
\footnotetext[1]{CY Cergy Paris Universit\'e, Laboratoire Analyse, G\'eom\'etrie, Mod\'elisation, F-95302 Cergy-Pontoise, France. Website: https://berthoumieujordan.wordpress.com \ E-mail: {\tt jordan.berthoumieu@cyu.fr}}
\maketitle

\begin{abstract}

In this article, we focus on the stability of dark solitons for a general one-dimensional nonlinear Schrödinger equation. More precisely, we prove the orbital stability of a chain of travelling waves whose speeds are well ordered, taken close to the speed of sound $c_s$ and such that the solitons are initially localized far away from each other. The proof relies on the arguments developed by F. Béthuel, P. Gravejat and D. Smets in~\cite{BetGrSm1} and first introduced in~\cite{MarMeTs1} by Y. Martel, F. Merle and T.-P. Tsai.
\end{abstract}

\section{Introduction}\label{section: intro}

We are interested in the defocusing nonlinear Schrödinger equation

\begin{equation}\label{NLS}\tag{$NLS$}
    i\partial_t \Psi +\partial_x^2 \Psi +\Psi f(|\Psi|^2)=0\quad\text{on }\R\times\R.
\end{equation}

This equation appears as a relevant model in condensed matter physics. In particular, it arises in the context of the Bose-Einstein condensation or of superfluidity (see~\cite{AbiHuMeNoPhTu,Coste1}), but also in nonlinear optics (see~\cite{KivsLut1,KivPeSt1}), when the natural condition at infinity is 

\begin{equation}\label{nonvanishing condition à l'infini}
    |\Psi(t,x)|\underset{|x|\rightarrow +\ii}{\longrightarrow}1.
\end{equation}

In the latter context, this condition expresses the presence of a nonzero background. It differs from the case of null condition at infinity, in the sense that the resulting dispersion relation is different.\\
In \eqref{NLS}, the nonlinearity $f$ can be taken equal to $f(\rho)=1-\rho$ and we then obtain the so-called Gross-Pitaevskii equation, but we can also take many other functions $f$ satisfying the condition $f(1)=0$. This provides possible alternative behaviours as enumerated by D. Chiron in~\cite{Chiron7}. In the sequel, we restrict our attention to the defocusing case, meaning that we assume that

\begin{equation}\label{f'(1) <0}
    f'(1)<0.
\end{equation}

If $\Psi$ does not vanish, we can lift it as $\Psi=\rho e^{i\varphi}$ where $\rho$ and $\varphi$ are as smooth as $f$ is. Setting the new variables $\eta=1-\rho^2$ and $v=-\partial_x \varphi$, we obtain the hydrodynamical form of the equation 

\begin{equation}\label{NLShydro}\tag{$NLS_{hy}$}
    \left\{
\begin{array}{l}
    \partial_t\eta =-2\partial_x \big(v(1-\eta)\big), \\
    \partial_t v =-\partial_x \bigg( f(1-\eta)-v^2-\dfrac{\partial_x^2 \eta}{2(1-\eta)}-\dfrac{(\partial_x\eta)^2}{4(1-\eta)^2} \bigg). \\
\end{array}
\right.
\end{equation}

The linearized system around the trivial solution $(\rho,v)=(1,0)$ reduces, in the long wave approximation, to the free wave equation with the sound speed

\begin{equation}\label{relation entre vitesse c_s et f'(1)}
c_s = \sqrt{-2 f'(1)}.
\end{equation}

The plane wave solutions of this linearized system satisfy indeed the dispersion relation
\begin{equation*}
    \omega(\xi)=\pm\sqrt{\xi^4+c_s^2 \xi^2}.
\end{equation*}

For $k\geq 0$, we shall write $\energysethydrok :=H^{k+1}(\R)\times H^k(\R)$ and endow this space with the associated euclidean norm 
\begin{equation*}
    \Vert (\eta,v)\Vert_{\mathcal{X}^k}^2=\Vert \eta\Vert_{H^{k+1}}^2 + \Vert v\Vert_{H^{k}}^2.
\end{equation*}

Let us also introduce the non-vanishing associated subset

\begin{equation*}
    \Nenergysethydrok:=\big\{(\eta,v)\in \energysethydrok\big| \max_\R\eta <1\big\}.
\end{equation*}

We will label $\energysethydro:=\energysethydrok$ and $\Nenergysethydro:=\Nenergysethydrok$ when $k=0$. This functional setting is related to several quantities, which are at least formally, conserved along the flow. These are the energy

\begin{equation}\label{expression de l'énergie}
    E(\eta,v):=\int_\R e(\eta,v):=\dfrac{1}{8}\int_\R\dfrac{(\partial_x\eta)^2}{1-\eta}+\dfrac{1}{2}\int_\R (1-\eta)v^2+\dfrac{1}{2}\int_\R F(1-\eta),
\end{equation}

with

\begin{equation*}\label{definition de F}
    F(r):=\int_r^1 f(\rho)d\rho,
\end{equation*}

and (the renormalized) momentum 

\begin{equation}\label{expression du moment}
    p(\eta,v):=\int_\R\pi(\eta,v):=\dfrac{1}{2}\int_\R \eta v .
\end{equation}

As it is mentionned in~\cite{Bert1}, the assumptions that we will make on the nonlinearity $f$ give to $\Nenergysethydro$ the status of the energy set of~\eqref{NLShydro}. Note also that, in view of the previous expression, $p$ is smooth on $\energysethydro$. Moreover, if $f\in\mathcal{C}^k(\R_+)$, $E$ is $\mathcal{C}^{k+1}$ on $\Nenergysethydro$. Concerning the existence of solutions to~\eqref{NLShydro}, we recall the following local wellposedness theorem.

\begin{thm}[Gallo~\cite{Gallo3}]\label{théoreme: LWP of NLShydro}
    We assume that $f\in\mathcal{C}^3(\R_+)$ is such that for all $\rho\in\R$,
\begin{equation}\label{hypothèse de croissance sur F minorant intermediaire}\tag{H1}
    \dfrac{c_s^2}{4} (1-\rho)^2 \leq F(\rho).
\end{equation}
Let $k\in\{0,1,2\}$ and let $(\eta_0,v_0)\in \Nenergysethydrok$. There exist $T_{\max} >0$ and a unique solution $(\eta,v)\in\mathcal{C}^0\big([0,T_{\max} ),\Nenergysethydrok\big)$ to~\eqref{NLShydro} with initial datum $(\eta_0,v_0)$. The maximal time $T_{\max}$ is continuous with respect to the initial datum and is characterized by
    \begin{equation*}
        \lim_{t\rightarrow T_{\max}^-}\sup_{x\in\R}\eta(t,x) = 1.
    \end{equation*}

Moreover, the flow map is continuous on $\Nenergysethydrok$ and the energy and the momentum are conserved along the flow.
\end{thm}

\subsection{Travelling wave solutions and minimizing property}\label{section: minimizing property of the travelling waves}
In this article, we aim at understanding the dynamics of a chain of dark solitons when their speeds are ordered and taken close to the sound of speed $c_s$. "Dark soliton" is a general term designating special solutions of~\eqref{NLS}. In our framework, they are travelling waves of nonzero speed\footnote{They are also labelled grey solitons by opposition with the black soliton with speed $c=0$.} of the form $\Psi(t,x)=u(x-ct)$. We plug this ansatz in~\eqref{NLS} and we obtain the ordinary differential equation satisfied by the profile $u$, which is
\begin{equation}\label{TWC}\tag{$TW_c$}
    -icu'+u''+uf(|u|^2)=0.
\end{equation}

Under suitable conditions on the nonlinearity $f\in\mathcal{C}^3(\R_+)$ there exist travelling waves of given momentum in a certain range $(0,\gq_*)$ (see~\cite{Bert1}) and they are orbitally stable. We assume that these conditions hold true in the sequel. Namely

\begin{itemize}
    \item For all $\rho\in\R$,
\begin{equation}\label{hypothèse de croissance sur F minorant intermediaire}\tag{H1}
    \dfrac{c_s^2}{4} (1-\rho)^2 \leq F(\rho).
\end{equation}
    \item There exist $M\geq 0$ and $q\in [2, +\ii)$ such that for all $\rho \geq 2$,
\begin{equation}\label{hypothèse de croissance sur F majorant}\tag{H2}
    F(\rho)\leq M|1-\rho|^q.
\end{equation}
    \item \begin{equation}\label{condition suffisante pour la stabilité orbitale sur f''(1)+6f'(1)>0}\tag{H3}
    f''(1)+3f'(1)\neq 0.
\end{equation}
\end{itemize}

This result relies on a variational argument which consists on solving the minimization problem

\begin{equation}\label{définition de la minimization curve Emin}
    E_{\min}(\mathfrak{p}):=\inf\big\lbrace E(v)\big| v\in\mathcal{N}\energyset, p(v)=\mathfrak{p}\big\rbrace ,
\end{equation}

where $E$ and $p$ denote the energy and the momentum associated with the formulas~\eqref{expression de l'énergie} and~\eqref{expression du moment}. The Euler-Lagrange equation 

\begin{equation}\label{nabla E(v)-cnabla p(v)=0}
    \nabla E(v)-c\nabla p(v)=0
\end{equation}
indeed reduces to~\eqref{TWC} where the speed $c$ is the Lagrange multiplier of the problem. The question of uniqueness of the minimizer goes beyond the scope of this article but we shall touch upon a partial result in this way in Section~\ref{section: chain of solitons}. For $\gp\in (0,\gq_*)$, we shall label $\gv_{c(\gp)}$ a dark soliton minimizing~\eqref{définition de la minimization curve Emin} and mention that, due to condition~\eqref{condition suffisante pour la stabilité orbitale sur f''(1)+6f'(1)>0} (see Theorem~5.1 and the remark just after in~\cite{Maris6}), $c(\gp)$ necessarily lies in $(-c_s,c_s)\setminus\{ 0\}$.

Since $c(\gp)\neq 0$, the solutions $\gv_{c(\gp)}$ do not vanish on $\R$ (see Remark~4.6 in~\cite{Bert1}). Plugging the Madelung transform in~\eqref{TWC}, we obtain the hydrodynamical version of~\eqref{TWC} satisfied by the variables $(\eta_{c},v_{c}):=(1-|\gv_c|^2,-\varphi_c')$,
\begin{equation}\label{TWChydro}\tag{$TW_{c,hy}$}
    \left\{
\begin{array}{l}
    \dfrac{c}{2}\eta_c=v_c(1-\eta_c), \\
    c v_c = f(1-\eta_c)-v_c^2-\dfrac{\eta_c''}{2(1-\eta_c)}-\dfrac{(\eta_c')^2}{4(1-\eta_c)^2} . \\
\end{array}
\right.
\end{equation}

For the sake of completeness, we will check in Subsection~\ref{subsection: exponential decay} the bound
\begin{equation}\label{max eta_c <1}
    \max_\R\eta_c <1 ,
\end{equation}
for solitons with speed $c$ close to $c_s$. In the sequel, we will also take into account the invariance by translation of the equation by setting

\begin{equation*}
    Q_c:=(\eta_c,v_c)\quad\text{and}\quad Q_{c,a}:=(\eta_{c,a},v_{c,a}):=\big(\eta_c(.-a),v_c(.-a)\big),
\end{equation*}

for $a\in\R$.

\subsection{Chain of solitons in the transonic regime}\label{section: chain of solitons}
In~\cite{Chiron7}, D. Chiron proved the existence of a continuous branch of travelling waves in the transonic limit, i.e. with speed close to $c_s$. We slightly improve this result in the next theorem.

\begin{thm}\label{théorème: il existe une branche C^1 de solitons proche de c_s}
    There exists a critical speed $c_0 >0$ such that for $c\in (c_0,c_s)$, there exists a non constant smooth solution $Q_c$ of \eqref{TWChydro}, that is unique up to translations.  Furthermore, the mapping $c\in (c_0,c_s)\mapsto Q_c$ belongs to $\mathcal{C}^2\big((c_0,c_s),\mathcal{NX}^2(\R)\big)$ and for $c\in (c_0,c_s)$, we have 
    \begin{equation}\label{théorème: hydrodynamique condition de grillakis sans le signe}
    \dfrac{d}{dc}\big(p(Q_c)\big) < 0.
\end{equation}

Moreover there exists $a_{d},K_{d}>0$ independent of $c\in (c_0,c_s)$ and $x\in\R$ such that,
    \begin{equation}\label{estimée décroissance exponentielle à tout ordre pour eta_c et v_c}
\sum_{\substack{0\leq k_1\leq 3\\ 0\leq k_2\leq 2 \\ 0\leq j\leq 2}}(c_s^2-c^2)^{j-1}\Big(|\partial_c^j\partial_x^{k_1} \eta_c(x)|+c^{1+2j+2k_2}|\partial_c^j\partial_x^{k_2} v_c(x)|\Big)\leq K_{d}e^{-a_{d}\sqrt{c_s^2-c^2} |x|}.
\end{equation}
\end{thm}

\begin{rem}
    We can show that the branch is actually $\mathcal{C}^3\left((c_0,c_s),\big(H^2(\R)\times L^2(\R)\big)\cap\big(\mathcal{C}_{loc}^5(\R)\big)^2\right)$. Nevertheless, the choice of the specific exponents in Theorem~\ref{théorème: il existe une branche C^1 de solitons proche de c_s} is sharp for the study in this article.
\end{rem}

\begin{rem}
    The restriction to positive speeds is not an arbitrary choice. Indeed note that if $Q_c$ is a hydrodynamical travelling wave of speed $c$, then $(\eta_c,-v_c)$ is a hydrodynamical travelling wave of speed $-c$ and our results extend naturally to these travelling waves.
\end{rem}

In~\cite{GriShSt1}, M. Grillakis, J. Shatah and W.A. Strauss gave some sharp conditions for orbital stability (and instability) for a general class of Hamiltonian systems. D. Chiron showed in~\cite{Chiron8} that this general result could be specified to travelling waves for nonlinear Schrödinger type equations. The aforementioned condition for orbital stability is precisely~\eqref{théorème: hydrodynamique condition de grillakis sans le signe}.

Combining Theorem~\ref{théorème: il existe une branche C^1 de solitons proche de c_s} with the existence of minimizing travelling waves stated in Section~\ref{section: minimizing property of the travelling waves}, we can exhibit a unique branch of solitons with transonic speed that are orbitally stable. This allows us to consider the orbital stability of a chain of solitons the speeds of which are taken in this regime. That\footnote{In the following definition and throughout this article, we will earmark $\gc,\ga$ for $N$-dimensional vectors and script letters for real numbers.} is why we introduce the set of well ordered admissible speeds

\begin{equation*}\label{définition ensemble des vitesses admissibles}
    \adm_N:=\big\{\gc:=(c_1,...,c_N)\in (c_0,c_s)^N\big| \ c_1 < ...< c_N\big\}.
\end{equation*}

Each speed is here chosen such that the corresponding travelling wave is unique and orbitally stable. Another crucial condition for the orbital stability to hold lies in the initial length between the solitons. This induces to define the set

\begin{equation*}\label{définition ensemble des positions a_k en L}
    \pos_N(L):=\big\{ \ga:=(a_1,...,a_N)\in\R^N\big| a_{k+1}-a_k >L,\ \forall k\in\{1,...,N-1\}\big\},
\end{equation*}

for $L>0$. For $\gc=(c_1,...,c_N)\in \adm_N$, and $\ga=(a_1,...,a_N)\in\R^N$, we finally define the sum of solitons
\begin{equation*}
    R_{\gc,\ga}:=(\eta_{\gc,\ga},v_{\gc,\ga})=\sum_{i=1}^N Q_{c_i,a_i}.
\end{equation*}

With these notations, we can state our main result.

\begin{thm}\label{théorème: stabilité orbitale}
    Let $\gc^*\in\adm_N$. There exists $\alpha_*,L_*,A_*,\tau_* >0$, such that the following holds. If $Q_0=(\eta_0,v_0)\in\Nenergysethydro$ is such that for some $\ga^0:=(a_1^0,...,a_N^0)\in \pos_N (L_0)$ with $L_0\geq L_*$,

\begin{equation*}
    \alpha_0:=\normX{Q_0-R_{\gc^*,\ga^0}}\leq \alpha_*,
\end{equation*}

then, the unique solution $Q(t)=\big(\eta(t),v(t)\big)$ to~\eqref{NLShydro} associated with the initial datum $(\eta_0,v_0)$ is globally defined and there exists $(\ga,\gc)\in\mathcal{C}^1(\R_+,\R^{2N})$ such that for any $t\in\R_+$,

\begin{equation*}
    \normX{Q(t)-R_{\gc^*,\ga(t)}}\leq A_* K(\alpha_0,L_0),
\end{equation*}

where $K(\alpha_0,L_0)= \alpha_0+e^{-\tau_* L_0 }$. Regarding the modulation parameters, we have

\begin{equation*}
    \normR{ \ga'(t)-\gc(t)}+\normR{ \gc'(t)}\leq A_* K(\alpha_0,L_0).
\end{equation*}

\end{thm}

\begin{rem}
    Since every norm on a finite dimensional Banach space are topologically equivalent, we use, here as in the sequel, one unique notation to designate every norm on $\R^N$ or $M_N(\R)\sim\R^{N^2}$. We shall indeed write $\normR{\mathfrak{x}}$ for any $\mathfrak{x}:=(x_1,...,x_N)\in\R^N$.
\end{rem}

\subsection{Sketch of the proof}\label{section: sketch of the proof}

The proof can be summarized in the following several steps.
\subsubsection{Coercivity around the soliton $Q_c$}\label{subsection: linearization + coercivité à la weinstein + preuve de régularité}
As it is stated in Section~\ref{section: minimizing property of the travelling waves}, being a solution to~\eqref{TWC} means that the Euler-Lagrange equation $\nabla (E-cp)(Q_c)=0$ is satisfied. We introduce the bilinear form $\mathcal{H}_c:=\nabla^2 (E-cp)(Q_c)$ defined on $H^2(\R)\times L^2(\R)$ and $H_c(\varepsilon):=\psLdeux{\mathcal{H}_c(\varepsilon)}{\varepsilon}$ the corresponding quadratic form, which is well-defined on $\energysethydro$. By linearizing~\eqref{TWChydro}, we obtain the following expression
\begin{align*}\label{calcul de mathcal(H_c)}
    \mathcal{H}_c=\begin{pmatrix} \mathcal{L}_c & -\dfrac{c}{2(1-\eta_c)}\\
    -\dfrac{c}{2(1-\eta_c)} & 1-\eta_c
    \end{pmatrix}
\end{align*}
with

\begin{equation*}
    \mathcal{L}_c(\varepsilon_\eta)=-\Big(\dfrac{\varepsilon_\eta'}{4(1-\eta_c)}\Big)'-\dfrac{c^2+2F(1-\eta_c)+2(1-\eta_c)f(1-\eta_c)+2(1-\eta_c)^2 f'(1-\eta_c)}{4(1-\eta_c)^2}\varepsilon_\eta.\\
\end{equation*}

The quadratic form $H_c$ is coercive under the following orthogonality constraints (see~\cite{Weinste1} for similar arguments).

\begin{prop}\label{proposition: coercité de H_c sous conditions d'orthogonalité}
    Let $c\in (0,c_s)$. There exists $l_c>0$ such that for all $\varepsilon\in \energysethydro$, satisfying the orthogonal conditions 
    \begin{equation}\label{proposition condition d'orthogonalité sur varepsilon}
        0=\psLdeuxLdeux{\varepsilon}{Q_c'}=\nabla p(Q_c).\varepsilon,
    \end{equation}
    then
\begin{equation}\label{coercivité H_c(varepsilon) geq l_c norm varepsilon}
    H_c(\varepsilon)\geq l_c\normX{\varepsilon}^2.
\end{equation}

\end{prop}

The proof of Proposition~\ref{proposition: coercité de H_c sous conditions d'orthogonalité} requires to deal with derivatives like $\partial_c Q_c$. This explains why we first need to verify (see Section~\ref{section: existence of a branch of smooth travelling waves}) that the branch of solitons of speed $c$ is smooth. Apart from the smoothness of the considered functions, the proof of Proposition~\ref{proposition: coercité de H_c sous conditions d'orthogonalité} is reminiscent from~\cite{BetGrSm1}. We obtain a formula for the essential spectrum that generalize the one in~\cite{BetGrSm1}, that is
\begin{equation*}
    \spess(\mathcal{H}_c)\subset \left[\dfrac{c_s^2-c^2}{1+c_s^2 + \sqrt{(1-c_s^2)^2+4c^2}},+\ii\right).
\end{equation*}

Furthermore, we observe that the linearized operator $\mathcal{H}_c$ owns a unique negative direction and a unique vanishing direction. Eliminating these directions by adding orthogonality conditions~\eqref{proposition condition d'orthogonalité sur varepsilon} eventually leads to~\eqref{coercivité H_c(varepsilon) geq l_c norm varepsilon}. We refer to the proof of Proposition~1 in~\cite{BetGrSm1} for more details. Spectral considerations allow us to precise the choice of $l_c$, that can be taken as follows
\begin{equation}\label{définition l_c}
l_c:=\inf_{\substack{\varepsilon\in\energysethydro\setminus \{0\}\\
    \varepsilon\text{ satisfies~\eqref{proposition condition d'orthogonalité sur varepsilon}.    }}}\dfrac{|H_c(\varepsilon)|}{\normX{\varepsilon}^2}>0.
\end{equation}

\subsubsection{Almost minimizing property of a sum of solitons}\label{subsection: almost minimizing property d'une somme de solitons}

A first consequence of Proposition~\ref{proposition: coercité de H_c sous conditions d'orthogonalité} is that, if $\varepsilon$ satisfies the orthogonal condition~\eqref{proposition condition d'orthogonalité sur varepsilon}, then by the minimizing property of $Q_{c,a}$, we get

\begin{equation}\label{propriété de minimization d'un soliton pour E-cp}
    (E-cp)(Q_{c,a}+\varepsilon)\geq (E-cp)(Q_{c,a})+l_c\normX{\varepsilon}^2 + \mathcal{O}\big(\normX{\varepsilon}^3\big).
\end{equation}

Now, we ought to find the same kind of estimate when we perturb $\Rca$ instead of $Q_{c,a}$. Let $\gc:=(c_1,...,c_N)\in\adm_N$ and $\ga\in\pos_N(L)$ with some $L >0$. Let us set
\begin{align}\label{définition mathcalO^perp_gc,ga(L)}
    \Operpca:=\left\{ \Rca +\varepsilon\in\energysethydro\bigg| \begin{array}{l}
    \psLdeuxLdeux{\varepsilon}{Q_{c_k,a_k}}=\nabla p(Q_{c_k,a_k}).\varepsilon = 0, \forall k\in\{1,...,N\}\\
    \ga\in\pos_N(L)
\end{array}
\right\},
\end{align}

and consider a perturbation $Q\in\Operpca$. Instead of considering $E-c p$ for a single soliton of speed $c$, let $\gc^*:=(c_1^*,...,c_N^*)$ and let us construct a function
\begin{equation}\label{définition de la fonctionnelle G}
    G(Q):=E(Q)-\sum_{k=1}^N c_k^* p_k(Q).
\end{equation}
that resembles $E-c_k^* p$ around each soliton $Q_{c_k^*}$, and that will eventually satisfy some coercivity property around the chain of solitons $\Rca$. We refer to~\eqref{définition: localized momentum autour de Q_c_k} for the definition of $p_k$.\newpage
With this goal in mind, we introduce the two following partitions of unity. Set $\tau >0$. For $k\in\{1,...,N\}$, we define the functions \begin{equation*}
    \Phi_k(x):=\dfrac{1}{2}\bigg(\tanh\Big(\tau\big(x-a_k+\frac{L}{4}\big)\Big)-\tanh\Big(\tau\big(x-a_k-\frac{L}{4}\big)\Big)\bigg),
\end{equation*}

and \begin{equation*}
    \Phi_{k,k+1}(x)=\left\{
\begin{array}{l}
    \dfrac{1}{2}\bigg(1-\tanh\Big(\tau\big(x-a_1+\frac{L}{4}\big)\Big)\bigg)\quad\text{if }k=0, \\
    \dfrac{1}{2}\bigg(\tanh\Big(\tau\big(x-a_k-\frac{L}{4}\big)\Big)-\tanh\Big(\tau\big(x-a_{k+1}+\frac{L}{4}\big)\Big)\bigg)\quad\text{if }k\in\{1,...,N-1\}, \\
    \dfrac{1}{2}\bigg(1+\tanh\Big(\tau\big(x-a_N-\frac{L}{4}\big)\Big)\bigg)\quad\text{if }k=N.\\
\end{array}
\right.
\end{equation*}

We have, by construction,

\begin{equation}\label{partition de l'unité Phi_k}
    \sum_{k=1}^N \Phi_k + \sum_{k=0}^N \Phi_{k,k+1} = 1,
\end{equation}

and we define

\begin{equation*}
    \varepsilon_k:=\varepsilon (.+a_k)\sqrt{\Phi_k(.+a_k)}\quad\text{and}\quad \varepsilon_{k,k+1}:=\varepsilon \sqrt{\Phi_{k,k+1}}.
\end{equation*}

In addition, choose $\tau_0>0$ and for $k\in\{1,...,N+1\}$, set
\begin{equation*}
    \chi_{k}(x)=\left\{
\begin{array}{l}
    1\quad\text{if }k=1, \\
    \dfrac{1}{2}\bigg(1+\tanh\left(\tau_0\Big(x-\frac{a_k+a_{k-1}}{2}\Big)\right)\bigg)\quad\text{if }k\in\{2,...,N\}, \\
    0\quad\text{if }k=N+1.\\
\end{array}
\right.
\end{equation*}

Note that we have $$\sum_{k=1}^N (\chi_{k}-\chi_{k+1}) = 1.$$

The choice of $\tau_0$ and $\tau$ shall be fixed later. Namely $\tau$ shall be fixed in the proof of Corollary~\ref{corollaire: controle de la fonctionnelle G} and $\tau_0$ in the proof of Proposition~\ref{proposition: monotonicity}. This choice shall only depend on the nonlinearity $f$ and $\gc^*$, that is why we allow us to shrink their value throughout the article. In order to focus on each solitons, we define 
\begin{equation}\label{définition: localized momentum autour de Q_c_k}
    p_k(Q):=\int_\R \pi(Q)(\chi_k-\chi_{k+1}).
\end{equation}

The quantity defined in~\eqref{définition: localized momentum autour de Q_c_k} can be interpreted as a  momentum localized around the soliton $Q_{c_k,a_k}$. Indeed, provided that the $a_k$ are sufficiently far away from each other, we can formally write that $\chi_k(x)-\chi_{k+1}(x)\sim 1$ for $x$ close to $a_k$ and $\chi_k(x)-\chi_{k+1}(x)\sim 0$ for $x$ close to $a_j$ for $j\neq k$. Plugging this into the expression of $p_k$, this could be regarded as $p_k(Q)=p(Q)$ if $Q\sim Q_{c_k,a_k}$. We also define
\begin{equation}\label{définition nu_c}
    \nu_\gc:=\min\Big\{\nu_{c_k}:=\sqrt{c_s^2-c_k^2}\Big|k\in \{1,...,N\}\Big\},
\end{equation}

and \begin{equation}\label{définition de l_c}
    l_\gc:=\min\big \{l_{c_k}\big|k\in \{1,...,N\}\big\},
\end{equation}
according to Proposition~\ref{proposition: coercité de H_c sous conditions d'orthogonalité}. We now write the precise approximation of the energy and the (localized) momentum of a perturbed chain of solitons $Q\in\Operpca$, when $\ga\in\pos_N(L)$. Henceforth, we introduce the notation $g=\grandOde{h}$ that means that there exists a constant $C>0$ that only depends on the parameter $\gc^*$ introduced above such that $|g|\leq C|h|$. 

\begin{prop}\label{proposition: développement de E en Q=R_c,a + varepsilon}
    For $Q\in\Operpca$, and $\tau >0$ such that $2\tau <\nu_\gc$, we have
    \begin{align}
        E&(Q)=\sum_{k=1}^N E(Q_{c_k})+\dfrac{1}{2}\bigg( \sum_{k=1}^N\nabla^2 E(Q_{c_k})(\varepsilon_k,\varepsilon_k)+\sum_{k=0}^N\nabla^2 E(0)(\varepsilon_{k,k+1},\varepsilon_{k,k+1})\bigg)\notag\\
        &+\mathcal{O}\Big(\Lambda(L,\gc)e^{-a_d\nu_{\gc}L}\Big)+\mathcal{O}\big(\normX{\varepsilon} \Lambda(L,\gc)^{\frac{1}{2}}e^{-a_d\nu_{\gc}L}\big)+\mathcal{O}(\tau\normX{\varepsilon}^2)+\mathcal{O}\Big(\normX{\varepsilon}^2 e^{-\frac{a_d\tau L}{2}} \Big)+\mathcal{R}_{\gc,\ga}(\varepsilon),\notag
    \end{align}
    where \begin{equation}\label{définition: mathcal R varepsilon}
        \mathcal{R}_{\gc,\ga}(\varepsilon)=\int_0^1 \dfrac{(1-t)^2}{2}\nabla^3 E( \Rca+t\varepsilon)(\varepsilon,\varepsilon,\varepsilon)dt,
    \end{equation}
    and
    \begin{equation*}\label{definition Lambda(L,c)}
        \Lambda(L,\gc)=\dfrac{1}{\nu_{\gc}}+L.
    \end{equation*}
\end{prop}

\begin{prop}\label{proposition: développement de p_k en Q=R_c,a + varepsilon}
    For $Q\in\Operpca$, and $\tau,\tau_0 >0$ such that $\tau_0<2\tau <\nu_\gc$, we have
    \begin{align*}
        p_k&(Q) = \  p(Q_{c_k})+\dfrac{1}{2}\Big(\nabla^2 p (Q_{c_k}).(\varepsilon_k,\varepsilon_k) + \nabla^2 p_k (0).(\varepsilon_{k,k+1},\varepsilon_{k,k+1})  + \nabla^2 p_k (0).(\varepsilon_{k-1,k},\varepsilon_{k-1,k}) \Big)\\
        & +\mathcal{O}\Big(\Lambda(L,\gc)(e^{-a_d\nu_{\gc}L}+e^{-a_d\tau_0 L})\Big)+\mathcal{O}\big(\normX{\varepsilon} \Lambda(L,\gc)^{\frac{1}{2}}(e^{-a_d\nu_{\gc}L}+e^{-a_d\tau_0 L})\big)+\mathcal{O}\Big(\normX{\varepsilon}^2 e^{-\frac{a_d\tau L}{2}} \Big).
    \end{align*}
\end{prop}

\subsubsection{Orthogonal decomposition}\label{subsection: orthogonal decomposition}
Now, we fix $\gc^*$ and set
\begin{equation}\label{définition tubular}
    \mathcal{U}_{\gc^*}(\alpha,L):=\left\{ Q\in\energysethydro\big| \inf_{\ga\in\pos_N(L)}\normX{Q-R_{\gc^*,\ga}} <\alpha \right\}.
\end{equation}
We decompose any function $Q\in\mathcal{U}_{\gc^*}(\alpha,L)$ as a chain of solitons, the speeds and localisations of which depend smoothly on $Q$, plus an extra perturbation. Here, the speeds and localisations are constructed implicitly such that the orthogonal conditions~\eqref{proposition condition d'orthogonalité sur varepsilon} hold true. We shall also see that these implicit functions, called \textbf{modulation parameters}, are controlled in accordance of the parameters $(\alpha,L)$ of the set $\mathcal{U}_{\gc^*}(\alpha,L)$. We first introduce a notation that will be used throughout the article. For real numbers $\alpha,\alpha',L$ and $L'$, we define

\begin{equation*}
    (\alpha,L)\prec (\alpha',L')\ \Longleftrightarrow\  \alpha\leq\alpha'\text{ and } L\geq L'.
\end{equation*}

\begin{prop}\label{proposition: décomposition orthogonale sans dynamique}
    Set $\gc^*=(c_1^*,...,c_N^*)\in\adm_N$. There exist constants $\alpha_1,L_1 ,K_1>0$, only depending on $\gc^*$, and two functions $$(\mathfrak{C},\mathfrak{A})=\big((c_1,...,c_N),(a_1,...,a_N)\big)\in\mathcal{C}^1\big( \mathcal{U}_{\gc^*}(\alpha_1,L_1),\R^{2N}\big)$$ such that
    for any $Q=(\eta,v)\in \mathcal{U}_{\gc^*}(\alpha_1,L_1)$, the perturbation

    \begin{equation}\label{décomposition de varepsilon sans dynamique}
        \varepsilon:=\varepsilon(Q):=Q-R_{\mathfrak{C}(Q),\mathfrak{A}(Q)}
    \end{equation}
satisfies the orthogonality conditions
\begin{equation}\label{condition d'orthogonalité sur varepsilon dans décomposition orthogonale}
    \psLdeuxLdeux{\varepsilon}{\partial_x Q_{c_k,a_k}}=\nabla p(Q_{c_k,a_k}).\varepsilon =0,
\end{equation}
for any $k\in\{1,...,N\}$. Moreover, if for some $(\alpha,L)\prec (\alpha_1,L_1)$, and some $\ga^*\in\pos_N(L)$, we have 
\begin{equation}\label{estimation normX{(ta,v)-R_{athfrak{C}(ta,v),athfrak{A}(eta,v)}}leq lpha}
    \normX{Q-R_{\gc^*,\ga^*}}\leq \alpha,
\end{equation}
then \begin{equation}\label{normX{varepsilon(eta,v)}+Vertathfrak{C}(ta,v)-gc^*Vert_1 + Vert athfrak{A}}
    \normX{\varepsilon(Q)}+\normR{ \mathfrak{C}(Q)-\gc^*} +\normR{ \mathfrak{A}(Q)-\ga^*} \leq K_1\alpha. 
\end{equation}
\end{prop}

Also, we obtain some uniform control (with respect to $Q=(\eta,v)$) on several elementary positive quantities, defined for $\gc\in \adm_N $, e.g.~\eqref{définition nu_c} and~\eqref{définition de l_c}, but also
\begin{equation*}\label{définition de mu_c}
    \mu_\gc:=\min\big\{ c_k-c_0\big|k\in \{1,...,N\}\big\},
\end{equation*}
and
\begin{equation*}\label{définition kappa_c}
    \kappa_{\gc}:=\min\left\{ - \dfrac{d}{dc}\Big(p(Q_c)\Big)_{|c=c_k}\Big| k\in\{1,...,N\}\right\}.
\end{equation*}

In this direction, we need to introduce first the following lemma, that provides a uniform bound from below between a chain of solitons and the constant function $1$. The proof of this lemma is in Appendix~A.

\begin{lem}\label{lemme: norm Lii de eta_c^* inférieur à delta}
    There exist $(\alpha_2,L_2)\prec (\alpha_1,L_1)$ and $\beta^*\in (0,1)$ such that for any $(\alpha,L)\prec (\alpha_2,L_2)$ and $\ga\in\pos_N(L)$, we have \begin{equation}
        \alpha+\beta^* < 1\quad\text{and}\quad \normLii{\eta_{\gc^*,\ga}}\leq \beta^*.
    \end{equation} 
\end{lem}

Now, we state the uniform bounds on $\mu_{\gc(Q)},\nu_{\gc(Q)},\kappa_{\gc(Q)}$ and $l_{\gc(Q)}$.

\begin{cor}\label{corollaire: décomposition orthogonale}
    Under the orthogonality conditions of Proposition~\ref{proposition: décomposition orthogonale sans dynamique}, and taking possibly a larger $L_2$ and a smaller $\alpha_2$, we have

    \begin{equation}\label{mathcalU_gc^*(alpha_1,L_1)subsetNenergysethydro}
        \mathcal{U}_{\gc^*}(\alpha_2,L_2)\subset\Nenergysethydro.
    \end{equation}

    Moreover, there exists $l_{*}>0$ such that for all $Q=(\eta,v)\in \mathcal{U}_{\gc^*}(\alpha_2,L)$ with $L\geq L_2$, we have
    \begin{equation}\label{1-etageq dfracmu_gc^8}
        1-\eta\geq \dfrac{1-\beta^*}{2},
    \end{equation}
    \begin{equation}\label{mathfrakA(eta,v)inpos_N(L-1)}
        \mathfrak{A}(Q)\in\pos_N(L-1),
    \end{equation}
and

    \begin{equation}\label{nu_athfrak(eta,v)geq dfracnu_gc^*}
        \nu_{\mathfrak{C}(Q)}\geq \dfrac{\nu_{\gc^*}}{2},\  \mu_{\mathfrak{C}(Q)}\geq \dfrac{\mu_{\gc^*}}{2},\ \kappa_{\mathfrak{C}(Q)}\geq \dfrac{\kappa_{\gc^*}}{2},\ l_{\mathfrak{C}(Q)}\geq l_*.
    \end{equation}

\end{cor}

The previous decomposition can be partially summarized as the fact that for any $L\geq L_2\geq L_1+1$,

\begin{equation}\label{mathcalUgc^*(alpha_1,L)subset mathcalO^perp}
    \mathcal{U}_{\gc^*}(\alpha_1,L)\subset \mathcal{U}^\perp_{\mathfrak{C}(Q),\mathfrak{A}(Q)}(L-1)\subset \mathcal{U}^\perp_{\mathfrak{C}(Q),\mathfrak{A}(Q)}(L_1).
\end{equation}

From now on, we impose $\tau_0<2\tau<\frac{\nu_{\gc^*}}{2}$. As a consequence of this choice, we can highlight the almost minimizing property of a chain of solitons $R_{\gc,\ga}$ for the functional $G$ defined in~\eqref{définition de la fonctionnelle G}.

\begin{cor}\label{corollaire: controle de la fonctionnelle G}
    There exists $\widetilde{l}_*>0$ such that for any $Q\in\mathcal{U}_{\gc^*}(\alpha,L)$ decomposed as $Q=R_{\mathfrak{C}(Q),\mathfrak{A}(Q)}+\varepsilon$ according to Proposition~\ref{proposition: décomposition orthogonale sans dynamique}, with $(\alpha,L)\prec (\alpha_2,L_2)$, we have

    \begin{equation*}        
    \sum_{k=1}^N \big(E(Q_{c_k^*})-c_k^* p(Q_{c_k^*})\big)+\dfrac{\widetilde{l}_*}{4}\normX{\varepsilon}^2+\mathcal{O}(\normX{\varepsilon}^3)+\mathcal{O}\big(\normR{\mathfrak{C}(Q)-\gc^*}^2\big)+\mathcal{O}\Big(Le^{-a_d\tau_0 L}\Big)\leq G(Q),
    \end{equation*}

    and \begin{equation*}
        G(Q)\leq\sum_{k=1}^N \big(E(Q_{c_k^*})-c_k^* p(Q_{c_k^*})\big)+\mathcal{O}(\normX{\varepsilon}^2)+\mathcal{O}\big(\normR{\mathfrak{C}(Q)-\gc^*}^2\big)+\mathcal{O}\Big(Le^{-a_d\tau_0 L}\Big).
    \end{equation*}
\end{cor}

\subsubsection{Evolution in time}\label{section: sketch of the proof evolution in time}
In the two following subsections, we deal with the time evolution of the perturbation which will eventually conclude the proof of the main result. The first step consists in implementing this evolution in the orthogonal decomposition of Section~\ref{subsection: orthogonal decomposition}. Whenever we take $(\alpha,L)\prec (\alpha_ 1,L_1)$, there exists, by local wellposedness, a solution $Q$ to~\eqref{NLShydro} in $\mathcal{C}^0\left([0,T],\Nenergysethydro\right)$ with $0 < T=T(\alpha,L) < T_{\max}$ such that for any $t\in [0,T]$, 
\begin{equation}
    Q(t)\in\mathcal{U}_{\gc^*}(\alpha,L).
\end{equation}

Then Proposition~\ref{proposition: décomposition orthogonale sans dynamique} enables us to define the functions \begin{equation*}
    \gc(t):=\big(c_1(t),...,c_N(t)\big):=\mathfrak{C}\big(\eta(t),v(t)\big)\text{ and }\ga(t):=\big(a_1(t),...,a_N(t)\big):=\mathfrak{A}\big(\eta(t),v(t)\big),
\end{equation*}

and also 
\begin{equation*}
    \varepsilon(t):=\big(\eta(t),v(t)\big)-R_{\gc(t),\ga(t)},
\end{equation*}
which depend continuously on $t$ and satisfy the orthogonality conditions~\eqref{condition d'orthogonalité sur varepsilon dans décomposition orthogonale}. Bear also in mind that Corollary~\ref{corollaire: décomposition orthogonale} holds, so that we have uniform bounds for $\nu_{\gc(t)},\mu_{\gc(t)}$, etc. We write $E(t):=E\big(Q(t)\big)$, $p_k(t):=p_k\big(Q(t)\big)$ and $G(t):=G\big(Q(t)\big)$ for $t\in [0,T]$. We shall see that up to taking further $(\alpha_2,L_2)$ we can improve the smoothness of $\gc,\ga$ and give a uniform control in time of these functions. This shall make a crucial use of the fact that the speed $\gc^*:=(c_1^*,...,c_N^*)\in\adm_N$ \textbf{is ordered}. In this sense, we define \begin{equation*}
    \sigma^*:=\min\left\{ c_{k+1}^*-c_k^*\ \big| \ k\in\{ 1,...,N-1\}\right\}>0.
\end{equation*}

\begin{prop}\label{proposition: controle de ga'(t)-gc(t) + gc'(t)}
    There exists $(\alpha_3,L_3)\prec(\alpha_2,L_2)$ such that if $(\alpha,L)\prec (\alpha_3,L_3)$, then $(\gc,\ga)\in\mathcal{C}^1([0,T],\adm_N\times\R^N)$ and we have for any $t\in [0,T]$
    \begin{equation}\label{Vert ga'(t)-gc(t)Vert_1 + Vert gc'(t)Vert_1 =mathcal(normXvarepsilon(t))}
        \normR{ \ga'(t)-\gc(t)} + \normR{\gc'(t)} =\mathcal{O}(\normX{\varepsilon(t)})+\mathcal{O}\big(Le^{-\frac{a_d\nu_{\gc^*}L}{2}}\big).
    \end{equation}

Moreover, we also have for any $t\in [0,T]$ and $k\in\{1,...,N-1\}$,
\begin{equation}\label{a_k+1(t)-a_k(t)geq a_k+1(0)-a_k(0)+igma^* t >L-1+sigma^* t}
    a_{k+1}(t)-a_k(t)\geq a_{k+1}(0)-a_k(0)+\sigma^* t >L-1+\sigma^* t,
\end{equation}

and for any $k\in\{1,...,N\}$,
\begin{equation}\label{sqrt c_s^2-ig(a_k'(t)big)^2geq dfracnu}
    \sqrt{c_s^2-\big(a_k'(t)\big)^2}\geq \dfrac{\nu_{\gc^*}}{2}.
\end{equation}

\end{prop}

Another crucial ingredient in the proof is a monotonicity formula for the momentum.
More precisely, we consider now for $k\in\{1,...,N\}$, the quantity 
\begin{equation*}\label{définition de p_k}
    \widetilde{p}_k(Q):=\int_\R \pi(Q)\chi_k,
\end{equation*}
and $\widetilde{p}_k(t):=\widetilde{p}_k\big(Q(t)\big)$.

\begin{rem}\label{remarque: interprétation de widetilde p_k}
    Whereas $p_k(Q)$ could be read as the amount of momentum provided by the single soliton $Q_{c_k}$, $\widetilde{p}_k(Q)$ approximates the amount of momentum provided by the sum of solitons $Q_{c_k,a_k},...,Q_{c_N,a_N}$. This claim is all the more consistent from the property that \begin{equation*}
        \widetilde{p}_k(Q)=\sum_{j=k}^N p_j(Q).
    \end{equation*}
In particular, $\widetilde{p}_1(Q)=p(Q)$.
\end{rem}

\begin{rem}
    Since the modulation parameter $\ga$ now depends on time, so do implicitly the functions $\chi_k$ for any $k\in\{1,...,N\}$.
\end{rem}

The conservation of the momentum implies that $\widetilde{p}_1$ is conserved. Regarding $\widetilde{p}_k$ for $k\geq 2$, and provided again that $c_1^*<...<c^*_N$, we can prove the following monotonicity result.
\begin{prop}\label{proposition: monotonicity}
    There exists $(\alpha_4,L_4)\prec (\alpha_3,L_3)$ such that, if $(\alpha,L)\prec (\alpha_4,L_4)$, $\widetilde{p}_k$ is differentiable and we have for any $t\in [0,T]$ and $k\in\{1,...,N\}$,

    \begin{equation}\label{formule de monotonie sur p_k}
        -\dfrac{d}{dt}\big(\widetilde{p}_k(t)\big)\leq \mathcal{O}\Big(Le^{-a_d\tau_0 (L-1+\sigma^* t)}\Big), 
    \end{equation}

    and \begin{equation}\label{formule de monotonie sur G}
        \dfrac{d}{dt}\big(\mathcal{G}(t)\big)\leq \mathcal{O}\Big(Le^{-a_d\tau_0 (L-1+\sigma^* t)}\Big).
    \end{equation}
\end{prop}

\subsubsection{Proof of the orbital stability}
\label{section: proof of the orbital stability}

\begin{figure}[!b]
    \centering
    \begin{tikzpicture}
\draw[->] (-1,0) -- (7,0) node[below]{$t$};

\filldraw[black] (0,0) circle (2pt) node[below]{$0$};

\filldraw[black] (5,0) circle (2pt) node[below]{$T^*$};

\filldraw[black] (6,0) circle (2pt) node[below]{$T^*+\theta$};

\draw (0,3.5) arc (90:270:0.3 and 1.5);
\draw[dashed] (0,0.5) arc (-90:90:0.3 and 1.5);
\draw (5,2) ellipse (0.3 and 1.5);

\draw[red] (6,2) ellipse (0.3 and 1.5);

\filldraw[black] (0,0.5) -- (5,0.5);
\filldraw[black] (0,3.5) -- (5,3.5);
\filldraw[red] (5,0.5) -- (6,0.5);
\filldraw[red] (5,3.5) -- (6,3.5);

\node[above] at (-1.5,3.18) {$\mathcal{U}_{\gc^*}(\alpha_5,L_0-2)$};

\draw[blue] (0,2.4) arc (90:270:0.15 and 0.4);
\draw[dashed,blue] (0,1.6) arc (-90:90:0.15 and 0.4);
\draw[blue] (5,2) ellipse (0.15 and 0.4);
\draw[blue] (6,2) ellipse (0.15 and 0.4);

\filldraw[dashed,blue] (0,1.6) -- (5,1.6);
\filldraw[dashed,blue] (5,1.6) -- (6,1.6);

\filldraw[blue,dashed] (0,2.4) -- (5,2.4);
\filldraw[blue,dashed] (5,2.4) -- (6,2.4);

\draw[black] (0,1.3) .. controls (1.5,2) and (2.5,1.55) .. (3,2.0) 
            .. controls (4,2.8) and (4.5,2.9) .. (5,3.5) -- (6,4.3);

\draw[blue] (0,1.76) .. controls (1.5,2.7) and (3.5,1.25) .. (3.8,1.7) 
            .. controls (4,1.9) and (4.5,2) .. (5,2.1);

\draw[blue] (5,2.1) -- (6.05,2.3);


\node[blue] at (-1.5,2) {$\mathcal{U}_{\gc^*}(\alpha_*,L_*)$};

\end{tikzpicture}
 
     \caption{Illustration of the continuation argument of Theorem~\ref{théorème: stabilité orbitale}.\ \small\textit{Here we represent the sets of the form $\mathcal{U}_{\gc^*}(\alpha,L)$ as cylinders with smaller size when $\alpha$ and $\frac{1}{L}$ are smaller. The black curve represents a solution $Q(t)$ with initial condition in $\mathcal{U}_{\gc^*}(\alpha_5,L_0-2)$, this solution must leave the neighborhood at some maximal time $T^*$. However, taking $(\alpha_*,L_*)$ such that the effects of the perturbative method are well controlled (see Proposition~\ref{proposition: controle de ga(t), c(t)-c^* et varepsilon(t)}), we manage to extend the solution $Q(t)$ for a larger time $T^*+\theta$ and such that the latter still lies in the initial neighborhood.}}
    \label{figure: cylindres}
\end{figure}
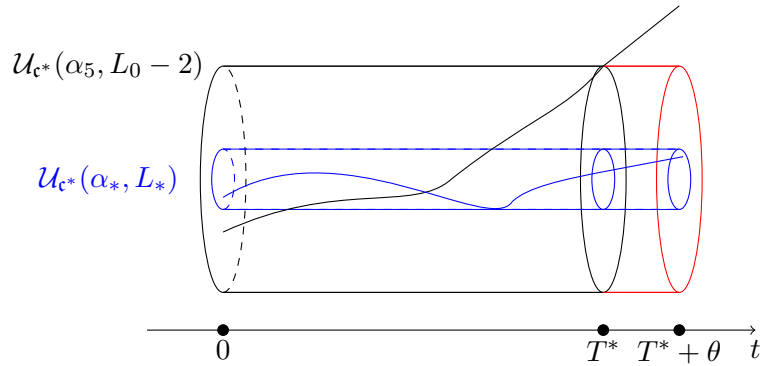

To conclude as for the orbital stability, we are going to use the control of the perturbation $\varepsilon(t)$ in terms of the distance between each solitons at time $t\in [0,T]$, the initial perturbation $\varepsilon(0)$ and $\normR{\gc(0)-\gc^*}$. That is the aim of the following proposition, the proof of which is at the end of Section~\ref{section: dynamics of the modulation parameters}.\\

\begin{prop}\label{proposition: controle de ga(t), c(t)-c^* et varepsilon(t)}
    There exists $(\alpha_5,L_5)\prec (\alpha_4,L_4)$ such that if $(\alpha,L)\prec (\alpha_5,L_5)$, we have for any $t\in [0,T]$,
    \begin{equation}\label{ga(t) in pos_N(L-1)}
        \ga(t)\in\pos_N(L-1),
    \end{equation}

\begin{equation}\label{normRgc(t)-gc^*=grandOdenormXvarepsilon(0)}
    \normR{\gc(t)-\gc^*}=\grandOde{\normX{\varepsilon(0)}^2}+\normR{\gc(0)-\gc^*}+\grandOde{\normX{\varepsilon(t)}^2}+\grandOde{Le^{-a_d\tau_0 L}},
\end{equation}
    
    \begin{equation}\label{normXvarepsilon(t)^2 leq mathcalO(normXvarepsilon(0)^2)}
        \normX{\varepsilon(t)}^2 \leq \mathcal{O}(\normX{\varepsilon(0)}^2)+\normR{\gc(0)-\gc^*}+\mathcal{O}\Big(Le^{-a_d\tau_0 L}\Big).
    \end{equation}
\end{prop}

We finally prove the main result, by a \textbf{continuation} argument. An illustration is displayed in Figure~\ref{figure: cylindres} to describe what the argument consists on.

\begin{proof}[Proof of Theorem~\ref{théorème: stabilité orbitale}]
Consider an initial datum $Q_0$ such that $\normX{Q_0-R_{\gc^*,\ga^0}}:=\alpha_0$ for some $\ga^0\in\pos_N(L_0)$ with $(\alpha_0,L_0)\prec (\alpha_*,L_*)$ and $(\alpha_*,L_*)\prec(\frac{\alpha_5}{2},L_5+2)$ to be precised later. Let $Q(t)$ be the corresponding solution locally existing on some interval $[0,T]$ such that for any $t\in [0,T]$, $$Q(t)\in\mathcal{U}_{\gc^*}(\alpha_5,L_0-2).$$ 
Since $(\alpha_5,L_0-2)\prec (\alpha_i,L_i)$ for any $i\in\{1,...,5\}$, we can exhibit $\ga(t),\gc(t),\varepsilon(t)$ such that all the estimates in the results in Section~\ref{section: sketch of the proof evolution in time} hold for any $t\in [0,T]$. Then we set

\begin{equation*}
    I:=\left\{ T\geq 0\big|\ \forall t\in [0,T], Q(t)\in \mathcal{U}_{\gc^*}(\alpha_5,L_0-2)\right\}\quad\text{and}\quad T^*:=\sup I.
\end{equation*}

Suppose by contradiction that $T^*<+\ii$. By~\eqref{mathcalU_gc^*(alpha_1,L_1)subsetNenergysethydro}, we have $0\leq T^*\leq T_{\max}$. In addition, $T^* >0$, by continuity of the flow. Now we proceed in two steps, first we show that $T^*<T_{\max}$ and then that there exists $\theta >0$ such that $0<T^*+\theta<T_{\max}$ and $T^*+\theta\in I$.\\
\underline{Step 1.} By contradiction, if $T^*=T_{\max}$ in case $T_{\max}<+\ii$, $\normLii{\eta(t)}$ would tend to 1 as $t$ tends to $T^*=T_{\max}$ by local well-posedness. However, for any $t\in [0,T^*)$, there exists $\ga^t\in\pos_N(L_0-2)$ such that $\normX{Q(t)-R_{\gc^*,\ga^t}}<\alpha_5$. Thus, taking possibly a smaller $\alpha_5$ to deal with the best multiplying constant in the one-dimensional Sobolev embedding, we can suppose that $\normLii{\eta(t)-\eta_{\gc^*,\ga^t}}<\alpha_5$. Since $L_0-2 \geq L_5$, Lemma~\ref{lemme: norm Lii de eta_c^* inférieur à delta} yields
\begin{align*}
     \normLii{\eta(t)}\leq \alpha_5 + \beta^* ,
\end{align*}
and we obtain a contradiction by passing to the limit $t\rightarrow T^*$.\\
\underline{Step 2.} We use~\eqref{normRgc(t)-gc^*=grandOdenormXvarepsilon(0)} and~\eqref{normXvarepsilon(t)^2 leq mathcalO(normXvarepsilon(0)^2)} in Proposition~\ref{proposition: controle de ga(t), c(t)-c^* et varepsilon(t)} to obtain

\begin{align}
    \normX{Q(t)-R_{\gc^*,\ga(t)}}&\leq \normX{\varepsilon(t)}+\normX{R_{\gc(t),\ga(t)}-R_{\gc^*,\ga(t)}}\notag\\
    &\leq \mathcal{O}\left( \normX{\varepsilon(0)}+\normR{\gc(0)-\gc^*}+L_0 e^{-a_d\tau_0 L_0}+\big|\gc(t)-\gc^*\big|\right)\notag\\
    &\leq K_6\left( \alpha_0+ e^{-a_d \frac{\tau_0}{2}L_0}    \right)\label{estimée finale K_6}.
\end{align}

By Proposition~\ref{proposition: controle de ga'(t)-gc(t) + gc'(t)} and the choice of $\tau_0$ before Corollary~\ref{corollaire: controle de la fonctionnelle G}, we deduce 
\begin{equation}\label{estimée finale K_7}
    \normR{\ga'(t)-\gc(t)}+\normR{\gc'(t)}\leq K_7\left(\alpha_0 + e^{-a_d \frac{\tau_0}{2}L_0} \right),
\end{equation}

for any $t\in [0,T^*)$. Using successively~\eqref{ga(t) in pos_N(L-1)} and then~\eqref{estimée finale K_6}, we infer that there exists $t^n\rightarrow T^*$ with $t^n<T^*$ such that for any $n$, \begin{align*}
    \inf_{\ga\in\pos_N(L_0-2)}\normX{Q(T^*)-R_{\gc^*,\ga}}&\leq \normX{Q(T^*)-R_{\gc^*,\ga(t^n)}}\\
    &\leq \normX{Q(T^*)-Q(t^n)}+\normX{Q(t^n)-R_{\gc^*,\ga(t^n)}}\\
    &\leq  \normX{Q(T^*)-Q(t^n)}+K_6\left( \alpha_0+ e^{-a_d \frac{\tau_0}{2}L_0}   \right).
\end{align*}

Now passing to the limit $n\rightarrow +\ii$, and imposing additionally 
\begin{equation*}
     K_6\left( \alpha_*+ e^{-a_d \frac{\tau_0}{2} L_*}   \right)<\dfrac{\alpha_5}{2},
\end{equation*} we obtain that $Q(T^*)\in\mathcal{U}_{\gc^*}\left(\frac{\alpha_5}{2},L_0-2\right)$. By local wellposedness, we show the same way that there exists $\theta>0$ such that $T^*+\theta <T_{\max}$ and $Q(T^*+\theta)\in\mathcal{U}_{\gc^*}\left(\alpha_5,L_0-2\right)$ which contradicts the maximality of $T^*$. We have shown that $T^*=+\ii$, thus the solution $Q(t)$ is globally defined and the estimates~\eqref{estimée finale K_6},~\eqref{estimée finale K_7} hold for any $t\in\R_+$. We conclude by taking $\tau_*=\frac{a_d\tau_0}{2}$ and $A^*$ larger than both $K_6$ and $K_7$.
\end{proof}

\numberwithin{equation}{section}

\section{Existence of a branch of smooth travelling waves}\label{section: existence of a branch of smooth travelling waves}To prove Theorem~\ref{théorème: il existe une branche C^1 de solitons proche de c_s}, we proceed in three steps. First we cite D. Chiron's result that provides a unique local branch of solutions of~\eqref{TWC} in the transonic limit.

\begin{thm}[Theorem~4 in~\cite{Chiron7}]\label{théoreme: Chiron existence branche soliton}
    Let $f\in\mathcal{C}^3(\R_+)$ and assume that~\eqref{condition suffisante pour la stabilité orbitale sur f''(1)+6f'(1)>0} holds. There exist $\delta>0$ and $0<c_0<c_s$ such that for any $c\in (c_0,c_s)$, there exists a solution $\gv_c$ to~\eqref{TWC} such that $\normLii{1-|\gv_c|}\leq \delta$. Morevover, this solution is unique up to a translation and a constant phase shift of modulus one.
    Furthermore, we have the asymptotic estimate
    \begin{equation}\label{équivalent moment p(gv_c)}
    p(\gv_c)\underset{c\rightarrow c_s}{\sim}\dfrac{6c_s\nu_c^3}{k^2},
\end{equation}
where $k:=2f''(1)+6f'(1)$.
\end{thm}

In Subsection~\ref{subsection: amélioration de Chiron pour hydro}, we provide a slight improvement of the previous theorem, showing that the branch $c\mapsto Q_c$ is actually smooth. Once we know that there exists locally a unique smooth hydrodynamical travelling wave for $c$ close to $c_s$, we show that it decays exponentially. Using this decay at infinity, we eventually improve the smoothness of the map $c\mapsto Q_c$ to the space $\mathcal{C}^2\big((c_0,c_s),\mathcal{NX}^2(\R)\big)$.

\subsection{A priori existence of a branch of travelling waves}
\label{subsection: amélioration de Chiron pour hydro}
In our framework, $f\in\mathcal{C}^3(\R_+)$, then by Theorem~\ref{théoreme: Chiron existence branche soliton}, we exhibit $c_0$ such that for $c\in (c_0,c_s)$ there exists a unique solution $\gv_c$ to~\eqref{TWC}, up to a translation and a constant phase shift. This travelling wave lies in the energy set and does not vanish. As a consequence, we can construct the corresponding hydrodynamical travelling wave, unique up to translations, which belongs to $\Nenergysethydro$. As explained in~\cite{Chiron7}, we can freeze the invariance by translation such that the resulting hydrodynamical travelling wave $(\eta_c,v_c)$ is even. For $c\in (c_0,c_s)$, we set $X(c,Q):=\nabla E(Q)-c\nabla p(Q)$. We can check from~\eqref{expression de l'énergie} that $X$ is a $\mathcal{C}^3$ function (because $f\in\mathcal{C}^3(\R_+)$) on $\R\times \left\{(\eta,v)\in H^2(\R)\times L^2(\R)\big| \eta,v\text{ are even and }\max_{\R}\eta <1\right\}$ the subset of the Banach space $ \R\times H^2\times L^2$. By definition, we deduce $\partial_2 X(c,Q_{c}).\varepsilon = \mathcal{H}_{c}(\varepsilon)$. Since $\ker(\mathcal{H}_{c})=\spanned{\partial_x Q_{c}}$ is spanned by an odd function, then $\ker\big(\partial_2 X(c,Q_{c})\big)=\{0\}$. Since $\mathcal{H}_c$ is self-adjoint, we infer from the Fredholm alternative that $\partial_2 X(c,Q_{c})$ is in fact invertible. Thus, by the implicit function theorem, there exists a neighborhood of $c$ on which $X(c,Q_c)=0$ owns a unique solution. By uniqueness of the solution (up to translation) of \eqref{TWChydro} at fixed speed, this solution is $Q_c$ and the map $c\mapsto Q_c\in\mathcal{C}^3\big((c_0,c_s),H^2(\R)\times L^2(\R)\big)$. Moreover $\eta_c$ and $v_c$ satisfy

\begin{equation}\label{equation -eta'' = g(eta)}
    -\partial_x^2\eta_c=\dfrac{1}{2}\mathcal{N}_c'(\eta_c),
\end{equation}

where \begin{equation}\label{expression de N_c}
    \mathcal{N}_c(x)=c^2 x^2 -4(1-x)F(1-x),
\end{equation}

and

\begin{equation}\label{expression de v_c en fonction de eta_c}
    v_c=\dfrac{c\eta_c}{2(1-\eta_c)}.
\end{equation}

By standard arguments on smooth dependence (with respect to $c$) for ODEs depending on a parameter, the branch is $\mathcal{C}^3$ on $(c_0,c_s)$ with values in $\big(\mathcal{C}_{loc}^5(\R)\big)^2$ (see for instance~\cite{QuefZui0}).

\subsection{Decay at infinity}
\label{subsection: exponential decay}
By Subsection~\ref{subsection: amélioration de Chiron pour hydro}, there exists a unique branch of solution to~\eqref{TWChydro} lying in the set\begin{equation}\label{définition de l'ensemble I}
    \mathcal{I}:=\mathcal{C}^3\Big((c_0,c_s),\big(H^2(\R)\times L^2(\R)\big)\cap\big(\mathcal{C}_{loc}^3(\R)\times\mathcal{C}_{loc}^2(\R)\big)\Big).
\end{equation}

The exponents in the definition of $\mathcal{I}$ are set in accordance with establishing~\eqref{estimée décroissance exponentielle à tout ordre pour eta_c et v_c}. We state additional equations that are also satisfied by $(\eta_c,v_c)\in\Nenergyset$ in order to study the exponential decay. Integrating~\eqref{equation -eta'' = g(eta)}, we obtain

\begin{equation}\label{equation -eta'^2=nu(eta)}
    -(\partial_x\eta_c)^2=\mathcal{N}_c(\eta_c).
\end{equation}

Moreover, note that \begin{equation}\label{expression de N + c'est un grand O de x^3}
    \mathcal{N}_{c_s}(x)=c_s^2 x^2 -4(1-x)F(1-x)=\mathcal{O}(x^3).
\end{equation}

We use a result established in~\cite{LinZhiw1} (see also~\cite{Chiron7}) that gives a sufficient and necessary criterion for the existence and uniqueness of a solution satisfying~\eqref{equation -eta'' = g(eta)}. It essentially adapts a general result concerning ordinary differential equations due to H. Berestycki and P.-L. Lions (see Theorem~5 in~\cite{BereLions1}).

\begin{prop}[\cite{LinZhiw1},\cite{Chiron7}]\label{prop: existence d'une solution eta qui tend vers 0 si il existe un zero z_c}
Let $c\in (c_0,c_s)$. There exists a unique (up to translations) non constant solution $\eta_c$ of~\eqref{equation -eta'' = g(eta)} vanishing at infinity if and only if there exists
$\xi_c\in (0,1)$ such that $\mathcal{N}_c(\xi_c)=0$, $\mathcal{N}_c(x) <0$ on $(0,\xi_c)\text{ and }\mathcal{N}'_c(\xi_c)>0$. In that case, $\eta_c$ is even, reaches its maximum in $0$ and $\eta_c(0)=\xi_c$.
\end{prop}

\begin{rem}
    The zero $\xi_c$ can be supposed minimal for this property, so we can see later that $\xi_c$ tends to 0 as $c$ tends to $c_s$.
\end{rem}

We first give a bound from below (equal to $c/c_s$ in the Gross-Pitaevskii case) on $|\gv_c|$.
\begin{prop}\label{proposition: min |gv_c|=delta c/c_s}
    There exists $\delta \in (0,1)$, independent of $c\in (c_0,c_s)$, such that
    \begin{equation*}\label{le minimum de |gv_c| est c/c_s}
    \min_\R |\gv_c|=|\gv_c(0)|=\sqrt{1-\eta_c(0)}\geq \delta \dfrac{c}{c_s}.
\end{equation*}
In terms of the hydrodynamic variables, this reads
\begin{equation}\label{le maximum de eta_c est 1- delta c/c_s}
    \max_\R \eta_c=\eta_c(0)=\xi_c\leq 1-\delta^2\dfrac{c^2}{c_s^2}.
\end{equation}
\end{prop}

\begin{proof}
The function $\eta_c$ achieves its maximum in $0$ (up to a translation). Furthermore, by~\eqref{equation -eta'^2=nu(eta)}, we have $0=\mathcal{N}_c\big(\eta_c(0)\big)$, i.e. $0=c^2\eta_c(0)^2-4\big(1-\eta_c(0)\big)F\big(1-\eta_c(0)\big)$, and using a second order Taylor's formula with remainder, we obtain
    \begin{equation}\label{preuve min |gv_c|=delta c/c_s: c^2/c_s^2=...}
        \dfrac{c^2}{c_s^2}=(1-\eta_c(0))\Big( 1+\dfrac{\eta_c(0)}{2c_s^2}\int_0^1(1-t)f''\big(1-t\eta_c(0)\big)dt\Big).
    \end{equation}

In particular, we get \begin{equation*}
    \dfrac{c^2}{c_s^2}\dfrac{1}{1+\Vert f''\Vert_{L^\ii([0,1])}/4c_s^2}\leq|\gv_c(0)|^2,
\end{equation*}

which provides the estimates in Proposition~\ref{le minimum de |gv_c| est c/c_s}.
\end{proof}

\begin{rem}
    In order to simplify the notations, we occasionally write $f_1(x)\lesssim f_2(x)$ that stands for $f_1(x)=\mathcal{O}\big(f_2(x)\big)$.
\end{rem}

\begin{lem}\label{lemme: régularité de Q_c en fonction de f + exponential decay for spaciale}
We have $\eta_c\in\mathcal{C}^{3}(\R)\cap H^{3}(\R)$ and $v_c\in\mathcal{C}^{2}(\R)\cap H^{2}(\R)$, and
\begin{equation}\label{inégalité décroissance exponentielle de eta_c}
    \sum_{\substack{0\leq k_1\leq 3 \\ 0\leq k_2\leq 2}} |\partial_x^{k_1}\eta_c(x)|+c^{1+2k_2}|\partial_x^{k_2} v_c(x)|\lesssim e^{-\nu_c|x|}.
\end{equation}

\end{lem}

\begin{proof}

\underline{Step 1: Exponential decay of $\eta_c$.} The exponential decay of $\eta_c$ is proved in~\cite{Chiron8} and is stated as follows,

\begin{equation*}
    |\eta_c(x)|\lesssim M_c e^{-\nu_c|x|},
\end{equation*}
where
\begin{equation}\label{expression de M_c}
    M_c=\xi_c e^{\int^{\xi_c}_0\frac{\mathcal{N}_{c_s}(\xi)}{\sqrt{-\xi^2\mathcal{N}_c(\xi)}\big(\sqrt{-\mathcal{N}_c(\xi)}+\sqrt{\nu_c^2\xi^2}\big)}d\xi }. 
\end{equation}

We need to improve the previous estimate. In order to do that, we shall handle the dependence in $c$ in the constant $M_c$. We first prove the following intermediate result. We define $\widetilde{k}:=-\frac{3}{k}$ with $k$ defined in Theorem~\ref{théoreme: Chiron existence branche soliton}.
\begin{claim}\label{claim: xi_c lesssim varepsilon^2} 
We have
\begin{equation}\label{développement asymptotique de xi_c en nu_c}
    \xi_c \underset{c\rightarrow c_s}{=}\widetilde{k}\nu_c^2+\mathcal{O}(\nu_c^4).
\end{equation}

Moreover, 
\begin{equation}\label{calcul de partial_c xi_c}
    \partial_c \xi_c=-\dfrac{2c\xi_c^2}{\mathcal{N}_c'(\xi_c)}\underset{c\rightarrow c_s}{=}-2c_s \widetilde{k}+\mathcal{O}(\nu_c^2),
\end{equation}

and \begin{equation*}
    \partial_c^2\xi_c \underset{c\rightarrow c_s}{\sim}-2\widetilde{k}.
\end{equation*}
\end{claim}

\begin{proof}
We show that the first positive zero $\xi_c\in (0,1)$ satisfies $\xi_c=-3\nu_c^2/k +\mathcal{O}( \nu_c^3)$. We first write a Taylor expansion of $\mathcal{N}_c$ at $0$:
\begin{equation}\label{développement asymptotique de n_c(xi) en 0}
    \mathcal{N}_c(\xi)\underset{\xi\rightarrow 0}{=}-\xi^2\big(\nu_c^2+\dfrac{k}{3} \xi + \mathcal{O}(\xi^2)\big).
\end{equation}

Now we show that $\xi_c$ strictly decreases to $0$ as $c\rightarrow c_s$, so we can plug $\xi_c\neq 0$ into~\eqref{développement asymptotique de n_c(xi) en 0} and get $0\underset{c\rightarrow c_s}{=} \nu_c^2+\dfrac{k}{3}\xi_c + \mathcal{O}(\xi_c^2)$. We would then deduce simultaneously that $\xi_c=\mathcal{O}(\nu_c^2)$ and~\eqref{développement asymptotique de xi_c en nu_c}. To do this, we consider $c\mapsto \xi_c =\eta_c(0)$ which is smooth (see Subsection~\ref{subsection: amélioration de Chiron pour hydro}), and differentiating $\mathcal{N}_c(\xi_c)=0$ leads to the first identity in~\eqref{calcul de partial_c xi_c}. Since $\mathcal{N}_c'(\xi_c)>0$, this means that $\xi_c$ strictly decreases with respect to $c$. Moreover it lies in $(0,1)$ for any $c$, so that there exists $\xi_*<1$ such that \begin{equation*}
    \xi_c\underset{c\rightarrow c_s}{\longrightarrow}\xi_*.
\end{equation*}
We assume by contradiction that $\xi_*\neq 0$. Passing to the limit $c\rightarrow c_s$ into $\mathcal{N}_c(\xi_c)=0$ and using~\eqref{hypothèse de croissance sur F minorant intermediaire}, we obtain $c_s^2 \xi_*^2\geq c_s^2(1-\xi_*)\xi_*^2$, then $\xi_*\geq 1$. Thus we must have $\xi_*=0$. To conclude the proof, we differentiate $\partial_c\xi_c$ with respect to $c$ and use~\eqref{développement asymptotique de xi_c en nu_c} and~\eqref{calcul de partial_c xi_c}. 
\end{proof}

We can control the integral in~\eqref{expression de M_c}, by using successively~\eqref{expression de N + c'est un grand O de x^3}, Proposition~\ref{prop: existence d'une solution eta qui tend vers 0 si il existe un zero z_c}, the substitution $\xi=\nu_c^2\eta$ and~\eqref{développement asymptotique de n_c(xi) en 0},
\begin{align*}
    \bigg|\int^{\xi_c}_0\dfrac{\mathcal{N}_{c_s}(\xi)}{\sqrt{-\xi^2\mathcal{N}_c(\xi)}(\sqrt{\mathcal{-N}_c(\xi)}+\sqrt{\nu_c^2\xi^2})}d\xi\bigg|& \leq \int_0^{\xi_c}\dfrac{|\mathcal{N}_{c_s}(\xi)|}{\sqrt{-\xi^4 \nu_c^2 \mathcal{N}_c(\xi)}}d\xi \lesssim \int_0^{\xi_c}\dfrac{\xi}{\nu_c\sqrt{-\mathcal{N}_c(\xi)}}d\xi\\
    &=\nu_c^3\int_0^{\frac{\xi_c}{\nu_c^2}}\dfrac{\eta d\eta}{\sqrt{-\mathcal{N}_c(\nu_c^2\eta)}}\underset{c\rightarrow c_s}{\sim}\int_0^{\widetilde{k}}\dfrac{d\eta}{\sqrt{1-\frac{\eta}{\widetilde{k}}}}=2\widetilde{k},
\end{align*}
where $\widetilde{k}=-\frac{3}{k}>0$. Therefore we conclude from the previous claim that,
\begin{equation}\label{|eta_c(x)|leq Knu_c^2 e^-nu_c |x|}
    |\eta_c(x)|\lesssim\nu_c^2 e^{-\nu_c |x|}.
\end{equation}

\underline{Step 2: Exponential decay of the remaining terms.} In view of~\eqref{expression de N_c}, $\mathcal{N}'_c\in\mathcal{C}^{3}((-\ii,1])$ so it is continuous. Since $\eta_c\in H^1(\R)$, by using the one-dimensional Sobolev embedding and~\eqref{equation -eta'' = g(eta)}, we obtain that $\eta_c\in\mathcal{C}^2(\R)$. In particular, by a standard bootstrap argument $\partial_x^2\eta_c\in\mathcal{C}^3(\R)$ i.e. $\eta_c\in\mathcal{C}^{5}(\R)$. Now we deal with the exponential decay. Using~\eqref{equation -eta'^2=nu(eta)} and~\eqref{expression de N + c'est un grand O de x^3}, we get $|\partial_x\eta_c|\leq \nu_c |\eta_c|+\grandOde{|\eta_c|^{\frac{3}{2}}}\lesssim |\eta_c|$. Regarding the other derivatives, we write a Taylor expansion of $\mathcal{N}_c'$, so that we obtain
    \begin{equation}\label{preuve: développement de taylor de n_c'(eta_c) à l'ordre 1}
        -\partial_x^2 \eta_c=\dfrac{\eta_c}{2}\int_0^1 \mathcal{N}_c''(t\eta_c)dt.
    \end{equation}

By~\eqref{|eta_c(x)|leq Knu_c^2 e^-nu_c |x|}, $\eta_c$ is uniformly bounded with respect to $c\in (0,c_s)$ and $x\in\R$. In view of \begin{equation*}
    \mathcal{N}_c''(t\eta_c)=2c^2+8f(1-t\eta_c)+4(1-t\eta_c)f'(1-t\eta_c),
\end{equation*} 

and since $f'$ is continuous, the integral in~\eqref{preuve: développement de taylor de n_c'(eta_c) à l'ordre 1} is uniformly bounded with respect to $x\in\R$ and $c\in (0,c_s)$. Therefore, we have shown $|\partial_x^2\eta_c|\lesssim |\eta_c|$ and thus the exponential decay of $\partial_x^2\eta_c$. We can differentiate~\eqref{equation -eta'' = g(eta)} with respect to $x$ and with a bootstrap argument, we control $\partial^3_x\eta_c$ and $\partial_x^{4}\eta_c$ in terms of $|\eta_c|$. Regarding the derivatives of $v_c$, by~\eqref{expression de v_c en fonction de eta_c} and Proposition~\ref{proposition: min |gv_c|=delta c/c_s}, we have $c|v_c|\lesssim |\eta_c|$ and 
\begin{equation*}
    \partial_x v_c=\dfrac{c}{2}\dfrac{\partial_x\eta_c}{(1-\eta_c)^2},
\end{equation*}
so that $c^3|\partial_x v_c|\lesssim |\partial_x\eta_c|$. By induction, we show that $c^{1+2k}|\partial_x^k v_c|\lesssim|\eta_c|$, for $k\in\{ 0,...,2\}$.
\end{proof}

Now we deal with the derivatives of $\eta_c$ and $v_c$ with respect to $c$. We obtain estimates similar to $(2.9)$ in~\cite{BetGrSm1}.

\begin{lem}\label{lemme: exponential decay} There exists $K_d>0$ and $a_d>0$ independent of $c$ and $x$ such that
    \begin{equation}
\sum_{\substack{0\leq k_1\leq 3 \\ 0\leq k_2\leq 2 \\ 1\leq j\leq 3}}\nu_c^{2(j-1)}\Big(|\partial_c^j\partial_x^{k_1} \eta_c(x)|+c^{1+2j+2k_2}|\partial_c^j\partial_x^{k_2} v_c(x)|\Big)\leq K_d e^{-a_d\nu_c |x|}.
\end{equation}

\end{lem}

\begin{proof}
First, we prove
    \begin{equation*}
        |\partial_c\eta_c|+c^3 |\partial_c v_c|\lesssim e^{-\frac{\nu_c x}{2}}.
    \end{equation*}
We compute
\begin{equation*}
    \partial_c v_c=\dfrac{v_c}{c}+\dfrac{c}{2}\dfrac{\partial_c\eta_c}{(1-\eta_c)^2},
\end{equation*}

then by Proposition~\ref{proposition: min |gv_c|=delta c/c_s},
we infer that $c^3\partial_c v_c$ enjoys the same decay than $c^2 |v_c|+|\partial_c\eta_c|$. We already handled the decay of $v_c$ in Lemma~\ref{lemme: régularité de Q_c en fonction de f + exponential decay for spaciale}, so it remains to deal with $\partial_c\eta_c$. We consider $K_c(x):=(\partial_c\eta_c(x))^2$. By~\eqref{equation -eta'' = g(eta)}, we obtain
\begin{equation*}
    K_c''=2(\partial_c \partial_x\eta_c)^2+2\partial_c\eta_c\partial_c\partial_x^2\eta_c \geq -K_c\mathcal{N}''_c(\eta_c)-4c\eta_c\partial_c \eta_c.
\end{equation*}

Since $\eta_c(x)$ tends to $0$ as $x$ tends to $\pm\ii$, we have
\begin{equation*}
    -\mathcal{N}''_c(\eta_c)\underset{|x|\rightarrow +\ii}{\longrightarrow}2\nu_c^2.
\end{equation*}
Then, there exists $R_* > 0$ such that for any $|x|\geq R_*$, $\nu_c^2\leq -\mathcal{N}''_c(\eta_c)$ hence $K_c''\geq \nu_c^2 K_c-4c\eta_c\partial_c\eta_c$. Thus \begin{equation*}
    K_c''-\nu_c^2 K_c\geq -4c\sqrt{2}\dfrac{\eta_c}{\nu_c}\sqrt{\dfrac{1}{2}}\nu_c\partial_c\eta_c.
\end{equation*}

Using the exponential decay of $\eta_c$ in~\eqref{|eta_c(x)|leq Knu_c^2 e^-nu_c |x|}, we get \begin{equation*}
    K_c''(x)-\dfrac{\nu_c^2}{2} K_c(x)\geq -\dfrac{8c^2\eta_c^2(x)}{\nu_c^2}\gtrsim -c^2\nu_c^2 e^{-2\nu_c |x|},
\end{equation*}

hence there exists a constant $C>0$ such that 
\begin{equation*}
    -K_c''(x)+\dfrac{\nu_c^2}{2}K_c(x)\leq C\nu_c^2 e^{-\frac{\nu_c}{\sqrt{2}} |x|}\leq C\nu_c^2 e^{-\frac{\nu_c}{2\sqrt{2}} |x|} .
\end{equation*}

We introduce a general lemma for the exponential decay of solutions to elliptic differential inequalities, the proof of which is at the end of this subsection.

\begin{lem}\label{lemme: -g''+omega^2 g à décroissance exponentielle implique g décroissance expo}
    Let $g$ be a non negative function in $\mathcal{C}^2(\R)$. Set $\omega >0$ and $0<a<1$ and assume that there exists $C>0$ such that for any $x\in\R$, we have $-g''(x)+\omega^2 g(x)\leq C\omega^2 e^{-a\omega |x|}$ and $e^{-2\omega |x|}\big(g(x)e^{\omega |x|}\big)'$ tends to $0$ as $|x|$ tends to $+\ii$. Then there exists $\widetilde{C}>0$ independent of $\omega$ such that for all $x\in\R$, 
    \begin{equation*}
        g(x)\leq \widetilde{C}e^{-a\omega |x|},
    \end{equation*}

with $\widetilde{C}=g(0)+\frac{C}{1-a^2}$       
\end{lem}
We apply this lemma with $a=1/2$, $\omega(c):=\frac{\nu_c}{\sqrt{2}}$ and $g=K_c$. It remains to check that

\begin{equation*}
    e^{-2\omega |x|}\big(K_c(x)e^{\omega |x|}\big)'\underset{|x|\rightarrow +\ii}{\longrightarrow}0.
\end{equation*}

Indeed, at fixed speed $c\in (c_0,c_s)$, $\partial_c\eta_c,\partial_c\partial_x\eta_c\in H^1(\R)$ (the branch of solitons lies in the set defined in~\eqref{définition de l'ensemble I}) so that $K_c$ and $K_c'$ tend to $0$ as $|x|$ tend to $+\ii$. Finally using~\eqref{calcul de partial_c xi_c}, we infer that $\big(\partial_c\eta_c(0)\big)^2$ is bounded by a constant that does not depend on $c$. Therefore, there exists $\widetilde{C}>0$ independent of $c$ and $x$ such that \begin{equation*}
    K_c(x)\leq \widetilde{C}e^{-\frac{\nu_c}{2\sqrt{2}}|x|},
\end{equation*}
which provides the exponential decay of $\partial_c\eta_c$.

Now we deal with the coupled derivative $\partial_c \partial_x\eta_c$. We deduce the exponential decay of $\partial_c \partial_x^2\eta_c$ from differentiating~\eqref{equation -eta'' = g(eta)} with respect to $c$ and using the exponential decay of $\partial_c\eta_c$ and $\eta_c$. The exponential decay of $\partial_c \partial_x\eta_c$ can be then derived from integrating $\partial_c \partial_x^2\eta_c$ between $x>0$ and $+\ii$ (and between $-\ii$ and $x<0$), since $\partial_c \partial_x\eta_c\in H^1(\R)$.

Next we set $L_c:=(\partial_c^2\eta_c)^2$. Similarly to $K_c$, we differentiate~\eqref{equation -eta'' = g(eta)} twice with respect to $c$, and using the exponential decay of $\eta_c$ and $\partial_c\eta_c$, we deduce
that there exists $a_d\in (0,1)$, $R_*>0$ and $C>0$ such that for $|x|\geq R_*$,
\begin{equation*}
    -L_c''+\nu_c^2 L_c\leq C e^{-a_d\nu_c|x|}.
\end{equation*}
We verify that the assumptions of Lemma~\ref{lemme: -g''+omega^2 g à décroissance exponentielle implique g décroissance expo} are satisfied, and thus we have the bound, up to taking a larger constant $C>0$,
\begin{equation*}
    |\partial_c^2\eta_c(x)|\leq \dfrac{C}{\nu_c^2} e^{-a_d\nu_c|x|},
\end{equation*}

where $C$ does not depend on $c$ because of the asymptotics for $\partial_c^2\xi_c$ in Claim~\ref{claim: xi_c lesssim varepsilon^2}. Since $f\in\mathcal{C}^3(\R_+)$, we can repeat the same type of arguments, and deduce the exponential decay of the remaining derivatives $\partial_c^j\partial_x^3\eta$ when $j\in\{0,...,2\}$ and $\partial_c^3\partial_x^k\eta_c$ when $k\in\{0,...,3\}$. We conclude the proof by deriving the decay of the derivatives $\partial^j_c\partial_x^k v_c$ from the formulae~\eqref{expression de v_c en fonction de eta_c}.
\end{proof}

We finish this section by the proof of Lemma~\ref{lemme: -g''+omega^2 g à décroissance exponentielle implique g décroissance expo}.

\begin{proof}
By hypothesis, we have for $x\geq 0$, \begin{equation*}
    -\Big(\big(g(x)e^{\omega x})'e^{-2\omega x}\Big)'e^{\omega x}=-g''(x)+\omega^2 g(x)\leq C\omega^2 e^{-a\omega x}.    
\end{equation*}

Integrating this equation between $y$ positive and $+\ii$ provides, by using the limit in the hypotheses,
\begin{equation*}
    \big(g(y)e^{\omega y}\big)'\leq \dfrac{C\omega}{a+1}e^{-(a-1)\omega y}.
\end{equation*}

We integrate between $0$ and $z$ positive and we get
\begin{equation*}
    g(z)\leq \Big(g(0)+\dfrac{C}{1-a^2}\big(e^{-(a-1)\omega y}-1\big)\Big)e^{-\omega z}\leq \widetilde{C}e^{-a\omega z}.          
\end{equation*}
The proof is similar for $z$ negative.

\end{proof}

\subsection{Proof of Theorem~\ref{théorème: il existe une branche C^1 de solitons proche de c_s}}\label{subsection: proof of theorem existence de travelling wave}
Using the decay at infinity and the a priori existence of a unique branch in $\mathcal{I}$, we now conclude and show that the branch is $\mathcal{C}^2\big((c_0,c_s),\mathcal{NX}^2(\R)\big)$. We proceed locally. Set $c_1\in (c_0,c_s)$ and $\delta >0$ such that $(c_1-\delta,c_1+\delta)\subset (c_0,c_s)$. We show that $c\mapsto\eta_c\in\mathcal{C}^0\big((c_1-\delta,c_1+\delta),H^3(\R)\big)$.
Let $\varepsilon >0$ and $c_2\in (c_1-\delta,c_1+\delta)$, we first have for all $R >0$, \begin{align*}
    \Vert\eta_{c_1}-\eta_{c_2}\Vert_{H^3}^2\leq \Vert \eta_{c_1}-\eta_{c_2}\Vert_{H^3(-R,R)}^2+\Vert \eta_{c_1}-\eta_{c_2}\Vert_{H^3(\R\setminus [-R,R])}^2.
\end{align*}

On the one hand, we have by Lemma~\ref{lemme: exponential decay}, \begin{align*}
    \Vert \eta_{c_1}-\eta_{c_2}\Vert_{H^3(\R\setminus [-R,R])}^2 & \leq\sum_{k=0}^3 2\big(\Vert\partial_x^k\eta_{c_1}\Vert_{L^2(\R\setminus [-R,R])}^2+\Vert\partial_x^k\eta_{c_2}\Vert_{L^2(\R\setminus [-R,R])}^2\big) \\
    &\leq 24K_d^2 \Vert e^{-2a_d\nu_{c_1+\delta}|x|}\Vert_ {L^2(\R\setminus [-R,R])}^2.
\end{align*}

We can find $R>0$ large enough such that $\Vert \eta_{c_1}-\eta_{c_2}\Vert_{H^5(\R\setminus [-R,R])}^2\leq \varepsilon^2/2$. On the other hand, by Lemma~\ref{lemme: régularité de Q_c en fonction de f + exponential decay for spaciale} with $j=1$ and $k\in\{0,...,3\}$ and doing a first order Taylor expansion with respect to $c$, we have in particular that 

\begin{equation*}
\big|\partial_x^k\eta_{c_1}(x)-\partial_x^k\eta_{c_2}(x)\big|= \left|\int_{c_1}^{c_2}\partial_c\partial_x^k\eta_{c}(x)dc\right|\leq K_d|c_1-c_2|.
\end{equation*}

Since 

\begin{align*}
    \Vert \eta_{c_1}-\eta_{c_2}\Vert_{H^3(-R,R)}^2 \leq \sum_{k=0}^3 2R \Vert \partial_x^k\eta_{c_1}-\partial_x^k\eta_{c_2}\Vert_{L^\ii(-R,R)}^2,
\end{align*}
for $c_2$ close enough to $c_1$, we have $\Vert \eta_{c_1}-\eta_{c_2}\Vert_{H^3(-R,R)}^2\leq \varepsilon^2/2$. In conclusion, we can find $\delta>0$ such that, if $|c_1-c_2|\leq \delta$, then $\Vert\eta_{c_1}-\eta_{c_2}\Vert_{H^3}\leq \varepsilon$. Using that we can control the third derivative with respect to $c$ in Lemma~\ref{lemme: exponential decay}, we deal the same way for the remaining derivatives $\partial_c^j\partial_x^k\eta_{c_1}$ for $j\leq 2$ and $k\leq 3$. The fact that $c\mapsto v_c\in\mathcal{C}^2\left((c_0,c_s),H^2(\R)\right)$ then follows from~\eqref{expression de v_c en fonction de eta_c}.

Finally, we prove that \begin{equation}\label{équivalent de p'(c) quand c proche c_s}
    \dfrac{d}{dc}\big(p(Q_c)\big) \underset{c\rightarrow c_s}{\sim}\dfrac{-18c_s^2\nu_c}{k^2}.
\end{equation}

As a consequence of~\eqref{equation -eta'^2=nu(eta)} and Proposition~\ref{prop: existence d'une solution eta qui tend vers 0 si il existe un zero z_c}, and doing the substitution $\xi=\eta_c(0)$ in the definition of the momentum~\eqref{expression du moment}, we get the following formula, and the following asymptotics already known (we refer to~\cite{Chiron7}):

\begin{equation}\label{expression de p(Q_c) avec xi_c le zéro}
        p(Q_c)=\dfrac{c}{2}\int_0^{\xi_c}\dfrac{\xi^2}{(1-\xi)\sqrt{-\mathcal{N}_c(\xi)}}d\xi=\mathcal{O}(\nu_c^3).
\end{equation}

Now we write
\begin{align}
    p(Q_c)=\underbrace{\int_0^{\xi_c}(g-\widetilde{g})(\xi,c)d\xi}_{I_c}+ \underbrace{\int_0^{\xi_c}\widetilde{g}(\xi,c)d\xi}_{II_c}\label{décomposition de p(gv_c) en g-gsing +gsing},
\end{align}

where \begin{equation*}
    g(\xi,c)=\dfrac{c}{2}\dfrac{\xi^2}{(1-\xi)\sqrt{-\mathcal{N}_c(\xi)}}\quad\text{and}\quad \widetilde{g}(\xi,c)=\dfrac{c}{2}\dfrac{\xi_c^2}{(1-\xi_c)\sqrt{(\xi_c-\xi)\mathcal{N}_c'(\xi_c)}}.
\end{equation*}

Using~\eqref{développement asymptotique de xi_c en nu_c}, we state a few asymptotics that will help us later. Recall that $\widetilde{k}:=-\frac{3}{k}$, we have
\begin{equation}\label{équivalent xi_c}
    \xi_c \underset{c\rightarrow c_s}{=}\widetilde{k}\nu_c^2 + \mathcal{O}(\nu_c^3),
\end{equation}

\begin{equation}
    \mathcal{N}_c'(\xi_c)\underset{c\rightarrow c_s}{=}\widetilde{k}\nu_c^4 + \mathcal{O}(\nu_c^6),
\end{equation}

\begin{equation}
\mathcal{N}_c''(\xi_c)\underset{c\rightarrow c_s}{=}4\nu_c^2 + \mathcal{O}(\nu_c^4).
\end{equation}

On the one hand, we compute

\begin{equation}\label{calcul et équivalent II_c}
    II_c=\int_0^{\xi_c}\widetilde{g}(\xi,c)d\xi =\dfrac{c\xi_c^{\frac{5}{2}}}{(1-\xi_c)\sqrt{\mathcal{N}'_c(\xi_c)}}\underset{c\rightarrow c_s}{\sim}c_s\widetilde{k}^2\nu_c^3,
\end{equation}

so that

\begin{align*}
    \dfrac{d}{dc}\big( II_c\big)  = \dfrac{1}{c}II_c + &\dfrac{\frac{5}{2}c\xi_c^{\frac{3}{2}}\partial_c\xi_c}{(1-\xi_c)\sqrt{\mathcal{N}'_c(\xi_c)}}+ \dfrac{c\xi_c^{\frac{5}{2}}\partial_c\xi_c}{(1-\xi_c)^2 \sqrt{\mathcal{N}'_c(\xi_c)}}\\
    &-\dfrac{2c^2\xi_c^{\frac{7}{2}}}{(1-\xi_c)\big(\mathcal{N}'_c(\xi_c)\big)^{\frac{3}{2}}}-\dfrac{c\xi_c^{\frac{5}{2}}\partial_c\xi_c\mathcal{N}_c''(\xi_c)}{2(1-\xi_c)\big(\mathcal{N}'_c(\xi_c)\big)^{\frac{3}{2}}}.
\end{align*}

We recall that, by Claim~\ref{claim: xi_c lesssim varepsilon^2},
\begin{equation*}
    \partial_c\xi_c =-\dfrac{2c\xi_c^2}{\mathcal{N}_c'(\xi_c)},
\end{equation*}

therefore we obtain the asymptotics

\begin{equation}\label{équivalent d/dc II_c}
    \dfrac{d}{dc}\big( II_c\big)\underset{\xi\rightarrow \xi_c}{=}-3c_s^2\widetilde{k}^2\nu_c +\grandOde{\nu_c^3}.
\end{equation}

On the other hand, we notice \begin{equation}
    g(\xi,c)-\widetilde{g}(\xi,c)\underset{\xi\rightarrow \xi_c}{\sim}\dfrac{c\xi_c\sqrt{\xi_c-\xi}(\xi_c-2)}{2\sqrt{\mathcal{N}'_c(\xi_c)}},
\end{equation}

so that the function $(g-\widetilde{g})(.,c)$ extends to $\xi_c$ by continuity. We infer, that 

\begin{align*}
    \dfrac{d}{dc}\big(I_c\big) & = \partial_c\xi_c \lim_{\xi\rightarrow\xi_c}(g-\widetilde{g})(\xi,c)+\int_0^{\xi_c}\partial_c(g-\widetilde{g})(\xi,c)d\xi\\
    & = 0+ \underbrace{\int_0^{\xi_c} l(\xi,c)d\xi}_{p_1(c)} + \underbrace{\int_0^{\xi_c} c\partial_c l(\xi,c)d\xi}_{p_2(c)},
\end{align*}

with \begin{equation*}
    l(\xi,c)=\dfrac{1}{c}\big(g(\xi,c)-g(\xi_c,c)\big)=\dfrac{1}{2}\Bigg(\dfrac{\xi^2}{(1-\xi)\sqrt{-\mathcal{N}_c(\xi)}}-\dfrac{\xi_c^2}{(1-\xi_c)\sqrt{(\xi_c-\xi)\mathcal{N}'_c(\xi_c)}}\Bigg).
\end{equation*}

We have, by~\eqref{expression de p(Q_c) avec xi_c le zéro} and~\eqref{calcul et équivalent II_c},
\begin{equation*}
    p_1(c)=\dfrac{p(Q_c)}{c}-\dfrac{II_c}{c}=\mathcal{O}(\nu_c^3),
\end{equation*}

and we derive, using the expression of $\partial_c\xi_c$, that
\begin{align*}
    &p_2 (c)  =\underbrace{\dfrac{c^2}{2}\int_0^{\xi_c}\Bigg(\dfrac{\xi^4}{(1-\xi)\big(-\mathcal{N}_c(\xi)\big)^\frac{3}{2}}-\dfrac{\xi_c^4}{(1-\xi_c)\big((\xi_c-\xi)\mathcal{N}'_c(\xi_c)\big)^\frac{3}{2}}\Bigg)d\xi }_{p_2^1(c)}\\
    & +\underbrace{c^2\bigg(\dfrac{\xi_c^3}{(1-\xi_c)\big(\mathcal{N}'_c(\xi_c)\big)^\frac{3}{2}}-\dfrac{\xi_c^4\mathcal{N}''_c(\xi_c)}{2(1-\xi_c)\big(\mathcal{N}'_c(\xi_c)\big)^\frac{5}{2}}+\dfrac{\xi_c^3}{\big(\mathcal{N}'_c(\xi_c)\big)^\frac{3}{2}}\Big(\dfrac{2}{1-\xi_c}+\frac{\xi_c}{(1-\xi_c)^2}\Big)\bigg)\int_0^{\xi_c}\dfrac{d\xi}{\sqrt{\xi_c-\xi}}}_{p_2^2(c)}.\\
\end{align*}

The second term reduces to \begin{equation}\label{équivalent de p_2^2}
p_2^2(c)=2c_s^2\widetilde{k}^2 \nu_c +\mathcal{O}(\nu_c^3).
\end{equation}

On the other hand, we decompose the first one as

\begin{align}
    p_2^1(c)&=\dfrac{c^2}{2}\int_0^{\xi_c}\Bigg(\dfrac{\xi^4}{(1-\xi)\big(-\mathcal{N}_c(\xi)\big)^\frac{3}{2}}-\dfrac{\xi^4}{(1-\xi_c)\big(-\mathcal{N}_c(\xi)\big)^\frac{3}{2}}\Bigg)d\xi \notag\\
    & + \dfrac{c^2}{2}\int_0^{\xi_c}\Bigg(\dfrac{\xi^4}{(1-\xi_c)\big(-\mathcal{N}_c(\xi)\big)^\frac{3}{2}}-\dfrac{\xi^4}{(1-\xi_c)\big((\xi_c-\xi)\big(\mathcal{N}'_c(\xi_c)\big)^\frac{3}{2}}\Bigg)d\xi\notag \\
    & +\dfrac{c^2}{2}\int_0^{\xi_c}\Bigg(\dfrac{\xi^4}{(1-\xi_c)\big((\xi_c-\xi)\mathcal{N}'_c(\xi_c)\big)^\frac{3}{2}}-\dfrac{\xi_c^4}{(1-\xi_c)\big((\xi_c-\xi)\mathcal{N}'_c(\xi_c)\big)^\frac{3}{2}}\Bigg)d\xi \notag\\
    & =-c_s^2\widetilde{k}^2 \nu_c +\mathcal{O}(\nu_c^3)\label{équivalent de p_2^1}.
\end{align}

Now adding~\eqref{équivalent d/dc II_c},~\eqref{équivalent de p_2^2} and~\eqref{équivalent de p_2^1} , we obtain, according to decomposition~\eqref{décomposition de p(gv_c) en g-gsing +gsing},

\begin{equation}
    \dfrac{d}{dc}\big(p(\gv_c)\big) \underset{c\rightarrow c_s}{\sim}-2c_s^2\nu_c\widetilde{k}^2.
\end{equation}

To conclude there exists a threshold $c_0$ such that for $c\in (c_0,c_s)$, condition~\eqref{théorème: hydrodynamique condition de grillakis sans le signe} is fulfilled.

\section{Perturbation of a sum of solitons and minimizing properties}\label{Perturbation of a sum of solitons and minimizing properties}
Before we get to the proof of Proposition~\ref{proposition: développement de E en Q=R_c,a + varepsilon}, we need to settle some tools that will simplify the computations. For $M\geq 2$, we define the following linear subspace of the set of multivariate real polynomials:

\begin{align*}
    \mathcal{P}_M=\mathcal{P}[X_1,...,X_M]:=\Big\{ \sum_{|\alpha|\leq m} p_\alpha X^\alpha \Big| m\in\N, p_{ke_i}=0\ ,\forall (k,i)\in \{0,...,m\}\times\{1,...,M\}\Big\},
\end{align*}

where $(e_1,...,e_M)$ is the canonical base of $\R^M$. This vector space is constructed such that the monomials are not contained in it and is also stable under the usual multiplication. One shall list a few of these remarkable polynomials. For instance, the elementary symmetric polynomials are contained in $\mathcal{P}_M$, and will be labelled

\begin{equation}
    S^{k,M}:=\sum_{1\leq j_1 <...<j_k\leq M} X_{j_1}...X_{j_k},
\end{equation}

for $k\in\{ 1,...,M\}$. To simplify, we will denote $P(\eta_{\gc,\ga}):=P(\eta_{c_1,a_1},...,\eta_{c_N,a_N})$ and $P(v_{\gc,\ga}):=P(v_{c_1,a_1},...,v_{c_N,a_N})$, for $P\in \mathcal{P}[X_1,...,X_N]$ and $P(Q_{\gc,\ga}):=P(\eta_{c_1,a_1},...,\eta_{c_N,a_N},v_{c_1,a_1},...,v_{c_N,a_N})$ for $P\in\mathcal{P}[X_1,...,X_N,Y_1,...,Y_N]$. Furthermore, for polynomials $P\in\mathcal{P}[X_{\sigma(1)},...,X_{\sigma(M)}]$ with $M<N$ and $\sigma\in\mathfrak{S}_M$, we will still use the notation $P(\eta_{\gc,\ga})$ to designate $$P(\eta_{\gc,\ga})=P(\eta_{c_{\sigma(1)},a_{\sigma(1)}},...,\eta_{c_{\sigma(M)},a_{\sigma(M)}}).$$

For $M\geq 2$, we can also define

\begin{equation*}
    B^M:=\sum_{k=1}^M\sum_{\substack{j=1 \\ j\neq k}}^M X^2_k X_j\in\mathcal{P}_M\quad\text{and}\quad \widetilde{B}^{2M}:=\sum_{k=1}^M\sum_{\substack{j=1 \\ j\neq k}}^M Y^2_k X_j\in\mathcal{P}_{2M}
\end{equation*} 

and
\begin{equation*}
    C^M:=\sum_{k=1}^M\sum_{\substack{j=1 \\ j\neq k}}^M X^3_k X_j\in\mathcal{P}_M,
\end{equation*}

and even \begin{equation}
    D^M:= \prod_{k=1}^M (1-X_k)-(1-\sum_{k=1}^M X_k) =\sum_{k=2}^{M} (-1)^{k} S^{k,M}\in\mathcal{P}_M,
\end{equation}

or \begin{equation}
    D_k^{M-1}:= \prod_{\substack{j=1 \\ j\neq k}}^M (1-X_j)-(1-\sum_{\substack{j=1 \\ j\neq k}}^M X_j)\in\mathcal{P}_{M-1}.
\end{equation}

Now we introduce a lemma (proven in Section~\ref{section: several polynomials}) which provides a control of the $L^p$-norms of such polynomials, when they are evaluated at functions that decay exponentially.

\begin{lem}\label{lemme: intégrale d'un polynome couplé est un O(Lexp(-Lnu)}
For any $M\geq 2$, $\ga\in\pos_M(L),\gb\in (\R_+^*)^M$, $P\in \mathcal{P}[X_1,...,X_M]$, and functions $(f_k)_{k\in\{1,...,M\}}$ such that 

\begin{equation}\label{f_k(x)=mathcal{O}(e^{-b_k |x|})}
    f_k(x)=\grandOde{e^{-b_k |x|}},
\end{equation}

we have for any $p\in [1,+\ii]$,

\begin{equation}\label{Vert Pbig(tau_{a_1}f_1,...,tau_{a_M}f_Mbig)Vert_{L^p}=mathcal{O}Big( L^{frac{1}{p}}e^{-min_k (b_k)L}Big)}
    \left\Vert P\big(\tau_{a_1}f_1,...,\tau_{a_M}f_M\big)\right\Vert_{L^p}=\mathcal{O}\bigg( \Big(\frac{2}{p\min_{k} (b_k)}+L\Big)^{\frac{1}{p}}e^{-\min_k (b_k)L}\bigg),
\end{equation}

where $\tau_{a_k}f_k:=f_k(.-a_k)$.

\end{lem}

\begin{rem}\label{remarque: on ne distingue pas l'évaluation des polynomes couplées en fonction de eta ou v}
    An important illustration of this lemma is to replace the functions $f_k$ by any derivatives or powers of $\eta_{c_k,a_k}$ or $v_{c_k,a_k}$ which satisfy~\eqref{f_k(x)=mathcal{O}(e^{-b_k |x|})} because of their decay.
\end{rem}

\begin{proof}[Proof of Proposition~\ref{proposition: développement de E en Q=R_c,a + varepsilon}]
We write
\begin{equation*}
    E(Q)=E(R_{\gc,\ga}+\varepsilon)=E(R_{\gc,\ga})+
    \nabla E(R_{\gc,\ga}).\varepsilon+\dfrac{\nabla^2 E(R_{\gc,\ga})(\varepsilon,\varepsilon)}{2}+\mathcal{R}_{\gc,\ga}(\varepsilon),
\end{equation*}

where the remainder $\mathcal{R}_{\gc,\ga}(\varepsilon)$ can be explicited as follows, \begin{equation*}\label{mathcal{R}(eta,v)=mathcal{O}(VertvarepsilonVert_X^3)}
    \mathcal{R}_{\gc,\ga}(\varepsilon)=\int_0^1 \dfrac{(1-t)^2}{2}\nabla^3 E( \Rca+t\varepsilon)(\varepsilon,\varepsilon,\varepsilon)dt.
\end{equation*}

We also compute the higher order derivatives
\begin{equation}\label{expression de nabla E(R_{gc,ga}).varepsilon}
    \nabla E(R_{\gc,\ga}).\varepsilon =\int_\R \bigg(\dfrac{\partial_x\eta_{\gc,\ga}\partial_x\varepsilon_\eta}{4(1-\eta_{\gc,\ga})}+\dfrac{(\partial_x\eta_{\gc,\ga})^2\varepsilon_\eta}{8(1-\eta_{\gc,\ga})}+(1-\eta_{\gc,\ga})v_{\gc,\ga}\varepsilon_v+\dfrac{1}{2}\big(f(1-\eta_{\gc,\ga})-v_{\gc,\ga}^2\big)\varepsilon_\eta\bigg),
\end{equation}
and 
\begin{align}\label{expression de nabla^2 E(R_{gc,ga}).varepsilon}
    \nabla^2 E(R_{\gc,\ga})(\varepsilon,\varepsilon)=&\int_\R \bigg( \dfrac{(\partial_x \eta_{\gc,\ga})^2\varepsilon_\eta^2}{4(1-\eta_{\gc,\ga})^3}+\dfrac{\partial_x\eta_{\gc,\ga}\varepsilon_\eta\partial_x\varepsilon_\eta}{2(1-\eta_{\gc,\ga})^2}+\dfrac{(\partial_x\varepsilon_\eta)^2}{4(1-\eta_{\gc,\ga})}\notag\\
    & +(1-\eta_{\gc,\ga})\varepsilon_v^2-2v_{\gc,\ga}\varepsilon_\eta\varepsilon_v-\dfrac{f'(1-\eta_{\gc,\ga})}{2}\varepsilon_\eta^2 \bigg).
\end{align}

\begin{claim}\label{claim: E(R_{gc,ga})=sum_{k=1}^N E(Q_{c_k})+mathcal{O}Big(Le^{frac{nu_{gc_*}L}{2}}Big)}
    \begin{equation*}
    E(R_{\gc,\ga})=\sum_{k=1}^N E(Q_{c_k})+\mathcal{O}\Big(\Lambda(L,\gc)e^{-a_d\nu_\gc L}\Big).
\end{equation*}
\end{claim}

\begin{proof}

Recall that
    \begin{equation}\label{E(R_gc,ga)=dfracint_dfracpartial_xeta_gc,ga)^2}
        E(R_{\gc,\ga})=\dfrac{1}{8}\int_\R\dfrac{(\partial_x\eta_{\gc,\ga})^2}{1-\eta_{\gc,\ga}}+\dfrac{1}{2}\int_\R (1-\eta_{\gc,\ga})v_{\gc,\ga}^2+\dfrac{1}{2}\int_\R F(1-\eta_{\gc,\ga}).
    \end{equation}

First, we study the kinetic energy, namely, the first and second terms in~\eqref{E(R_gc,ga)=dfracint_dfracpartial_xeta_gc,ga)^2}, by writing 
\begin{equation}\label{(1-eta_gc,ga)v_c,ga^2}
    (1-\eta_{\gc,\ga})v_{\gc,\ga}^2=\sum_{k=1}^N (1-\eta_{c_k,a_k})v_{c_k,a_k}^2 -B^{2N}(Q_{\gc,\ga})+2S^{2,N}(v_{\gc,\ga})-2S^{2,N}(v_{\gc,\ga})S^{1,N}(\eta_{\gc,\ga}).
\end{equation}

Using in addition Proposition~\ref{proposition: min |gv_c|=delta c/c_s}, we deal with the first term of $E(\Rca)$ by writing \begin{align}
    \dfrac{(\partial_x \eta_{\gc,\ga})^2}{1-\eta_{\gc,\ga}}&=\sum_{k=1}^N \dfrac{(\partial_x \eta_{c_k,a_k})^2}{(1-\eta_{c_k,a_k})\Big(1-\sum_{\substack{l=1\\ l\neq k}}^N\frac{\eta_{c_l,a_l}}{1-\eta_{c_k,a_k}}\Big)} + \dfrac{2}{1-\eta_{\gc,\ga}}S^{2,N}(\partial_x\eta_{\gc,\ga})\notag\\
    &=\sum_{k=1}^N \dfrac{(\partial_x \eta_{c_k,a_k})^2}{1-\eta_{c_k,a_k}} + \mathcal{O}\big(\widetilde{B}^{2N}(\eta_{\gc,\ga},\partial_x\eta_{\gc,\ga})\big)+\mathcal{O}\big(S^{2,N}(\partial_x\eta_{\gc,\ga})\big)\label{sum_k=1^N dfrac(partial_x eta_c_k,a_k)^2}.
\end{align}

As a consequence of Lemma~\ref{lemme: F,f, f'(1-eta_c)} in the appendix, we obtain
\begin{equation}\label{F(1-eta_{gc,ga})=sum_{k=1}^N F(1-eta_{c_k,a_k})+ mathcal{O}ig((S^{2,N}+B^N+C^Netagc,a)big)}
    F(1-\eta_{\gc,\ga})=\sum_{k=1}^N  F(1-\eta_{c_k,a_k})+ \mathcal{O}\big(S^{2,N}(\eta_{\gc,\ga})\big)+\mathcal{O}\big(S^{3,N}(\eta_{\gc,\ga})\big) +\mathcal{O}\big(B^N(\eta_{\gc,\ga})\big)+\mathcal{O}\big(C^N(\eta_{\gc,\ga})\big).
\end{equation}

We integrate~\eqref{(1-eta_gc,ga)v_c,ga^2},~\eqref{sum_k=1^N dfrac(partial_x eta_c_k,a_k)^2} and~\eqref{F(1-eta_{gc,ga})=sum_{k=1}^N F(1-eta_{c_k,a_k})+ mathcal{O}ig((S^{2,N}+B^N+C^Netagc,a)big)} on $\R$. Therefore, by exponential decay~\eqref{estimée décroissance exponentielle à tout ordre pour eta_c et v_c} and using Lemma~\ref{lemme: intégrale d'un polynome couplé est un O(Lexp(-Lnu)} with $p=1$, we conclude the proof of the claim.
\end{proof}

Now we deal with $\nabla E(\Rca).\varepsilon$. As a matter of example, we only handle the term associated with the nonlinearity $f$ and we refer to~\cite{BetGrSm1} concerning the other terms. By Lemma~\ref{lemme: F,f, f'(1-eta_c)}, we have 
\begin{equation*}
    f(1-\eta_{\gc,\ga})\varepsilon_\eta  = \sum_{k=1}^N f(1-\eta_{c_k,a_k}) \varepsilon_\eta 
    +\mathcal{O}\Big(\varepsilon_\eta B^N(\eta_{\gc,\ga})\Big)+\mathcal{O}\Big(\varepsilon_\eta S^{2,N}(\eta_{\gc,\ga})\Big).
\end{equation*}

Integrating the previous equation on $\R$, using the Cauchy-Schwarz inequality and eventually Lemma~\ref{lemme: intégrale d'un polynome couplé est un O(Lexp(-Lnu)} with $p=2$, leads to
\begin{align*}
    \int_\R f(1-\eta_{\gc,\ga})\varepsilon_\eta & = \sum_{k=1}^N \int_\R f(1-\eta_{c_k,a_k}) \varepsilon_\eta 
    +\mathcal{O}\Big(\normLdeux{\varepsilon_\eta}\normLdeux{S^{2,N}(\eta_{\gc,\ga})}  \Big)+\mathcal{O}\Big(\normLdeux{\varepsilon_\eta}\normLdeux{B^N(\eta_{\gc,\ga})}  \Big)\\
    & = \sum_{k=1}^N \int_\R f(1-\eta_{c_k,a_k}) \varepsilon_\eta+\mathcal{O}\Big(\Vert\varepsilon\Vert_\mathcal{X} \Lambda(L,\gc)^{\frac{1}{2}}e^{-a_d\nu_\gc L}\Big).
\end{align*}

The other terms can be dealt with similarly, so that we infer
\begin{equation*}
    \nabla E(R_{\gc,\ga}).\varepsilon = \sum_{k=1}^N \nabla E(Q_{c_k,a_k}).\varepsilon + \mathcal{O}\big(\Vert\varepsilon\Vert_\mathcal{X} \Lambda(L,\gc)^{\frac{1}{2}}e^{-a_d\nu_\gc L}\big).
\end{equation*}

Since $Q\in \Operpca$ and $Q_{c_k,a_k}$ solves~\eqref{nabla E(v)-cnabla p(v)=0} with $c=c_k$, we have for any $k\in\{1,...,N\}$,

\begin{equation*}
    \nabla E(Q_{c_k,a_k}).\varepsilon =\big(\nabla E(Q_{c_k,a_k})-c_k\nabla p(Q_{c_k,a_k})\big).\varepsilon =0,
\end{equation*}

so that 
\begin{equation}\label{nabla E(R_{gc,ga}).varepsilon = mathcal{O}big(VertvarepsilonVert_X L^{rac{1}{2}}e^{-frac{nu_{gc^*}L}{2}}big)}
    \nabla E(R_{\gc,\ga}).\varepsilon = \mathcal{O}\big(\Vert\varepsilon\Vert_\mathcal{X} \Lambda(L,\gc)^{\frac{1}{2}}e^{-a_d\nu_\gc L}\big).
\end{equation}

Now we handle the quadratic term, we first define

\begin{equation*}\label{nabla^2 E(R_gc,ga).(varepsilon,varepsilon)=:int_R J(R_gc,ga,varepsilon)}
    \nabla^2 E(R_{\gc,\ga}).(\varepsilon,\varepsilon)=:\int_\R J(R_{\gc,\ga},\varepsilon),
\end{equation*}

and we decompose according to the partition of the unity~\eqref{partition de l'unité Phi_k},

\begin{equation*}
    \int_\R J(R_{\gc,\ga},\varepsilon)=\sum_{k=1}^N \int_\R J(R_{\gc,\ga},\varepsilon)\Phi_k + \sum_{k=0}^N \int_\R J(R_{\gc,\ga},\varepsilon)\Phi_{k,k+1}.
\end{equation*}

We write 

\begin{claim}\label{claim: int_R J(R_gc,ga,varepsilon)Phi_k =}
    For $\tau$ such that $2\tau<\nu_\gc$, we have

\begin{equation*}
    \int_\R J(R_{\gc,\ga},\varepsilon)\Phi_k = \int_\R J(Q_{c_k,a_k},\varepsilon)\Phi_k + \mathcal{O}\Big( \normX{\varepsilon}^2 e^{-\frac{a_d\tau L}{2}}\Big)
\end{equation*}

and 

\begin{equation*}
    \int_\R J(R_{\gc,\ga},\varepsilon)\Phi_{k,k+1}=\int_\R J(0,\varepsilon)\Phi_{k,k+1}+\mathcal{O}\Big(\normX{\varepsilon}^2 e^{-\frac{a_d\tau L}{2}} \Big).
\end{equation*}

\end{claim}

\begin{proof}
    Like we did for $\nabla E(\Rca).\varepsilon$, we only deal with the term in which the nonlinearity $f$ intervenes. In view of Lemma~\ref{lemme: F,f, f'(1-eta_c)}, we write

\begin{equation*}
        f'(1-\eta_{\gc,\ga})\varepsilon_\eta^2\Phi_k = f'(1-\eta_{c_k,a_k})\varepsilon_\eta^2\Phi_k + \mathcal{O}\Big(A_k^N\big(\eta_{c_1,a_1},...,\eta_{c_{k-1},a_{k-1}},\Phi_k,\eta_{c_{k+1},a_{k+1}},...,\eta_{c_{N},a_{N}}\big)\varepsilon_\eta^2\Big),
\end{equation*}

where $G^{N}_k=X_k\sum_{\substack{j=1 \\ j\neq k}}^N X_j\in\mathcal{P}[X_1,...,X_N]$. Since $\Phi_k(x)=\mathcal{O}\big(e^{-2\tau|x-a_k-\frac{L}{4}|}\big)$ and $\eta_{c_j,a_j}(x)=\mathcal{O}\big(e^{-\nu_\gc|x-a_j|}\big)$ for any $j\neq k$, and since $(a_1,...,a_{k-1},a_k+L/4,a_{k+1},...,a_N)\in\pos_N\big(\frac{3L}{4}\big)$, we can apply Lemma~\ref{lemme: intégrale d'un polynome couplé est un O(Lexp(-Lnu)} with $P=A_k^N$ and $p=+\ii$, and we deduce

\begin{equation}\label{int_R f'(1-eta_gc,gavarepsilon_eta^2Phi_k}
        \int_\R f'(1-\eta_{\gc,\ga})\varepsilon_\eta^2\Phi_k= \int_\R f'(1-\eta_{c_k,a_k})\varepsilon_\eta^2\Phi_k + \mathcal{O}\Big(\normX{\varepsilon}^2 e^{-\frac{3a_d\min(2\tau,\nu_\gc)L}{4}}\Big).
    \end{equation}

In a similar manner, we write the following Taylor expansion of the first order,

\begin{align*}
    f'(1-\eta_{\gc,\ga})\varepsilon_\eta^2\Phi_{k,k+1} &=f'(1)\varepsilon_\eta^2\Phi_{k,k+1}\\
    &-\int_0^1\varepsilon_\eta^2\Big(\eta_{c_k,a_k}\Phi_{k,k+1} + A_k^N(\eta_{c_1,a_1},...,\Phi_{k,k+1},...,\eta_{c_N,a_N})\Big) f''(1-t\eta_{\gc,\ga})dt.
\end{align*}

We have for $\widetilde{k}\leq k,\eta_{c_{\widetilde{k}},a_{\widetilde{k}}}(x)=\mathcal{O}\big(e^{-\nu_\gc(x-a_{\widetilde{k}})^+}\big)$ and $\Phi_{k,k+1}(x)=\mathcal{O}\big(e^{-2\tau(x-a_k-\frac{L}{4})^-}\big)$, and for $\widetilde{k}\geq k+1,\eta_{c_{\widetilde{k}},a_{\widetilde{k}}}(x)=\mathcal{O}\big(e^{-\nu_\gc(x-a_{\widetilde{k}})^-}\big)$ and $\Phi_{k,k+1}(x)=\mathcal{O}\big(e^{-2\tau(x-a_{k+1}+\frac{L}{4})^+}\big)$. Now using Lemma~\ref{lemme: intégrale d'un polynome couplé est un O(Lexp(-Lnu)} with $P=S^{2,2}\in\mathcal{P}[X_1,X_2]$ and $p=+\ii$, we get

\begin{equation}\label{int_R varepsilon_eta^2eta_c_k,a_kPhi_k,k+1}
    \int_\R \varepsilon_\eta^2\eta_{c_k,a_k}\Phi_{k,k+1}=\mathcal{O}\big(\normLdeux{\varepsilon_\eta}^2\normLii{\eta_{c_k,a_k}\Phi_{k,k+1}}\big)=\mathcal{O}\Big(\normX{\varepsilon}^2 e^{-\frac{a_d\min(2\tau,\nu_\gc)L}{4}}\Big).
\end{equation}

On the other hand, we also use Lemma~\ref{lemme: intégrale d'un polynome couplé est un O(Lexp(-Lnu)} with $P=A_k^N$ and $p=+\ii$, so that 

\begin{equation}\label{int_R varepsilon_eta^2 A_k^N(eta_c_1,a_1,...,Phi_k,k+1,...,eta_c_N,a_N}
    \int_\R \varepsilon_\eta^2 A_k^N(\eta_{c_1,a_1},...,\Phi_{k,k+1},...,\eta_{c_N,a_N})=\mathcal{O}\Big(\normX{\varepsilon}^2 e^{-\frac{a_d\min(2\tau,\nu_\gc)L}{4}}\Big).
\end{equation}

As a consequence of both~\eqref{int_R varepsilon_eta^2eta_c_k,a_kPhi_k,k+1} and~\eqref{int_R varepsilon_eta^2 A_k^N(eta_c_1,a_1,...,Phi_k,k+1,...,eta_c_N,a_N}, we derive

\begin{equation}\label{int_R f'(1-eta_c,ga)var}
    \int_\R f'(1-\eta_{\gc,\ga})\varepsilon_\eta^2\Phi_{k,k+1} =\int_\R f'(1)\varepsilon_\eta^2\Phi_{k,k+1} + \mathcal{O}\Big(\normX{\varepsilon}^2 e^{-\frac{a_d\min(2\tau,\nu_\gc)L}{4}}\Big).
\end{equation}

We use the same type of consideration for the other terms, then by hypothesis on $\tau$, and combining~\eqref{int_R f'(1-eta_gc,gavarepsilon_eta^2Phi_k} and~\eqref{int_R f'(1-eta_c,ga)var}, it eventually proves Claim~\ref{claim: int_R J(R_gc,ga,varepsilon)Phi_k =}:

\begin{equation*}
    \int_\R J(R_{\gc,\ga},\varepsilon)\Phi_{k,k+1}=\int_\R J(0,\varepsilon)\Phi_{k,k+1}+\mathcal{O}\Big(\normX{\varepsilon}^2 e^{-\frac{a_d\tau L}{2}} \Big).
\end{equation*}
\end{proof}

The rest of the argument does not depend on the nonlinearity $f$, but only on the construction of the partition of the unity~\eqref{partition de l'unité Phi_k}, so we just take $\tau$ like in Claim~\ref{claim: int_R J(R_gc,ga,varepsilon)Phi_k =} and deduce that
\begin{equation}\label{nabla^2 E(R_{gc,ga}).(varepsilon,varepsilon)=nabla^2 E(Q_{c_k}).(varepsilon_k,varepsilon_k)}
    \nabla^2 E(R_{\gc,\ga}).(\varepsilon,\varepsilon)=\nabla^2 E(Q_{c_k}).(\varepsilon_k,\varepsilon_k)+\nabla^2 E(0).(\varepsilon_{k,k+1},\varepsilon_{k,k+1}) + \mathcal{O}(\tau\normX{\varepsilon}^2)+\mathcal{O}\Big(\normX{\varepsilon}^2 e^{-\frac{a_d\tau L}{2}} \Big).
\end{equation}

Now summing the estimates on each order of the derivative of the energy, we conclude the proof of Proposition~\ref{proposition: développement de E en Q=R_c,a + varepsilon}.
\end{proof}

Contrary to the previous proofs in this section, the proofs of Proposition~\ref{proposition: développement de p_k en Q=R_c,a + varepsilon} does not depend on the shape of the nonlinearity $f$ and can be written with the methods from~\cite{BetGrSm1}.

\section{Orthogonal decomposition and control on the functional $G$}\label{section: orthogonal decomposition and control on the functional }
The proof of Proposition~\ref{proposition: décomposition orthogonale sans dynamique} follows the lines of the proof of the analogous Proposition~2 in~\cite{BetGrSm1}, provided the nonlinearity $f$ satisfies the suitable properties. The only crucial property to check is the exponential decay of the travelling waves, and of its first spatial derivatives, which is provided by Theorem~\ref{théorème: il existe une branche C^1 de solitons proche de c_s}. We refer to~\cite{BetGrSm1} for more details. Throughout this section and to suit the framework of Subsection~\ref{subsection: almost minimizing property d'une somme de solitons}, we shall use the notation $(\gc,\ga):=\big(\mathfrak{C}(Q),\mathfrak{A}(Q)\big)$ with $\mathfrak{C},\mathfrak{A}$ and $\varepsilon$ the functions given by Proposition~\ref{proposition: décomposition orthogonale sans dynamique}.

\begin{proof}[Proof of Corollary~\ref{corollaire: décomposition orthogonale}]
As for Proposition~\ref{proposition: décomposition orthogonale sans dynamique}, there are similarities with~\cite{BetGrSm1}, in particular concerning the proof of the first and second inequalities in~\eqref{nu_athfrak(eta,v)geq dfracnu_gc^*} for which we refer to~\cite{BetGrSm1}. Take $L\geq L_2$. Since $Q\in\mathcal{U}_{\gc^*}(\alpha_2,L)$, there exists $\ga^*\in\pos_N(L)$ such that $\normX{Q-R_{\gc^*,\ga^*}}<\alpha_2$. Taking $\alpha_2$ small enough,~\eqref{normX{varepsilon(eta,v)}+Vertathfrak{C}(ta,v)-gc^*Vert_1 + Vert athfrak{A}} and Lemma~\ref{lemme: norm Lii de eta_c^* inférieur à delta} yield to
\begin{equation*}
    1-\eta(x)\geq 1-\normLii{\eta_{\gc^*,\ga^*}}-\alpha_2\geq\dfrac{1-\beta^*}{2}.
\end{equation*}
This implies~\eqref{mathcalU_gc^*(alpha_1,L_1)subsetNenergysethydro} and~\eqref{1-etageq dfracmu_gc^8}. We deduce also~\eqref{mathfrakA(eta,v)inpos_N(L-1)} from taking $\alpha_2$ small enough in~\eqref{normX{varepsilon(eta,v)}+Vertathfrak{C}(ta,v)-gc^*Vert_1 + Vert athfrak{A}}.

Regarding $\kappa_\gc$, we first use that $c\mapsto Q_c$ is $\mathcal{C}^2$ by writing that, for any $k\in\{1,...,N\}$,
\begin{equation*}
    \left| \dfrac{d}{dc}\Big(p(Q_c)\Big)_{|c=c_k^*}-\dfrac{d}{dc}\Big(p(Q_c)\Big)_{|c=c_k}\right|\leq \left\Vert \dfrac{d^2}{dc^2}\Big(p(Q_c)\Big)\right\Vert_{L^\ii([c_k^*,c_k])}|c_k^*-c_k|,
\end{equation*}

thus \begin{equation*}
    -\dfrac{d}{dc}\Big(p(Q_c)\Big)_{|c=c_k}\geq \kappa_{\gc^*} -\left\Vert \dfrac{d^2}{dc^2}\Big(p(Q_c)\Big)\right\Vert_{L^\ii([c_k^*,c_k])}|c_k^*-c_k|,
\end{equation*}
so it remains to bound $\Vert \frac{d^2}{dc^2}\big(p(Q_c)\big)\Vert_{L^\ii([c_k^*,c_k])}$ by a constant that only depends on $\gc^*$. The exponential decay in Theorem~\ref{théorème: il existe une branche C^1 de solitons proche de c_s} and Lemma~\ref{lemme: estimation différentielle E et p} provide a positive constant $\widetilde{K}$ such that for any $c\in[c_k^*,c_k]$, \begin{align}
    \left|\dfrac{d^2}{dc^2}\Big(p(Q_c)\Big)\right|&\leq \left|\nabla p(Q_c).\partial_c^2 Q_c\right|+\left|\nabla^2 p(Q_c).(\partial_c Q_c,\partial_c Q_c)\right|\leq \widetilde{K}\label{borne de d/dc^2 p(Q_c)},
\end{align}

where $\widetilde{K}$ can be explicited in terms of $\gc^*,\nu_{\gc^*}$ and $\mu_{\gc^*}$, provided that $|\gc-\gc^*|\leq \frac{\mu_{\gc^*}}{2}$. We then deduce the third inequality in~\eqref{nu_athfrak(eta,v)geq dfracnu_gc^*} by taking $\alpha_2$ sufficiently small so that the previous condition is satisfied and such that $\kappa_{\gc^*}-\widetilde{K}K_1\alpha_2\geq \frac{\kappa_{\gc^*}}{2}$.\\
Now we prove the last inequality in~\eqref{nu_athfrak(eta,v)geq dfracnu_gc^*}. According to definition~\eqref{définition l_c}, take $\varepsilon\in\energysethydro$ satisfying~\eqref{proposition condition d'orthogonalité sur varepsilon} with $c=c_k$, and write its orthogonal decomposition
\begin{equation}
    \varepsilon = \lambda \partial_x Q_{c_k^*} +\mu \nabla p(Q_{c_k^*}) +\varepsilon^*,
\end{equation}

where $\lambda,\mu\in\R$, $\varepsilon^*$ satisfies the orthogonal conditions~\eqref{proposition condition d'orthogonalité sur varepsilon} with $c=c_k^*$, and where $\nabla p(Q_{c_k^*})$ is identified to its representing vector in $L^2(\R)\times L^2(\R)$. Since $\varepsilon$ satisfies~\eqref{proposition condition d'orthogonalité sur varepsilon}, we compute 
\begin{equation*}
     \left\{
\begin{array}{l}
    \psLdeuxLdeux{\varepsilon}{\partial_x Q_{c_k^*}-\partial_x Q_{c_k}}=\psLdeuxLdeux{\varepsilon}{\partial_x Q_{c_k^*}}=\lambda\normLdeux{\partial_x Q_{c_k^*}}^2, \\
    \psLdeuxLdeux{\varepsilon}{\nabla p(Q_{c_k^*})-\nabla p(Q_{c_k})}=\psLdeuxLdeux{\varepsilon}{\nabla p(Q_{c_k^*})}=\mu\normLdeux{\nabla p(Q_{c_k^*})}^2. \\
\end{array}
\right.
\end{equation*}

Since the quantities $\partial_x Q_{c_k^*}$ and $\nabla p(Q_{c_k^*})$ are different from zero, using Lemmas~\ref{lemme: c,a mapsto R_c,a est lipschitz continue} and~\ref{lemme: estimation différentielle E et p} for the momentum, there exists a constant $C_*>0$, only depending on $\gc^*$ such that 
\begin{equation}\label{|lambda|+|mu|leq C_*normXvarepsilonnormRgc^*-gc}
    |\lambda|+|\mu|\leq C_*\normX{\varepsilon}\normR{\gc^*-\gc}.
\end{equation}

On the other hand, noticing in addition that $\psLdeuxLdeux{\nabla p(Q_{c_k^*})}{\partial_x Q_{c_k^*}}=0$, we have 
\begin{align}
    H_{c_k}(\varepsilon)&=\lambda^2 H_{c_k}(\partial_x Q_{c_k^*})+\mu^2 H_{c_k}\big(\nabla p(Q_{c_k^*})\big)+H_{c_k}(\varepsilon^*) \notag \\ 
    &+2\lambda\psLdeuxLdeux{\mathcal{H}_{c_k}(\partial_x Q_{c_k^*})}{\varepsilon^*}+2\mu \psLdeuxLdeux{\mathcal{H}_{c_k}\big(\nabla p(Q_{c_k^*})\big)}{\varepsilon^*}\label{développement du calcul de H_c_k(varepsilon)}\\
    &+2\lambda\mu \psLdeuxLdeux{\mathcal{H}_{c_k}(\partial_x Q_{c_k^*})}{\nabla p(Q_{c_k^*})}\notag.
\end{align}

We notice, by definition of $\mathcal{H}_c$, that
\begin{equation}\label{mathcalH_c_k=(c_k^*-c_k)nabla^2 p}
    \mathcal{H}_{c_k}=(c_k^*-c_k)\nabla^2 p(Q_{c_k}) +\mathcal{H}_{c_k^*}.
\end{equation}

Taking benefit of~\eqref{mathcalH_c_k=(c_k^*-c_k)nabla^2 p}, there exists a positive constant $\widetilde{C}_*$ such that \begin{equation}\label{H_c(varepsilon)geq H_c_k^*arepsilon^*)-normRgc-gc^*}
    H_c(\varepsilon)\geq H_{c_k^*}(\varepsilon^*)-\normR{\gc-\gc^*}\normX{\varepsilon^*}^2 -\widetilde{C}_*\left( (\lambda^2+\mu^2)(\normR{\gc-\gc^*}+1)+(|\lambda|+|\mu|)\normX{\varepsilon^*} \right).
\end{equation}

Using the control~\eqref{|lambda|+|mu|leq C_*normXvarepsilonnormRgc^*-gc}, the orthogonality conditions on $\varepsilon^*$, the coercivity of $H_{c_k^*}$, and up to taking a larger constant $\widetilde{C}_*$, we obtain
\begin{equation*}
    H_c(\varepsilon)\geq \big(l_{c_k^*}-\normR{\gc-\gc^*}\big)\normX{\varepsilon^*}^2-\widetilde{C}_*\normR{\gc-\gc^*}\normX{\varepsilon}\big( \normR{\gc-\gc^*}\normX{\varepsilon} +\normX{\varepsilon^*} \big).
\end{equation*}

Moreover, we verify that $\normX{\varepsilon}^2=\grandOde{\lambda^2+\mu^2}+\grandOde{(|\lambda|+|\mu|)\normLdeux{\varepsilon^*}}+\normX{\varepsilon^*}^2$. Up to taking a smaller $\alpha_2$ so that $|\gc^*-\gc|$ is small enough, we infer from~\eqref{|lambda|+|mu|leq C_*normXvarepsilonnormRgc^*-gc} that $\normX{\varepsilon}=\grandOde{\normX{\varepsilon^*}}$ . Plugging this into~\eqref{H_c(varepsilon)geq H_c_k^*arepsilon^*)-normRgc-gc^*}, we deduce that there exists $l_*>0$, only depending on $\gc^*$ such that for any $\varepsilon\in\energysethydro$ satisfying~\eqref{proposition condition d'orthogonalité sur varepsilon} and any $k\in\{1,...,N\}$, \begin{equation*}
    H_{c_k}(\varepsilon)\geq  l_*\normX{\varepsilon}^2,
\end{equation*}
which leads to the last inequality in~\eqref{nu_athfrak(eta,v)geq dfracnu_gc^*}.
\end{proof}

We have imposed that $\tau_0<2\tau<\frac{\nu_{\gc^*}}{2}$. Before we pass to the proof of Corollary~\ref{corollaire: controle de la fonctionnelle G}, we combine Proposition~\ref{proposition: développement de E en Q=R_c,a + varepsilon} and~\ref{proposition: développement de p_k en Q=R_c,a + varepsilon} and deduce the following corollary concerning the functional $G$ defined in~\eqref{définition de la fonctionnelle G}.

\begin{cor}\label{corollaire: développement asymptotique de G en Q= R_c,a + varepsilon}
Let $L\geq L_2$. For $Q\in\mathcal{U}^\perp_{\gc,\ga}(L)$, we have
    \begin{align*}
        G(Q)=&\sum_{k=1}^N\big( E(Q_{c_k^*})-c_k^* p(Q_{c_k^*})\big) +\dfrac{1}{2}\bigg(\sum_{k=1}^N H_{c_k}(\varepsilon_k) + \sum_{k=0}^N H_0^k(\varepsilon_{k,k+1},\varepsilon_{k,k+1})\bigg)\\
        &+\mathcal{O}\big(\normR{ \gc-\gc^*}^2\big)+\mathcal{O}\Big(\normX {\varepsilon}^2\big(\tau+e^{-\frac{a_d\tau L}{2}}\big)\Big)+ \mathcal{R}_{\gc,\ga}(\varepsilon)+\grandOde{\normX{\varepsilon}^4}\\
        &+\mathcal{O}\Big(\Lambda(L,\gc)(e^{-a_d\nu_{\gc}L}+e^{-a_d\tau_0 L})\Big)+\mathcal{O}\big(\normX{\varepsilon} \Lambda(L,\gc)^{\frac{1}{2}}(e^{-a_d\nu_{\gc}L}+e^{-a_d\tau_0 L})\big),
    \end{align*}
where 
$$H_0^{k}(\varepsilon_{k,k+1})=\left\{
\begin{array}{l}
    \nabla^2 (E-c_1 p_1)(0).(\varepsilon_{0,1},\varepsilon_{0,1})\quad\text{if }k=0, \\
    \nabla^2(E-c_k p_k-c_{k+1}p_{k+1})(0).(\varepsilon_{k,k+1},\varepsilon_{k,k+1})\quad\text{if }k\in\{1,...,N-1\}, \\
    \nabla^2(E-c_N p_N)(0).(\varepsilon_{N,N+1},\varepsilon_{N,N+1})\quad\text{if }k=N.\\
\end{array}
\right.$$

\end{cor}

\begin{proof}
We just combine Proposition~\ref{proposition: développement de E en Q=R_c,a + varepsilon} and~\ref{proposition: développement de p_k en Q=R_c,a + varepsilon}. We first deal with the combination of the second order terms appearing in the expression of both previous propositions. We recognize the desired terms $H_{c_k}(\varepsilon_k)$ and $H_0^k(\varepsilon_{k,k+1})$ except that the speeds are $c_k^*$ instead of $c_k$. This can be overcome by using the same kind of decomposition than in~\eqref{mathcalH_c_k=(c_k^*-c_k)nabla^2 p}. Indeed, we write that $\nabla^2(E-c_k^*p)=\nabla^2(E-c_k p)+(c_k-c_k^*)\nabla^2 p$, therefore we deal with the terms involving $\nabla^2 p$ by using Lemma~\ref{lemme: estimation différentielle E et p}, by using that $\normX{\varepsilon_k}\lesssim \normX{\varepsilon}$ and by using the standard Young's inequality. This can be summarized in the following inequalities: $\big|(c_k-c_k^*)\nabla^2 p(Q_{c_k}).(\varepsilon_k,\varepsilon_k)\big|\lesssim |c_k-c_k^*|\normX{\varepsilon}^2\lesssim \normR{\gc-\gc^*}^2+\normX{\varepsilon}^4$. As for the terms involving $\nabla^2 p_k$, the same decomposition than previously and a straightforward computation of the second derivative of $p_k$ leads to the same control. It remains to study the term $\sum_{k=1}^N \big(E(Q_{c_k})-c_k^* p(Q_{c_k})\big)$. We use a second order Taylor expansion between $Q_{c_k}$ and $Q_{c_k^*}$. The first order term cancels because of~\eqref{nabla E(v)-cnabla p(v)=0}, and the second order term can be dealt by using Lemma~\ref{lemme: c,a mapsto R_c,a est lipschitz continue} and~\ref{lemme: estimation différentielle E et p}, hence the resulting $\sum_{k=1}^N \big(E(Q_{c_k^*})-c_k^* p(Q_{c_k^*})\big)+\grandOde{|\gc-\gc^*|^2}$. 
\end{proof}

Recall that $H_c$ is coercive whenever $\varepsilon$ satisfies some orthogonality conditions. We make now use of this property to state some almost coercivity property along the $\varepsilon_k$.

\begin{cor}\label{corollaire: almost coercivité le long des varepsilon_k}
    Let $L\geq L_2$. For $Q\in\mathcal{U}^\perp_{\gc,\ga}(L)$, there exists $\widetilde{l}_{\gc} >0$ proportional\footnote{As usual in this framework, the proportionality constant only depends on $\gc^*$.} to $\min(l_{\gc},\nu_{\gc})$ such that for any $k\in\{0,...,N\}$, 
\begin{equation}
    H_0^k(\varepsilon_{k,k+1})\geq \widetilde{l}_{\gc}\normX{\varepsilon_{k,k+1}}^2,
\end{equation}
and for any $k\in\{1,...,N\}$,
    \begin{equation}
        H_{c_k}(\varepsilon_k)\geq \widetilde{l}_{\gc}\normX{\varepsilon_k}^2+ \mathcal{O}\Big(\normX{\varepsilon}^2 \Lambda(L,\gc)^{\frac{1}{2}}e^{-\frac{a_d\tau L}{2}}\Big).
    \end{equation}

Moreover,
\begin{equation}\label{sum_{k=1}^N normX{varepsilon_k}^2 + sum_{k=0}^N normX{varepsilon_{k,k+1}}^2eqormX{arepsilon}^2}
    \sum_{k=1}^N \normX{\varepsilon_k}^2 + \sum_{k=0}^N \normX{\varepsilon_{k,k+1}}^2\geq\normX{\varepsilon}^2 .
\end{equation}
\end{cor}

Regarding the proof of Corollary~\ref{corollaire: almost coercivité le long des varepsilon_k}, we just add that it is reminiscent of the fact that the quadratic form $\nabla^2 (E-c_sp)(0)$ is nonnegative. We refer to the proof of Lemma~1 in~\cite{BetGrSm1} for more details. Finally, we combine the almost coercivity property of Corollary~\ref{corollaire: almost coercivité le long des varepsilon_k} with Proposition~\ref{proposition: décomposition orthogonale sans dynamique} and we derive the proof of Corollary~\ref{corollaire: controle de la fonctionnelle G}.

\begin{proof}[Proof of Corollary~\ref{corollaire: controle de la fonctionnelle G}]
    We first prove the lower bound. We show that $\mathcal{R}_{\gc,\ga}(\varepsilon)=\grandOde{\normX{\varepsilon}^3}$. For that, we need to investigate the dependence in $\gc,\ga$ of each terms controlling $\nabla^3 E(R_{\gc,\ga}+t\varepsilon)$ in Lemma~\ref{lemme: estimation différentielle E et p} (applied with $l=3$). Let us begin with the term $\normX{R_{\gc,\ga}+t\varepsilon}$. Since $Q\in\mathcal{U}_{\gc^*}(\alpha,L)$, there exists $\ga^*\in\pos_N(L)$, such that $\normX{Q-R_{\gc^*,\ga^*}}\leq \alpha$. Invoking~\eqref{normX{varepsilon(eta,v)}+Vertathfrak{C}(ta,v)-gc^*Vert_1 + Vert athfrak{A}} and taking possibly a smaller $\alpha_2$ such that $\normR{\gc-\gc^*}\leq K_1\alpha_2 < \delta^*$ with $\delta^*$ from Lemma~\ref{lemme: c,a mapsto R_c,a est lipschitz continue} leads to
\begin{align*}
    \normX{R_{\gc,\ga}+t\varepsilon}&\leq K_{lip}\big(|\gc-\gc^*|+|\ga-\ga^*|\big)+\normX{R_{\gc^*,\ga^*}}+\normX{\varepsilon}\\
    &\leq M^* + \sum_{k=1}^N \normX{Q_{c_k^*}},
\end{align*}

where $M^*:=\max(1,K_{lip})K_1\alpha_2$ only depends on $\gc^*$. Secondly, using~\eqref{1-etageq dfracmu_gc^8}, we write for any $x\in\R$,\begin{align*}
    \big|1-\big(R_{\gc,\ga}(x)+t\varepsilon(x)\big)\big|\geq 1-\eta(x)\geq \dfrac{1-\beta^*}{2}.
\end{align*}

Finally, using that $f'''$ is continuous, we can write the Taylor expansion of $f''$ to the first order, for any $x\in\R$,
\begin{align*}
    \big|f''\big(1-\eta_{\gc,\ga}(x)-t\varepsilon(x)\big)\big|\leq \big|f''(1-\eta_{\gc^*,\ga^*})\big| + M^* \Vert f'''\Vert_{L^\ii([-M^*,M^*])}.
\end{align*}

Passing to the $L^\ii$-norm, the right-hand side of the previous inequality no longer depends on $\ga^*$, hence $\normLii{f''(1-\eta_{\gc,\ga}-t\varepsilon)}$ is bounded by a constant that only depends on $\gc^*$.

Secondly, by~\eqref{mathcalUgc^*(alpha_1,L)subset mathcalO^perp}, and up to taking a larger $L_2$, Corollary~\ref{corollaire: almost coercivité le long des varepsilon_k} applies and using additionally the property~\eqref{nu_athfrak(eta,v)geq dfracnu_gc^*}, there exists $\widetilde{l}_{*}>0$ that only depends on $\gc^*$ such that
    \begin{align}\label{sum_{k=0}^N big(H_{c_k}(varepsilon_k)+H_0^k(varepsilon_{k,k+1})big)}
        \sum_{k=1}^N H_{c_k}(\varepsilon_k)+\sum_{k=0}^N  H_0^k(\varepsilon_{k,k+1})&\geq \widetilde{l}_{\gc}\left( \sum_{k=1}^N \normX{\varepsilon_k}^2+\sum_{k=0}^N \normX{\varepsilon_{k,k+1}}^2\right) + \mathcal{O}\Big(\normX{\varepsilon}^2 \Lambda(L_1,\gc)^{\frac{1}{2}}e^{-\frac{a_d\tau L_1}{2}}\Big)\notag\\
        &\geq \widetilde{l}_{*}\normX{\varepsilon}^2+\mathcal{O}\Big(\normX{\varepsilon}^2 \Lambda(L_1,\gc)^{\frac{1}{2}}e^{-\frac{a_d\tau L_1}{2}}\Big).
    \end{align}

Now, plugging~\eqref{sum_{k=0}^N big(H_{c_k}(varepsilon_k)+H_0^k(varepsilon_{k,k+1})big)} and~\eqref{nu_athfrak(eta,v)geq dfracnu_gc^*} in Corollary~\ref{corollaire: développement asymptotique de G en Q= R_c,a + varepsilon} yields to
\begin{align*}
    G(Q)&\geq \sum_{k=1}^N \big(E(Q_{c_k^*})-c_k^* p(Q_{c_k^*})\big)+ \dfrac{\widetilde{l}_{*}}{2}\normX{\varepsilon}^2 + \mathcal{O}\big(\normR{ \gc-\gc^*}^2\big)\\
    &+\mathcal{O}\Big(L(e^{-a_d\nu_{\gc}L}+e^{-a_d\tau_0 L})\Big)+\mathcal{O}\Big(\normX {\varepsilon}^2\big(\tau+e^{-\frac{a_d\tau L_1}{2}}+L_1^{\frac{1}{2}} e^{-\frac{a_d\tau L_1}{2}}\big)\Big)+ \mathcal{R}_{\gc,\ga}(\varepsilon),
\end{align*}

where we have also used that $\normX{\varepsilon}$ is bounded, and that $L^{\frac{1}{2}}$ and $\Lambda(L,\gc)$ are controlled by $\mathcal{O}(L)$. As a consequence of the both previous steps, we fix the values of $\tau$ and $\tau_0$ such that $\tau_0<2\tau <\frac{\nu_{\gc^*}}{2}\leq \nu_\gc$, and $L_2$ large enough such that 
\begin{equation*}
    \mathcal{O}\Big(\normX {\varepsilon}^2\big(\tau+e^{-\frac{a_d\tau L_2}{2}}+L_2e^{-\frac{a_d\tau L_2}{2}}\big)\Big)\leq \dfrac{\widetilde{l}_{*}}{4}\normX{\varepsilon}^2,
\end{equation*}

and we can replace 
$\mathcal{O}\Big(L(e^{-a_d\nu_{\gc}L}+e^{-a_d\tau_0 L})\Big)$ by $\mathcal{O}\Big(Le^{-a_d\tau_0 L}\Big)$. Thus, we get \begin{equation*}
    G(Q)\geq \sum_{k=1}^N \big(E(Q_{c_k^*})-c_k^* p(Q_{c_k^*})\big)+ \dfrac{\widetilde{l}_{*}}{4}\normX{\varepsilon}^2 + \grandOde{\normR{\gc-\gc^*}^2}+\mathcal{O}\big(\normX{\varepsilon}^3\big)+\mathcal{O}\Big(Le^{-a_d\tau_0 L})\Big).
\end{equation*}

Let us now tackle the upper bound. It remains to bound $H_0$ and $H_{c_k}$ from above. By construction of $\Phi_k$, we have that $\normX{\varepsilon_k}\lesssim\normX{\varepsilon}$, and by the Cauchy-Schwarz inequality, $\big|H_{c_k}(\varepsilon_k)\big|\lesssim \normLdeux{\mathcal{H}_{c_k}(\varepsilon_k)}\normLdeux{\varepsilon}$. Regarding $\nabla^2 E (Q_{c_k^*})$, we have
\begin{align*}
    \nabla^2 E(Q_{c_k})(\varepsilon_k,\varepsilon_k) &  =\nabla^2 E(Q_{c_k^*})(\varepsilon_k,\varepsilon_k)+\int_0^1 \nabla^3 E\big( (1-t)Q_{c_k^*}+tQ_{c_k}\big).(Q_{c_k}-Q_{c_k^*},\varepsilon_k,\varepsilon_k)dt.
\end{align*}
Similarly to $\mathcal{R}_{\gc,\ga}(\varepsilon)$, we use Lemma~\ref{lemme: estimation différentielle E et p} with $l=2,3$ to bound uniformly in $\gc$ the operator norm of $\nabla^2 E(Q_{c_k^*})$ and $\nabla^3 E\big( (1-t)Q_{c_k^*}+tQ_{c_k}\big)$. We deduce

\begin{align*}
    \nabla^2 E(Q_{c_k})(\varepsilon_k,\varepsilon_k) &=\grandOde{\normX{\varepsilon}^2}.
\end{align*}

The same way, we derive $\nabla^2 p(Q_{c_k})(\varepsilon_k,\varepsilon_k)=\grandOde{\normX{\varepsilon}^2}$. As a consequence, we obtain that $H_{c_k}(\varepsilon_k),H^0_{c_k}(\varepsilon_k)=\grandOde{\normX{\varepsilon}^2}$. We conclude by combining the previous bounds with Corollary~\ref{corollaire: développement asymptotique de G en Q= R_c,a + varepsilon}.
\end{proof}

\section{Dynamics of the modulation parameters}
\label{section: dynamics of the modulation parameters}
In this section, we prove Proposition~\ref{proposition: controle de ga'(t)-gc(t) + gc'(t)}. We refer to the beginning of Subsection~\ref{section: sketch of the proof evolution in time} for the definition of the time-dependent functions $\varepsilon,\gc,\ga$.

\begin{proof}[Proof of Proposition~\ref{proposition: controle de ga'(t)-gc(t) + gc'(t)}]
First of all, we consider that the initial condition $(\eta_0,v_0)\in\mathcal{NX}^2(\R)$ so that, in view of~\eqref{NLShydro}, the solution exists locally in time and belong to $\mathcal{C}^1([0,T],\Nenergyset)$. By composition with the functions $\mathfrak{C},\mathfrak{A}$ exhibited in Proposition~\ref{proposition: décomposition orthogonale sans dynamique}, we infer that $(\gc,\ga)\in\mathcal{C}^1([0,T],\R^{2N})$. Therefore, plugging $\varepsilon=(\eta,v)-R_{\gc,\ga}$ in~\eqref{NLShydro} and using that each soliton $(\eta_{c_k},v_{c_k})$ satisfies~\eqref{TWChydro} with $c=c_k$, we obtain the equations

\begin{align*}
    \partial_t\varepsilon_\eta =&\  2\partial_x (\varepsilon_\eta v_{\gc,\ga}-\varepsilon_v(1-\eta_{\gc,\ga})+\varepsilon_\eta\varepsilon_v)+2\sum_{k=1}^N\sum_{\substack{j=1 \\ j\neq k}}^N \partial_x (\eta_{c_k,a_k}v_{c_j,a_j})\\
    & -\sum_{k=1}^N c_k' \partial_c\eta_{c_k,a_k} + \sum_{k=1}^N (a_k'-c_k) \partial_x\eta_{c_k,a_k},
\end{align*}

and

\begin{align*}
    \partial_t\varepsilon_v  = &\ \partial_x\Big(2v_{\gc,\ga}\varepsilon_v +\varepsilon_v^2+\varepsilon_\eta\int_0^1 f'(1-t\eta_{\gc,\ga})dt\Big)-\partial_x\Big( f(1-\eta_{\gc,\ga})-\sum_{k=1}^N f(1-\eta_{c_k,a_k})\Big)\\
    & +\partial_x\bigg(\dfrac{\partial_x^2(\eta_{\gc,\ga}+\varepsilon_\eta)}{2(1-\eta_{\gc,\ga}-\varepsilon_\eta)}-\sum_{k=1}^N \dfrac{\partial_x^2\eta_{c_k,a_k}}{2(1-\eta_{c_k,a_k})}\bigg)+\partial_x\bigg(\dfrac{\big(\partial_x\eta_{\gc,\ga}+\partial_x\varepsilon_\eta\big)^2}{4(1-\eta_{\gc,\ga}-\varepsilon_\eta)^2}-\sum_{k=1}^N \dfrac{(\partial_x\eta_{c_k,a_k})^2}{4(1-\eta_{c_k,a_k})^2}\bigg)\\
    &-\sum_{k=1}^N\sum_{\substack{j=1 \\ j\neq k}}^N \partial_x (v_{c_k,a_k}v_{c_j,a_j})+\sum_{k=1}^N (a_k'-c_k)\partial_x v_{c_k,a_k}-\sum_{k=1}^N c_k'\partial_c v_{c_k,a_k}.
\end{align*}

Now by differentiating the orthogonality conditions~\eqref{condition d'orthogonalité sur varepsilon dans décomposition orthogonale} on $\varepsilon(t)$ with respect to $t$, and plugging both previous equations in it, we get for any $k\in \{1,...,N\}$,

\begin{align}\label{beginpmatrix gc'(t) ga'(t)-gc(t)}
    M(t)\begin{pmatrix}
        \ga'(t)-\gc(t)\\
        \gc'(t)
    \end{pmatrix}=\Phi(t),
\end{align}

with \begin{align*}
    & M_{k,j}(t)=\psLdeuxLdeux{\partial_x Q_{c_j,a_j}}{\partial_x Q_{c_k,a_k}}-\delta_{j,k}\psLdeuxLdeux{\partial^2_x Q_{c_k,a_k}}{\varepsilon}\\
    & M_{k,j+N}(t)=-\psLdeuxLdeux{\partial_c Q_{c_j,a_j}}{\partial_x Q_{c_k,a_k}}+\delta_{j,k}\psLdeuxLdeux{\partial_x\partial_c Q_{c_k,a_k}}{\varepsilon}\\
    & M_{k+N,j}(t)=\nabla p(Q_{c_k,a_k}).\partial_x Q_{c_j,a_j}-\delta_{j,k}\nabla p (\partial_x Q_{c_k,a_k}).\varepsilon,\\
    & M_{k+N,j+N}(t)=-\nabla p(Q_{c_k,a_k}).\partial_c Q_{c_j,a_j}+\delta_{j,k}\nabla p (\partial_c Q_{c_k,a_k}).\varepsilon\\
\end{align*}

\begin{align*}
    &\Phi_k  (t) =  \psLdeuxLdeux{ \begin{pmatrix}
        2v_{\gc,\ga} +2\varepsilon_v+c_k & -2(1-\eta_{\gc,\ga})\\
        \int_0^1 f'(1-t\eta_{\gc,\ga})dt & 2v_{\gc,\ga} +2\varepsilon_v+c_k
    \end{pmatrix}\varepsilon}{\partial_x^2 Q_{c_k,a_k}}\\
    & +\sum_{j=1}^N \sum_{\substack{l=1 \\ l\neq j}}^N\psLdeuxLdeux{ 
         \begin{pmatrix}2\eta_{c_j,a_j}v_{c_l,a_l}\\
        -v_{c_j,a_j}v_{c_l,a_l}
    \end{pmatrix}}{\partial_x^2 Q_{c_k,a_k}}+\psLdeux{f(1-\eta_{\gc,\ga})-\sum_{j=1}^N f(1-\eta_{c_j,a_j})}{\partial_x^2 v_{c_k,a_k}}\\
    & +\psLdeux{\dfrac{\partial_x^2(\eta_{\gc,\ga}+\varepsilon_\eta)}{2(1-\eta_{\gc,\ga}-\varepsilon_\eta)}-\sum_{k=1}^N \dfrac{\partial_x^2\eta_{c_k,a_k}}{2(1-\eta_{c_k,a_k})}+\dfrac{\big(\partial_x\eta_{\gc,\ga}+\partial_x\varepsilon_\eta\big)^2}{4(1-\eta_{\gc,\ga}-\varepsilon_\eta)^2}-\sum_{k=1}^N \dfrac{(\partial_x\eta_{c_k,a_k})^2}{4(1-\eta_{c_k,a_k})^2}}{\partial_x^2 v_{c_k,a_k}},\\
\end{align*}

and 

\begin{align*}
    &\Phi_{k+N}(t)= 2\nabla p(\partial_x Q_{c_k,a_k}).\begin{pmatrix}
        (2v_{\gc,\ga} +2\varepsilon_v+c_k)\varepsilon_\eta -2(1-\eta_{\gc,\ga})\varepsilon_v\\
        \int_0^1  f'(1-t\eta_{\gc,\ga})\varepsilon_\eta dt + (2v_{\gc,\ga} +2\varepsilon_v+c_k)\varepsilon_v
    \end{pmatrix}\\
    & +\sum_{j=1}^N \sum_{\substack{l=1 \\ l\neq j}}^N 
         2\nabla p (\partial_x Q_{c_k,a_k}).\begin{pmatrix}2\eta_{c_j,a_j}v_{c_l,a_l}\\
        -v_{c_j,a_j}v_{c_l,a_l}
    \end{pmatrix}+\psLdeux{f(1-\eta_{\gc,\ga})-\sum_{j=1}^N f(1-\eta_{c_j,a_j})}{\partial_x \eta_{c_k,a_k}}\\
    & +\psLdeux{\dfrac{\partial_x^2(\eta_{\gc,\ga}+\varepsilon_\eta)}{2(1-\eta_{\gc,\ga}-\varepsilon_\eta)}-\sum_{k=1}^N \dfrac{\partial_x^2\eta_{c_k,a_k}}{2(1-\eta_{c_k,a_k})}+\dfrac{\big(\partial_x\eta_{\gc,\ga}+\partial_x\varepsilon_\eta\big)^2}{4(1-\eta_{\gc,\ga}-\varepsilon_\eta)^2}-\sum_{k=1}^N \dfrac{(\partial_x\eta_{c_k,a_k})^2}{4(1-\eta_{c_k,a_k})^2}}{\partial_x \eta_{c_k,a_k}}.\\
\end{align*}

We decompose $M(t)=D(t)+H(t)$ where all the entries of $D(t)$ are zero except the coefficients
\begin{align*}
    & D_{k,k}(t):=\normLdeuxLdeux{\partial_x Q_{c_k}}^2,\\
    & D_{k,k+N}(t):=-\psLdeuxLdeux{\partial_c Q_{c_k,a_k}}{\partial_x Q_{c_k,a_k}},\\
    & D_{k+N,k}(t):=\nabla p(Q_{c_k,a_k}).\partial_x Q_{c_k,a_k}.\\
    & D_{k+N,k+N}(t):=-\nabla p(Q_{c_k}).\partial_c Q_{c_k},\\
\end{align*}

\begin{claim}
    The matrix $D(t)$ is a diagonal matrix. Furthermore, for $\alpha$ small enough, $D(t)$ is invertible and $\left| D(t)^{-1}\right|$ is bounded by a constant that only depends on $\gc^*$.
\end{claim}

\begin{proof}
First we show that the coefficients above and below the diagonal are zero. The coefficients below the diagonal are zero by doing an integration by part. The coefficients above the diagonal reads \begin{equation*}
    -\psLdeuxLdeux{\partial_c Q_{c_k,a_k}}{\partial_x Q_{c_k,a_k}}=-\int_\R\Big( \partial_c\eta_{c_k}\partial_x\eta_{c_k}+\partial_c v_{c_k}\partial_x v_{c_k}\Big).
\end{equation*}

Since $\eta_c$ and $v_c$ are even, then the derivatives with respect to $x$ are odd and the derivatives with respect to $c$ are even. Thus the previous coefficient is zero by integrating on the real line an odd function. We obtain $D(t)=\text{diag}\big(D_{1,1}(t),...,D_{2N,2N}(t)\big)$. For $\alpha\leq \alpha_1$,~\eqref{nu_athfrak(eta,v)geq dfracnu_gc^*} imply that for every $k\in\{1,...,N\}$ and $t\in [0,T]$, $c_k(t)\in (c_0,c_s)$. Moreover $D_{k+N,k+N}(t)=-\frac{d}{dc}\big(p(Q_c)\big)_{|c=c_k}<0$ by Theorem~\ref{théorème: il existe une branche C^1 de solitons proche de c_s}. Then the diagonal coefficients of $D(t)$ do not vanish. Therefore, $D(t)$ is invertible and it remains to bound from below its operator norm. By Lemma~\ref{lemme: c,a mapsto R_c,a est lipschitz continue}, we have for any $k\in\{1,...,N\}$ and $t\in [0,T]$,
    \begin{equation}
        \normLdeux{\partial_x Q_{c_{k}^*}}\leq K_* |\gc^*-\gc| +\sqrt{D_{k,k}(t)}.
    \end{equation}
    and by~\eqref{nu_athfrak(eta,v)geq dfracnu_gc^*},
    \begin{equation}
        \dfrac{\kappa_{\gc^*}}{2}\leq D_{k+N,k+N}(t).
    \end{equation}

There exists $\alpha_3\leq \alpha_2$ sufficiently small such that 
\begin{equation}
    K_* K_1 \alpha_3 \leq\min\left\{ \dfrac{1}{2}\normLdeux{\partial_x Q_{c_{k}^*}},\dfrac{\kappa_{\gc^*}}{2}\right\}.
\end{equation}

Since $Q(t)\in\mathcal{U}_{\gc^*}(\alpha_1,L^1)$ for any $t\in [0,T]$,~\eqref{normX{varepsilon(eta,v)}+Vertathfrak{C}(ta,v)-gc^*Vert_1 + Vert athfrak{A}} applies and we obtain the desired uniform below bound (with respect to $t$ and $k$) for $D_{k,k}(t)$ and $D_{k+N,k+N}(t)$.
\end{proof}

Now, since $Q_{\gc,\ga}=\mathcal{O}(1)$ (by exponential decay and~\eqref{nu_athfrak(eta,v)geq dfracnu_gc^*}),
\begin{equation*}
    H_{j,j}(t)=-\psLdeuxLdeux{\partial^2_x Q_{c_j,a_j}}{\varepsilon}=\mathcal{O}(\normX{\varepsilon}).
\end{equation*}
By Lemma~\ref{lemme: intégrale d'un polynome couplé est un O(Lexp(-Lnu)}, since $\ga(t)\in\pos_N(L)$ with $L\geq L_1$, we have, for $j\neq k\in\{1,...,N\}$,
\begin{equation*}
    H_{k,j}(t)=\psLdeuxLdeux{\partial_x Q_{c_j,a_j}}{\partial_x Q_{c_k,a_k}}=\grandOde{Le^{-\frac{a_d\nu_{\gc^*}L}{2}}},
\end{equation*}

where we used that $\Lambda(L,\gc)$ is controlled by $\grandOde{L}$, by~\eqref{nu_athfrak(eta,v)geq dfracnu_gc^*}. We deal with the other terms the same way and we infer that $|H(t)|=\mathcal{O}(\normX{\varepsilon})+\mathcal{O}\left(Le^{-\frac{a_d\nu_{\gc^*} L}{2}}\right)$. Moreover~\eqref{normX{varepsilon(eta,v)}+Vertathfrak{C}(ta,v)-gc^*Vert_1 + Vert athfrak{A}} is satisfied and taking possibly a smaller $\alpha_3$ and $L_3\geq L_2$ large enough, then for $\alpha\leq \alpha_3$ and $L\geq L_3$, $\left|D(t)^{-1}H(t)\right|<1$, so that 
$M(t)=D(t)\big( I+D(t)^{-1}H(t)\big)$ is invertible and $\left| M(t)^{-1}\right|$ is bounded by a constant that only depends on $\gc^*$.

Now we show that $\left|\Phi(t)\right|=\mathcal{O}(\normX{\varepsilon})+\mathcal{O}\left(Le^{-\frac{a_d\nu_{\gc^*} L}{2}}\right)$. This is the same argument than for $H(t)$, since there only appear terms that involve $\varepsilon$ or multivariate polynomials in $\mathcal{P}_N$ evaluated in solitons. For sake of completeness, we deal with the terms involving the nonlinearity $f$. By Lemma~\ref{lemme: F,f, f'(1-eta_c)}, we get
\begin{align*}
    \bigg|\bigg\langle f(1-\eta_{\gc,\ga})-\sum_{j=1}^N f(1-\eta_{c_j,a_j}),\partial^2_x v_{c_k,a_k}\bigg\rangle_{L^2}\bigg|&=\grandOde{\normLun{(S^{2,N}+B^N)(\eta_{\gc,\ga})}\normLii{\partial_x^2 v_{c_k}}}\\
    &=\mathcal{O}\Big(Le^{-\frac{a_d\nu_{\gc^*} L}{2}}\Big).
\end{align*}

We use that $f'$ is continuous and the fact that $\eta_{\gc,\ga}=\mathcal{O}(1)$ in order to deduce \begin{equation*}
    \psLdeux{\varepsilon_\eta \int_0^1 f'(1-t\eta_{\gc,\ga})dt}{\partial_x^2 v_{c_k,a_k}}=\grandOde{\normX{\varepsilon}}.
\end{equation*}

We have shown~\eqref{Vert ga'(t)-gc(t)Vert_1 + Vert gc'(t)Vert_1 =mathcal(normXvarepsilon(t))}, provided that the initial condition $Q_0$ lies in $\mathcal{NX}^2(\R)$. Now taking $Q_0\in\Nenergyset$, we approach by density this initial value by a sequence in $\mathcal{NX}^2(\R)$ and we argue as in~\cite{BetGrSm1}. The only point where this density approach differs from the Gross-Pitaevskii case is the nonlinearity $f$, which appears in the equation of $\partial_t\varepsilon_v$ and in the expression of $\Phi(t)$. By hypothesis, $f$ is at least $\mathcal{C}^1$ and in view of the expression of $M(t)$ and $\Phi(t)$, this is sufficient to pass to the limit in~\eqref{beginpmatrix gc'(t) ga'(t)-gc(t)} along the approaching sequence. Moreover the convergence is in $\mathcal{C}^0([0,T],\R^{2N})$. This eventually shows that  $\gc,\ga\in\mathcal{C}^1\big([0,T],\R^{2N}\big)$ and that for any $t\in [0,T]$, \begin{equation}\label{preuve: ga'(t)-gc(t) + gc'(t) leq O...}
    \left| \ga'(t)-\gc(t)\right| + \left| \gc'(t)\right| = \mathcal{O}\Big(Le^{-\frac{a_d\nu_{\gc^*} L}{2}}\Big)+\grandOde{\normX{\varepsilon(t)}}.
\end{equation}

Finally, we derive~\eqref{a_k+1(t)-a_k(t)geq a_k+1(0)-a_k(0)+igma^* t >L-1+sigma^* t} and~\eqref{sqrt c_s^2-ig(a_k'(t)big)^2geq dfracnu} from~\eqref{preuve: ga'(t)-gc(t) + gc'(t) leq O...}, by using Proposition~\ref{proposition: décomposition orthogonale sans dynamique}, and taking possibly a further $(\alpha_3,L_3)$.
\end{proof}

\section{Monotonicity and final estimates}
\label{section: monotonie et final estimates}
This section is dedicated to establishing the almost monotonicity in time of the quantity defined for $k\in\{1,...,N\}$, by
\begin{equation}
    \widetilde{p}_k(\eta,v)=\dfrac{1}{2}\int_\R \chi_k\eta v,
\end{equation}
where
\begin{equation*}
    \chi_{k}(x)=\left\{
\begin{array}{l}
    1\quad\text{if }k=1, \\
    \dfrac{1}{2}\bigg(1+\tanh\left(\tau_0\Big(x-\frac{a_k+a_{k-1}}{2}\Big)\right)\bigg)\quad\text{if }k\in\{2,...,N\}. \\
\end{array}
\right.
\end{equation*}

Before the proof, we state a conservation type formula for the momentum. Like in Section~\ref{section: dynamics of the modulation parameters}, this proof also relies on an approximation of the initial condition by a sequence in $\mathcal{NX}^2(\R)$ and we refer to the proof of Corollary 3.1 in~\cite{BetGrSm1} for more details.

\begin{lem}\label{lemme: d/dt ( int psi eta v)}
Let $(\eta,v)\in\mathcal{C}^0([0,T],\Nenergyset)$ be a solution of~\eqref{NLShydro} and $\widetilde{\chi}\in\mathcal{C}^0\big([0,T],\mathcal{C}^3_b(\R)\big)\cap\mathcal{C}^1\big([0,T],\mathcal{C}^0_b(\R)\big)$, then $t\mapsto\psLdeux{\widetilde{\chi}}{\eta v}$ is differentiable and its derivative is
\begin{align}
    \dfrac{d}{dt}\left(\int_\R \widetilde{\chi}\eta v\right)=&\int_\R ·\partial_t\widetilde{\chi} \eta v+ \int_\R \partial_x\widetilde{\chi}\Big((1-2\eta)v^2 + \widetilde{F}(\eta)+\dfrac{(3-2\eta)(\partial_x\eta)^2}{4(1-\eta)^2}\Big) \notag \\ 
    & + \dfrac{1}{2}\int_\R \partial_x^3\widetilde{\chi} \big(\eta + \ln(1-\eta)\big)\label{d/dt ( int psi eta v)},
\end{align}
where $\widetilde{F}(\rho)=\rho f(1-\rho)-F(1-\rho)$.
\end{lem}

We notice that applying Lemma~\ref{lemme: d/dt ( int psi eta v)} to $\widetilde{\chi}\equiv 1$ implies the conservation of the momentum. In this sense, we also mention that the case $k=1$ in~\eqref{formule de monotonie sur p_k} is a consequence of Remark~\ref{remarque: interprétation de widetilde p_k}. For $k\in\{2,...,N\}$, we apply Lemma~\ref{lemme: d/dt ( int psi eta v)} to $\widetilde{\chi}(x,t):=\chi\left(\tau_0\big(x-X_k(t)\big)\right)$ where $\chi(x)=\frac{1}{2}(1+\tanh(x))$ and $X_k(t)=\frac{1}{2}\big(a_k(t)+a_{k-1}(t)\big)$. Therefore we can differentiate $\widetilde{p}_k$ and it leads to \begin{align*}
    \dfrac{d}{dt}\widetilde{p}_k(t)=\int_\R \bigg(\tau_0\chi'\Big(\tau_0\big(x-&X_k(t)\big)\Big)\Big(\widetilde{q}\big((\eta(t,x),v(t,x))\big)-X_k'(t)\eta v\Big)\\
    & + \dfrac{\tau_0^3}{2}\chi'''\Big(\tau_0\big(x-X_k(t)\big)\Big)\big(\eta(t,x)+\ln(1-\eta(t,x)\big)\bigg)dx,
\end{align*}

where $$\widetilde{q}(\eta,v)=(1-2\eta)v^2+\widetilde{F}(\eta)+\dfrac{(3-2\eta)(\partial_x\eta)^2}{4(1-\eta)^2}.$$

According to the arguments in~\cite{MarMeTs1}, we decompose this integral into two complementary parts. Namely, we decompose the line as an interval $I$ focusing on the area where $\chi$ varies and a second part which takes into account the whereabouts of the solitons in which $\chi'$ is small. This study can be summarized as follows.

\begin{lem}\label{lemme: decomposition du viriel sur I et I^c}
Given an interval $I$, there exist two positive constants $C$ and $C_{\ln}$ that only depends on $\gc^*$ such that we have
\begin{equation}\label{decomposition du viriel sur I et I^c}
    \dfrac{d}{dt}\widetilde{p}_k(t)\geq -C \sup_{x\in I^c}\chi'\Big(\tau_0\big(x-X_k(t)\big)\Big) + \tau_0\int_I q(\eta,v)\chi'\Big(\tau_0\big(x-X_k(t)\big)\Big)dx,
\end{equation}

where \begin{equation*}
    q(\eta,v)=(1-2\eta)v^2-\sqrt{c_s^2-\dfrac{\nu_{\gc^*}^2}{4}}|\eta v| +\left(-\int_0^1 rf'(1-r\eta)dr -\tau_0^2 C_{\ln}\right)\eta^2 .
\end{equation*}

\end{lem}

\begin{proof}
To obtain the second term in the right-hand side of~\eqref{decomposition du viriel sur I et I^c}, we first use~\eqref{sqrt c_s^2-ig(a_k'(t)big)^2geq dfracnu} to get $X_k'(t)^2\leq c_s^2-\frac{\nu_{\gc^*}^2}{4}$ and also write a first order Taylor expansion between $\eta$ and 0 which reads
\begin{equation}\label{first order taylor widetilde F entre 0 et eta}
    \widetilde{F}(\eta)=-\eta^2\int_0^1 r f'(1-r\eta)dr.
\end{equation}
Moreover, by~\eqref{1-etageq dfracmu_gc^8}, $\eta\leq \frac{1+\beta^*}{2}<1$, so that $\frac{(3-2\eta)(\partial_x\eta)^2}{4(1-\eta)^2}\geq 0$. Some elementary analysis on the function $\xi\mapsto \xi +\ln(1-\xi)$ provides the constant $C_{\ln}$ such that for any $\xi\in[0,\frac{1+\beta^*}{2}]\subset [0,1)$, we have $\big|\xi+\ln(1-\xi)|\leq \frac{C_{\ln}}{4}\xi^2$. Combining this and the fact that $|\chi'''|\leq 8\chi'$ yields to the right-hand side term in~\eqref{decomposition du viriel sur I et I^c}.\\
Concerning the complementary set $I^c$, we just bound, up to a constant, the integrand of $\frac{d}{dt}\widetilde{p}_k(t)$ by the energy density $e(\eta,v)$. This can be dealt as in~\cite{BetGrSm1}, apart from the terms involving $\eta^2$. For these terms, we use~\eqref{hypothèse de croissance sur F minorant intermediaire} and~\eqref{first order taylor widetilde F entre 0 et eta}, to get
\begin{equation*}
    |\widetilde{F}(\eta)+ \eta + \ln(1-\eta)|\leq \bigg(2\Vert f'\Vert_{L^\ii([0,2])}+C_{\ln}\bigg) \dfrac{ F(1-\eta)}{c_s^2},
\end{equation*} 
and the conclusion follows.
\end{proof}

By continuity of $f'$, there exists $\gd\in (0,\frac{1}{2}]$ such that \begin{equation}\label{définition f^delta}
    m_\gd:=\min_{|x|\leq \gd}\big(-f'(1+x)\big)>0,
\end{equation}
because of~\eqref{f'(1) <0}. We can now conclude the proof of Proposition~\ref{proposition: monotonicity}.

\begin{proof}[Proof of Proposition~\ref{proposition: monotonicity}]
We apply Lemma~\ref{lemme: decomposition du viriel sur I et I^c} to 
$$I=I_k(L):=\left[X_k(t)-\dfrac{L-1+\sigma^* t}{4},X_k(t)+\dfrac{L-1+\sigma^*t}{4}\right],$$
with some suitable choice of $(\alpha,L)$ that is exhibited later.

     \underline{Step 1: Nonnegativity on $I_k(L)$.} According to the orthogonal decomposition of $\eta$, we write for any $t\in [0,T]$, $\eta(t,x)=\eta_{\gc(t),\ga(t)}(x)+\varepsilon_\eta(t,x)$. By~\eqref{normX{varepsilon(eta,v)}+Vertathfrak{C}(ta,v)-gc^*Vert_1 + Vert athfrak{A}} and the one-dimensional Sobolev embedding, there exists $\alpha_4$ small enough such that for $\alpha\leq \alpha_4$ and $(t,x)\in [0,T]\times I_k(L)$, 
    \begin{equation*}
        \big|\eta(t,x)\big|\leq \sum_{j=1}^N \big|\eta_{c_j(t),a_j(t)}(x)\big|+\dfrac{\gd}{2}.
    \end{equation*}
    By exponential decay, we also have for any $j\in\{1,...,N\}$, $\big|\eta_{c_j(t),a_j(t)}(x)\big|\leq K_d e^{-a_d\nu_{c_j(t)}|x-a_j(t)|}$. By virtue of~\eqref{a_k+1(t)-a_k(t)geq a_k+1(0)-a_k(0)+igma^* t >L-1+sigma^* t}, we have for any $j$, and $(t,x)\in [0,T]\times I_k(L)$,
    \begin{align*}
        |x-a_j(t)|&\geq \big| X_k(t)-a_j(t)\big|-\big|X_k(t)-x|\geq \dfrac{L-1}{4}.
    \end{align*}
Therefore, using in addition~\eqref{nu_athfrak(eta,v)geq dfracnu_gc^*}, there exists $L_4$ large enough such that for $L\geq L_4$, for any $(t,x)\in [0,T]\times I_k(L)$, $\big|\eta(t,x)\big|\leq \gd$. Now labelling $\widetilde{\sigma}_0(\eta):=-\int_0^1 rf'(1-r\eta)dr-\tau_0^2 C_{\ln}$, we first notice that $\widetilde{\sigma}_0(\eta)\geq \frac{m_\gd}{2}+\grandOde{\tau_0^2}$. We impose that $\tau_0$ is small enough so that $\widetilde{\sigma}_0(\eta)>0$. We apply Gauss' reduction method to compute
\begin{equation*}
    q(\eta,v)=\widetilde{\sigma}_0(\eta)\left( \eta-\dfrac{\sigma_0}{2\widetilde{\sigma}_0(\eta)}v\right)^2+\left(1-2\eta-\dfrac{\sigma_0^2}{4\widetilde{\sigma}_0(\eta)}\right)v^2,
\end{equation*}
where $\sigma_0:=\sqrt{c_s^2-\frac{\nu_{\gc^*}^2}{4}}$. Now, we claim the following, and conclude the proof.
\begin{claim}\label{claim: Up to taking smaller parameters gd and tau_0}
    Up to taking smaller parameters $\gd$ and $\tau_0$, we have for any $(t,x)\in [0,T]\times I_k(L)$, \begin{equation*}
        \eta(t,x)\leq \gd\text{ and } 1-2\eta-\dfrac{\sigma_0^2}{4\widetilde{\sigma}_0(\eta)}\geq 0.
    \end{equation*}
\end{claim}

We then deduce that the integral in~\eqref{decomposition du viriel sur I et I^c} is positive.

\underline{Step 2: Conclusion.} Since $\chi'(y)\leq 2e^{-2|y|}$ for any $y\in\R$, we infer that for $x\notin I_k(L)$, $\chi'\big(\tau_0(x-X(t))\big)\leq 2e^{-\tau_0\frac{L-1+\sigma^* t}{2}}$. Combining the estimate on $I_k(L)^c$ with Step 1, we conclude that up to reducing the value of $\gd$ defined in~\eqref{définition f^delta} and the value of $\tau_0$, there exists suitable $(\alpha_4,L_4)$ such that for any $(\alpha,L)\prec (\alpha_4,L_4)$, and any $t\in [0,T]$, the monotonicity formula~\eqref{formule de monotonie sur p_k} can be derived from~\eqref{decomposition du viriel sur I et I^c}. To finish, we just state that~\eqref{formule de monotonie sur G} is a straightforward corollary of the monotonicity formula on $\widetilde{p}_k$.
\end{proof}

To finish this section, we prove Claim~\ref{claim: Up to taking smaller parameters gd and tau_0}.
\begin{proof}[Proof of Claim~\ref{claim: Up to taking smaller parameters gd and tau_0}]
The way we have constructed $(\alpha_4,L_4)$ only depends on $\gd$. Then we always can take $(\alpha_4,L_4)$ such that $|\eta|\leq \gd$. Secondly,  \begin{equation*}
    1-2\eta-\dfrac{\sigma_0^2}{4\widetilde{\sigma}_0(\eta)}\geq 0 \Longleftrightarrow \eta\leq\dfrac{1}{2}-\dfrac{\sigma_0^2}{8\widetilde{\sigma}_0(\eta)}.
\end{equation*}
Since $\eta$ is bounded and $f''$ is continuous, and using~\eqref{relation entre vitesse c_s et f'(1)}, we have 
\begin{equation*}
    \dfrac{1}{2}-\dfrac{\sigma_0^2}{8\widetilde{\sigma}_0(\eta)} =\dfrac{1}{2}-\dfrac{c_s^2-\frac{\nu_{\gc^*}}{4}}{8\Big(\frac{-f'(1)}{2}+\grandOde{\eta}+\grandOde{\tau_0^2}\Big)}=\dfrac{\nu_{\gc^*}^2}{8c_s^2}+\grandOde{\gd^2}+\grandOde{\tau_0^2}.
\end{equation*}

Taking possibly smaller $\gd$ and $\tau_0$, we obtain \begin{equation*}
    \eta\leq\gd\leq \dfrac{1}{2}-\dfrac{\sigma_0^2}{8\widetilde{\sigma}_0(\eta)}.
\end{equation*}
\end{proof}

Here we give the proof of Proposition~\ref{proposition: controle de ga(t), c(t)-c^* et varepsilon(t)}.
\begin{proof}
Assertion~\eqref{ga(t) in pos_N(L-1)} is a straightforward deduction from~\eqref{a_k+1(t)-a_k(t)geq a_k+1(0)-a_k(0)+igma^* t >L-1+sigma^* t}. Now we prove simultaneously the remaining assertions. Integrating~\eqref{formule de monotonie sur G}, yields the upper bound \begin{equation}\label{mathcalG(t)-athcalG(0)leq}
        \mathcal{G}(t)-\mathcal{G}(0)\leq \grandOde{Le^{-a_d\tau_0 L}}.
    \end{equation}
Combining~\eqref{mathcalG(t)-athcalG(0)leq}, the upper bound on $\mathcal{G}(0)$ and the lower bound on $\mathcal{G}(t)$ in Corollary~\ref{corollaire: controle de la fonctionnelle G} then provides, for any $t\in [0,T]$,
    \begin{equation}\label{normXvarepsilon(t)^2leq grandOdenormXvarepsilon(0)^2+grandOdenormRgc(0)-gc^*^2}
        \normX{\varepsilon(t)}^2\leq \grandOde{\normX{\varepsilon(0)}^2}+\grandOde{\normR{\gc(0)-\gc^*}^2} +\grandOde{\normX{\varepsilon(t)}^3}+\grandOde{\normR{\gc(t)-\gc^*}^2} +\grandOde{Le^{-a_d\tau_0 L}}.
    \end{equation}
    
Now, we write
\begin{equation}
    \normR{\gc(t)-\gc^*}\leq \normR{\gc(t)-\gc(0)}+\normR{\gc(0)-\gc^*}.
\end{equation}
Arguing the same way than in the proof of Corollary~\ref{corollaire: décomposition orthogonale} (see~\eqref{borne de d/dc^2 p(Q_c)}), we obtain that $\frac{d^2}{dc^2}\big(p(Q_c)\big)$ is uniformly bounded on $[c_k(t),c_k(0)]$ (or $[c_k(0),c_k(t)]$). Then, a Taylor expansion and~\eqref{nu_athfrak(eta,v)geq dfracnu_gc^*} provides $\alpha_5>0$ small enough such that if $\alpha\leq \alpha_5$, we have uniformly in $t$,
\begin{equation}
    \left| p(Q_{c_k(0)})-p(Q_{c_k(t)})\right| \geq \dfrac{\kappa_{\gc^*}}{2}|c_k(t)-c_k(0)|+\grandOde{\normR{c_k(t)-c_k(0)}^2}, 
\end{equation}
hence
\begin{equation}\label{normRgc(t)-gc^*eq grandOdeleft| p(Q_c_k(0))-p(Q_c_k(t))}
    \normR{\gc(t)-\gc^*}\leq \grandOde{\left| p(Q_{c_k(0)})-p(Q_{c_k(t)})\right|}+\normR{\gc(0)-\gc^*}.
\end{equation}

From the definition of $\widetilde{p}$ and Proposition~\ref{proposition: développement de p_k en Q=R_c,a + varepsilon} and the choice of $\tau_0$ made in the proof of Corollary~\ref{corollaire: controle de la fonctionnelle G}, we have for any $k\in\{1,...,N\}$ and any $t\in [0,T]$, \begin{equation}
    \widetilde{p}_k(t)-\widetilde{p}_{k+1}(t)=p_k(t)=p(Q_{c_k(t)})+\grandOde{\normX{\varepsilon(t)}^2}+\grandOde{Le^{-a_d \tau_0 L}},
\end{equation}
hence \begin{equation}\label{big|p(Q_c_k(0))-p(Q_c_k(t))big|leq sum_=k^k+1}
    \big|p(Q_{c_k(0)})-p(Q_{c_k(t)})\big|\leq \sum_{l=k}^{k+1} \big|\widetilde{p}_l(0)-\widetilde{p}_l(t)\big|+\grandOde{\normX{\varepsilon(t)}^2}+\grandOde{\normX{\varepsilon(0)}^2}+\grandOde{Le^{-a_d \tau_0 L}}.
\end{equation}

By~\eqref{formule de monotonie sur p_k},
\begin{equation*}
    \big|\widetilde{p}_l(0)-\widetilde{p}_l(t)\big|=\grandOde{Le^{-a_d \tau_0 L}},
\end{equation*}
therefore, by plugging the previous estimate in~\eqref{big|p(Q_c_k(0))-p(Q_c_k(t))big|leq sum_=k^k+1} and by using~\eqref{normRgc(t)-gc^*eq grandOdeleft| p(Q_c_k(0))-p(Q_c_k(t))}, we derive~\eqref{normRgc(t)-gc^*=grandOdenormXvarepsilon(0)} . Furthermore, taking possibly $\alpha_5$ smaller, there also exists $L_5$ large enough such that if $L\geq L_5$,~\eqref{normXvarepsilon(t)^2leq grandOdenormXvarepsilon(0)^2+grandOdenormRgc(0)-gc^*^2} provides~\eqref{normXvarepsilon(t)^2 leq mathcalO(normXvarepsilon(0)^2)}.
\end{proof}

\appendix

\section{Useful estimates}
\label{appendix}
In the following lemmas, we shall make a crucial use of the exponential decay of the solitons.
\begin{lem}\label{lemme: c,a mapsto R_c,a est lipschitz continue}
Set $\gc^*:=(c_1^*,...,c_N^*)$. Then there exists $K_{lip}>0$ only depending on $\gc^*$ such that for any $\ga:=(a_1,...,a_N),(a_1^*,...,a_N^*)\in\R^N$ and $\gc\in\mathcal{B}(\gc^*,\delta^*)$ with $\delta^* = \min(\frac{\mu_{\gc^*}}{2},\frac{\nu_{\gc^*}}{2})$, we have
\begin{equation*}
    \normX{R_{\gc^*,\ga^*}-R_{\gc,\ga}}\leq K_{lip}\big( \normR{\gc^*-\gc}+\normR{\ga^*-\ga}\big).
\end{equation*}
\end{lem}

\begin{proof}
    By invariance under translation, we can assume $\ga=0$. In order to obtain some Lipschitz estimate, we infer the following control, using the assumptions on $\delta^*$ and~\eqref{estimée décroissance exponentielle à tout ordre pour eta_c et v_c}. For $l \in\{0,1\}$, we have
\begin{align*}
    \normR{\nabla_{(a,c)}\partial_x^l\eta_{c_k,a_k}(x)}+\normR{\nabla_{(a,c)}v_{c_k,a_k}(x)} 
    &\leq K_d(1+\nu_{c_k}^2)\Big(1+\dfrac{1}{c_k^3}\Big)e^{-a_d\nu_{c_k}|x-a_k|}\\
    &\leq K_d(1+\nu_{c_k^*-\delta^*}^2)\Big(1+\dfrac{1}{\mu_{\gc^*}^3}\Big)e^{-a_d\nu_{c_k^*+\delta^*}|x-a_k|}
\end{align*}
Therefore, $\normX{\nabla_{(c,a)}Q_{c_k,a_k}}$ is bounded by a constant $K_{lip}$, depending only on $\gc^*$, for any $(c,a)$ and this leads to the desired estimate.
\end{proof}

\begin{rem}\label{remarque: lemma lipschitz marche toujours avec normes plus restrictives}
    We can use similarly the exponential decay of Subsection~\ref{subsection: exponential decay} in order to show that Lemma~\ref{lemme: c,a mapsto R_c,a est lipschitz continue} still holds if we take the norms $\Vert .\Vert_{\mathcal{X}^k}$ with $k\in\{1,...,4\}$.
\end{rem}

Here we give a proof of Lemma~\ref{lemme: norm Lii de eta_c^* inférieur à delta}. We define
\begin{equation*}
    \ge_{\gc}:=\max\left\{ \normLii{\eta_{c_k}}\big| k\in\{1,...,N\}\right\}<1,
\end{equation*}
by~\eqref{max eta_c <1}.

\begin{proof}[Proof of Lemma~\ref{lemme: norm Lii de eta_c^* inférieur à delta}]
    We define for $k\in\{1,...,N \}$, $L>0$ and $\ga=(a_1,...,a_N)\in\pos_N(L)$, the sets
\begin{equation*}
    J_k:=\left\{ x\in\R\Big| |x-a_k|<\frac{L}{3}\right\}.
\end{equation*}
By construction $J_k\cap J_j=\emptyset$. Set $J=\cup_{k=1}^N J_k$. For all $k\in\{1,...,N\}$ and $x\in J_k$, $\left| \eta_{c_k^*,a_k}(x)\right|\leq \ge_{\gc^*}< 1$. Then
\begin{equation*}
    \left|\sum_{k=1}^N \eta_{c_k^*,a_k}(x)\mathds{1}_{J_k}(x)\right|\leq \ge_{\gc^*}.
\end{equation*}

By exponential decay~\eqref{estimée décroissance exponentielle à tout ordre pour eta_c et v_c}, we write for any $x\in\R$,
\begin{align*}
    |\eta_{\gc^*,\ga}(x)|&\leq \left|\sum_{k=1}^N \eta_{c_k^*,a_k}(x)\mathds{1}_{J_k}(x)\right| +\left|\sum_{k=1}^N \sum_{\substack{j=1 \\j\neq k}}^N\eta_{c_j^*,a_j}(x)\mathds{1}_{J_k}(x)\right|+\left|\sum_{k=1}^N \eta_{c_k^*,a_k}(x)\mathds{1}_{J^c}(x)\right|\\
    &\leq \ge_{\gc^*} + K_d\sum_{k=1}^N \sum_{\substack{j=1 \\j\neq k}}^N e^{-\nu_{\gc^*}|x-a_j|}\mathds{1}_{J_k}(x) + K_d\sum_{k=1}^N e^{-\nu_{\gc^*}|x-a_k|}\mathds{1}_{J^c}(x).
\end{align*}

Since $\ga\in\pos_N(L)$, for $x\in J_j$ with $j\neq k$, 
\begin{equation*}
    e^{-\nu_{\gc^*}|x-a_k|}\leq e^{-\frac{2\nu_{\gc^*}L}{3}}.
\end{equation*}

Moreover, for $x\in J^c$, and any $k\in\{1,...,N \}$, 
\begin{equation*}
    e^{-\nu_{\gc^*}|x-a_k|}\leq e^{-\frac{\nu_{\gc^*}L}{3}}.
\end{equation*}

We conclude the proof by taking $L_6$ large enough such that
\begin{equation*}
    \beta^*:=\ge_{\gc^*}+ N(N-1)K_de^{-\frac{2\nu_{\gc^*}L_6}{3}}+K_dN e^{-\frac{\nu_{\gc^*}L_6}{3}} <1,
\end{equation*}
and choose $\alpha_6<1-\beta^*$.
\end{proof}

We deduce from~\eqref{expression de l'énergie} and~\eqref{expression du moment} the following result that provides a control on the operator $\mathcal{X}$-norm of the first derivative of the energy and the momentum. 
\begin{lem}\label{lemme: estimation différentielle E et p}
For $Q,\varepsilon\in\energysethydro$, we have 
\begin{equation*}
    \left|\nabla p(Q).\varepsilon\right|\leq \normLdeuxLdeux{Q}\normX{\varepsilon}\quad\text{and}\quad \left|\nabla^2 p(Q).(\varepsilon,\varepsilon)\right|\leq \normX{\varepsilon}^2.
\end{equation*}

Moreover there exists a constant $K_E>0$, and integers $i_E,\widetilde{i}_E$ such that for $Q=(\eta,v)\in\Nenergysethydro$ and $\varepsilon\in\energysethydro$, we have for $l\in\{1,2,3\}$, 
\begin{equation*}
    \left|\nabla^l E(Q).(\varepsilon_l)\right|\leq K_E\left( 1 +\dfrac{1}{\inf_\R(1-\eta)^{i_E}}+\normX{Q}^{\widetilde{i}_E}+\normLii{f^{(l-1)}(1-\eta)}\right)\normX{\varepsilon}^l,
\end{equation*}

where $\varepsilon_l$ designates $\varepsilon,(\varepsilon,\varepsilon)$ or $(\varepsilon,\varepsilon,\varepsilon)$ according to the value of $l$.
\end{lem}

\begin{proof}
    The two first derivatives of the energy can be directly deduced from the expressions~\eqref{expression de nabla E(R_{gc,ga}).varepsilon},~\eqref{expression de nabla^2 E(R_{gc,ga}).varepsilon}. Regarding the third one, we compute
\begin{equation*}
        \nabla^3 E(Q)(\varepsilon,\varepsilon,\varepsilon)=\int_\R \left(\dfrac{3}{4(1-\eta)^2}+\dfrac{3}{2(1-\eta)^3}+\dfrac{3}{8(1-\eta)^4}+\dfrac{f''(1-\eta)}{2}\right)\varepsilon_\eta^3-\int_\R \varepsilon_\eta^2\varepsilon_v,
    \end{equation*}
    and deduce the desired estimate.
\end{proof}

\section{Exponential decay for multivariate polynomials}\label{section: several polynomials}
The following lemma gives an explicit control on the $L^p$-norm of a product of decaying and translated exponential functions. Its proof can be found in~\cite{BetGrSm1}.

\begin{lem}\label{lemme: controle de norme L^p de deux exponentielles}
    Let $(a,b)\in\R^2$ with $a<b$, $(\nu_a,\nu_b)\in(\R_+^*)^2$ and set $y^\pm=\max(\pm y,0)$, then 

    \begin{equation*}
        \normLp{ e^{-\nu_a(.-a)^+}e^{-\nu_b(.-b)^-}}\leq \bigg(\dfrac{2}{p\min(\nu_a,\nu_b)}+b-a\bigg)^{\frac{1}{p}}e^{-\min(\nu_a,\nu_b)(b-a)},
    \end{equation*}

    for all $p\in [1,...,+\ii]$ and with the convention $\frac{1}{p}=0$ if $p=+\ii$.
\end{lem}

For $M\geq 2$, we recall the following linear subspace of the set of multivariate real polynomials:

\begin{align*}
    \mathcal{P}_M=\mathcal{P}[X_1,...,X_M]:=\Big\{ \sum_{|\alpha|\leq m} p_\alpha X^\alpha \Big| m\in\N, p_{ke_i}=0\ ,\forall (k,i)\in \{0,...,m\}\times\{1,...,M\}\Big\}.
\end{align*}

In light of the previous lemma, plugging functions that decay exponentially into such polynomials, we can give the following proof of Lemma~\ref{lemme: intégrale d'un polynome couplé est un O(Lexp(-Lnu)}.

\begin{proof}[Proof of Lemma~\ref{lemme: intégrale d'un polynome couplé est un O(Lexp(-Lnu)}]
Writing $P=\sum_{|\alpha|\leq m} p_\alpha X^\alpha\in\mathcal{P}[X_1,...,X_M]$, we have for every $\alpha=(\alpha_1,...,\alpha_M)$ at least two indexes $i\neq j$ such that $\alpha_{i},\alpha_{j}\geq 1$. Since the remaining functions  $(\tau_{a_k}f_k)^{\alpha_k}$ ($k\notin\{i,j \}$) are bounded by a constant $K>0$, we obtain

\begin{align*}
    \normLp{P\big(\tau_{a_1}f_1,...,\tau_{a_M}f_M\big)}&\leq \sum_{|\alpha|\leq m} |p_\alpha|\normLp{(\tau_{a_1}f_1)^{\alpha_1}...(\tau_{a_M}f_M)^{\alpha_M}}\\
    &\leq K\normLp{(\tau_{a_i}f_i)^{\alpha_{i}}(\tau_{a_j}f_j)^{\alpha_{j}}}\\
    & = \mathcal{O}\big( \normLp{e^{-\alpha_i b_i |x-a_i|} e^{-\alpha_j b_j |x-a_j|}}\big).
\end{align*}

Then, since $\alpha_i,\alpha_j\geq 1$ and $y^\pm\leq |y|$ for any $y\in\R$, according to Lemma~\ref{lemme: controle de norme L^p de deux exponentielles}, we have

\begin{equation}
    \normLp{ P\big(\tau_{a_1}(f_1),...,\tau_{a_M}(f_M)\big)}=\grandOde{ \Big(\frac{2}{p\min_{k} (b_k)}+L\Big)^{\frac{1}{p}}e^{-\min_k (b_k)L}}.
\end{equation}
\end{proof}

We now give some asymptotic developments of nonlinear quantities for the chain of solitons, in terms of polynomials of $\mathcal{P}_M$.
\begin{lem}\label{lemme: F,f, f'(1-eta_c)}
\begin{equation*}\label{f(1-eta_c)=sum f(1-eta_ck)+O(S+B+C)}
    F(1-\eta_{\gc,\ga})=\sum_{k=1}^N F(1-\eta_{c_k,a_k})+\mathcal{O}\big(S^{2,N}(\eta_{\gc,\ga})\big)+\mathcal{O}\big(S^{3,N}(\eta_{\gc,\ga})\big) +\mathcal{O}\big(B^N(\eta_{\gc,\ga})\big)+\mathcal{O}\big(C^N(\eta_{\gc,\ga})\big) 
\end{equation*}

\begin{equation*}\label{f(1-eta_c)=sum f(1-eta_ck)+O(S+B)}
    f(1-\eta_{\gc,\ga})=\sum_{k=1}^N f(1-\eta_{c_k,a_k})+\mathcal{O}\big((S^{2,N}+B^N)(\eta_{\gc,\ga})\big),
\end{equation*}

and for any $k\in\{1,...,N\}$

\begin{equation*}\label{f'(1-eta_c)=f'(1-eta_ck)+O(sum eta_k sans j)}
    f'(1-\eta_{\gc,\ga})=f'(1-\eta_{c_k,a_k})+\mathcal{O}\Big(\sum_{\substack{j=1\\j\neq k}}^N \eta_{c_j,a_j}\Big). 
\end{equation*}

\begin{proof}

We write for any $x\in\R$

\begin{equation}
F(1+x)=-\dfrac{x^2}{2}f'(1)-\dfrac{x^3}{2}\int_0^1 (1-t)^2 f''(1+tx)dt,
\end{equation}

hence, by the multinomial formulae,

\begin{align*}
    F(1-\eta_{\gc,\ga})-\sum_{k=1}^N F(1-&\eta_{ c_k,a_k})= -\dfrac{f'(1)}{2}\Big(\eta_{\gc,\ga}^2-\sum_{k=1}^N\eta_{c_k,a_k}^2 \Big)\\
    &-\Big(\eta_{\gc,\ga}^3-\sum_{k=1}^N\eta_{c_k,a_k}^3 \Big)\int_0^1\dfrac{(1-t)^2}{2}f''(1-\eta_{\gc,\ga}t)dt\\
    &- \sum_{k=1}^N \eta_{c_k,a_k}^3\int_0^1 \dfrac{(1-t)^2}{2}  \big(f''(1-\eta_{\gc,\ga}t)-f''(1-\eta_{c_k,a_k}t)\big)dt\\
    &=-f'(1)S^{2,N}(\eta_{\gc,\ga})-(3B^N+6S^{3,N})(\eta_{\gc,\ga})\int_0^1\dfrac{(1-t)^2}{2}f''(1-\eta_{\gc,\ga}t)dt\\
    & - \sum_{k=1}^N \eta_{c_k,a_k}^3\int_0^1 \dfrac{(1-t)^2}{2}  \big(f''(1-\eta_{\gc,\ga}t)-f''(1-\eta_{c_k,a_k}t)\big)dt.\\
\end{align*}

On the other hand, we have for $x\in\R$

\begin{equation}
f(1+x)=xf'(1)+x^2\int_0^1 (1-t)f''(1+tx)dt,
\end{equation}
hence

\begin{align*}
    f(1-\eta_{\gc,\ga})-\sum_{k=1}^N f(1-\eta_{c_k,a_k})
    &=\Big(\eta_{\gc,\ga}^2-\sum_{k=1}^N\eta_{c_k,a_k}^2 \Big)\int_0^1 (1-t) f''(1-t\eta_{\gc,\ga})dt \\
    &\ + \sum_{k=1}^N \eta_{c_k,a_k}^2 \int_0^1 (1-t)\big(f''(1-t\eta_{\gc,\ga})-f''(1-t\eta_{c_k,a_k})\big)dt.
\end{align*}

By Lemma~\ref{lemme: exponential decay}, for any $k\in\{ 1,...,N\}$, $\eta_{c_k,a_k}$ is bounded uniformly with respect to $\gc,\ga$ and $x$. Then by continuity of $f''$ and $f'''$, there exists $M$ independent of $\gc,\ga$ and $x,t$ such that for $i=1,2$, $\big| \int_0^1 \frac{(1-t)^i}{i} f''(1-t\eta_{\gc,\ga})dt\big|\leq M$ and

$$\big|f''(1-t\eta_{\gc,\ga})-f''(1-t\eta_{c_k,a_k})\big|\leq M t |\eta_{\gc,\ga}-\eta_{c_k,a_k}|\leq Mt\Big|\sum_{\substack{j=1 \\j\neq k}}^N \eta_{c_j,a_j}\Big|,$$

then we obtain

\begin{equation}
    F(1-\eta_{\gc,\ga})-\sum_{k=1}^N F(1-\eta_{c_k,a_k})= \mathcal{O}\big(S^{2,N}(\eta_{\gc,\ga})\big)+\mathcal{O}\big(S^{3,N}(\eta_{\gc,\ga})\big) +\mathcal{O}\big(B^N(\eta_{\gc,\ga})\big)+\mathcal{O}\big(C^N(\eta_{\gc,\ga})\big),
\end{equation}

whereas 
\begin{equation}
    f(1-\eta_{\gc,\ga})-\sum_{k=1}^N f(1-\eta_{c_k,a_k})= \mathcal{O}\big(S^{2,N}(\eta_{\gc,\ga})\big) +\mathcal{O}\big(B^N(\eta_{\gc,\ga})\big).
\end{equation}

Dealing the same way than previously with a Taylor expansion of order $1$, we also get the last expression in Lemma~\ref{f'(1-eta_c)=f'(1-eta_ck)+O(sum eta_k sans j)}.
\end{proof}

\end{lem}

\bibliographystyle{plain}
\bibliography{Biblio}

\end{document}